\documentclass[aop]{imsart}

\RequirePackage{amsthm,amsmath,amsfonts,amssymb}
\RequirePackage{mathtools}
\RequirePackage{mathrsfs}
\RequirePackage{xcolor}
\RequirePackage{tikz}
\usetikzlibrary{arrows.meta}
\RequirePackage[authoryear]{natbib}
\RequirePackage[colorlinks,citecolor=blue,urlcolor=blue,linkcolor=blue]{hyperref}

\startlocaldefs

\theoremstyle{plain}
\newtheorem{theorem}{Theorem}[section]
\newtheorem{proposition}[theorem]{Proposition}
\newtheorem{lemma}[theorem]{Lemma}

\newtheorem{definition}[theorem]{Definition}
\theoremstyle{definition}

\numberwithin{equation}{section}

\newcommand{\R}{\mathbb{R}}
\newcommand{\E}{\mathbf{E}}
\newcommand{\PP}{\mathbf{P}}
\newcommand{\1}{\mathbf{1}}
\DeclareMathOperator{\Var}{Var}
\DeclareMathOperator{\Cov}{Cov}
\DeclareMathOperator{\tr}{tr}
\newcommand{\ip}[2]{\langle #1,#2\rangle}
\newcommand{\Ip}[2]{\Big\langle #1,#2\Big\rangle}
\newcommand{\norm}[1]{\lVert #1\rVert}
\newcommand{\abs}[1]{\lvert #1\rvert}
\newcommand{\LapS}{\Delta_\Sigma}
\newcommand{\eps}{\varepsilon}
\newcommand{\J}{\mathsf{J}}
\newcommand{\Lop}{\mathsf{L}}
\newcommand{\fa}{\mathfrak{a}}
\newcommand{\Fc}{\mathcal{F}}
\endlocaldefs

\begin{document}

\begin{frontmatter}
\title{High-Dimensional Gaussian and Bootstrap Approximation:\\ Near $1/\sqrt{n}$ Rates with General Covariance}
\runtitle{Near $1/\sqrt{n}$ rates with general covariance}

\begin{aug}
\author[A]{\fnms{Miles E.}~\snm{Lopes}\ead[label=e1]{melopes@ucdavis.edu}}
\author[A]{\fnms{Ricardo}~\snm{Masini}\ead[label=e2]{rmasini@ucdavis.edu}}
\address[A]{Department of Statistics and Data Science, University of California, Davis\printead[presep={,\ }]{e1,e2}}
\end{aug}

\begin{abstract}
We consider Berry-Esseen bounds over rectangles for sums of $n$~random vectors in $\R^p$ that are i.i.d.~A longstanding open question in high-dimensional probability and statistics is whether bounds of order $1/\sqrt n$ are attainable up to factors of $\log(pn)$ when the population covariance matrix is \emph{unrestricted}. We answer this question affirmatively, provided that the random vectors have sub-Gaussian entries with non-zero variance. In addition, we provide analogous results for  bootstrap approximation and Gaussian comparison, which make use of a special discrepancy measure between covariance matrices that enables improvement with respect to prior work. Taking an approach based on interpolation, the novel aspects of the proofs deal with bounding the interpolation derivatives in a way that is not sensitive to the covariance structure. 
\end{abstract}



\end{frontmatter}

\section{Introduction}\label{sec:intro}
A central topic in high-dimensional probability and statistics is to develop upper bounds on the distance 
\begin{equation}\label{eqn:firstdistance}
\mathsf{D}_n=\sup_{A\in\mathscr{R}}\big|\PP(S_n\in A)-\PP(G\in A)\big|,
\end{equation}
where $S_n=\frac{1}{\sqrt n}\sum_{i=1}^n X_i$ is a normalized sum of centered i.i.d.~random vectors in $\R^p$, $G$ is a centered Gaussian vector with covariance matrix $\E(X_1X_1^{\top})$, and $\mathscr{R}$ is the class of all rectangles in $\R^p$. This measure of distance has attracted major interest because it strikes a special balance between theory and application. On one hand, it is known that accurate Gaussian approximation with respect to $\mathsf{D}_n$ is possible when $p\gg n$, while on the other hand, a wide range of high-dimensional inference procedures can be analyzed in terms of probabilities of the form $\PP(S_n\in A)$. We refer to~\citep{BCCHK18,CCKK23rev,CL25} for overviews of these developments.\\

\noindent\textbf{Prior work.} Much of the contemporary interest in Berry-Esseen type bounds for $\mathsf{D}_n$ was initiated in~\cite{CCK13}, where it was shown that if $\mathscr{R}$ is replaced with a certain subclass of rectangles, and if the entries of $X_1$ have sufficiently light tails, then the corresponding distance is at most of order $n^{-1/8}\log^{7/8}(pn)$. This bound was pathbreaking in two key respects, because it featured mild logarithmic dependence on the dimension $p$, while also imposing no restriction on the off-diagonal entries of $\E(X_1X_1^{\top})$. Moreover, the paper included similar results for Gaussian comparison and bootstrap approximation, which are topics that will be discussed in Section~\ref{sec:main}.

In the subsequent years, there was a cascade of theoretical advances seeking to extend and sharpen the results in~\cite{CCK13}. However, before summarizing this line of work, a  basic point to keep in mind is that for many years, there was not a clear consensus about whether near $n^{-1/2}$ rates of Gaussian approximation over rectangles are even \emph{possible} in a general sense when $p\gg n$. Indeed, at different times, conjectures were made that $n^{-1/6}$ or $n^{-1/4}$ rates of approximation may be optimal up to logarithmic factors when $\Sigma$ is unrestricted. Also, many contributions were made in improving the initial rates from near $n^{-1/8}$ to near $n^{-1/4}$~\citep{CCK17,DZ20,KMB21,CCKK22}, as well as in identifying conditions for ensuring asymptotically consistent Gaussian approximation when $p\gg n$~\citep{DL21,KMB21,KP24,DL25}.

Broadly speaking, most of the progress that has been made in bounding $\mathsf{D}_n$ at faster than near $n^{-1/4}$ rates has been based on additional structural assumptions. For example, the paper~\citep{FK21} obtained near $n^{-1/3}$ rates when the smallest eigenvalue $\lambda_{\min}$ of the correlation matrix of $X_1$ satisfies $\lambda_{\min}\gtrsim 1$, and showed that this can be improved to near $n^{-1/2}$ when $X_1$ also has a log-concave density. In turn, the achievability of $n^{-1/2}$ rates without a log-concave density when $\lambda_{\min}\gtrsim 1$ was shown in~\citep{Lop22,KR20,CCK23,AL26}. 
Furthermore, the paper~\citep{FKLZ23} obtained a similar result when a non-degeneracy assumption is satisfied by all of the $3\times3$ correlation submatrices, and the paper~\citep{Koi24} showed that Edgeworth-based approximations can yield near $n^{-1}$ rates when $\lambda_{\min}\gtrsim 1$ and $X_1$ has a Stein kernel. Meanwhile, it was shown in~\cite{FK24} that a log-concave density for $X_1$ alone is sufficient to ensure near $n^{-1/2}$ rates for general covariance matrices.  In a different direction, the papers~\citep{LLM20,LY22} allow for certain degenerate covariance matrices, but instead work with a subclass of rectangles and assume the entries of $X_1$ satisfy a weak variance decay condition. Lastly, and most recently, the paper~\citep{GZ26} made progress \emph{without} extra structure, showing that near $n^{-1/3}$ rates are achievable for general covariance matrices. \\

\noindent\textbf{Summary of results.} Despite all of the advances just described, the question of whether $\mathsf{D}_n$ can be bounded at near $n^{-1/2}$ rates for general covariance matrices has remained open, and the purpose of this work is to answer it affirmatively. More specifically, our first main result in Theorem~\ref{thm:main} shows that if the entries of $X_1$ are normalized to have unit variance and if the sub-Gaussian condition $\|X_{1j}\|_{\psi_2}\leq B$ holds for all $j\leq p$, then $\mathsf{D}_n$ is at most of order $Bn^{-1/2}\log^3(pn)\log(n)$.

As an extension of our Gaussian approximation result, we consider the related problem of comparing $p$-dimensional Gaussian vectors $\tilde  G\sim\mathcal{N}(0,\tilde \Sigma)$ and $ G\sim\mathcal{N}(0,\Sigma)$ over rectangles. Here, we show that the comparison can be controlled in terms of a special discrepancy measure  $d(\tilde  \Sigma,\Sigma)$, which is defined in~\eqref{eq:bt-disc} and allows for $\tilde\Sigma$ and $\Sigma$ to be degenerate.
In this notation, our Gaussian comparison result in Theorem~\ref{thm:bt} shows that if the entries of $G$ are normalized to have unit variance, and if $d(\tilde\Sigma,\Sigma)\leq \delta$, then the distance $\sup_{A\in\mathscr{R}}|\PP(\tilde  G\in A)-\PP( G\in A)|$ is at most of order $\delta \log(p)\log(1+1/\delta)$.
 To our knowledge, the discrepancy measure $d(\tilde\Sigma,\Sigma)$ is novel in this context, and it is an important technical ingredient that enables Theorem~\ref{thm:bt} to improve upon other recent Gaussian comparison results in high dimensions. (See the remarks after that result for further discussion.) Also, in handling degenerate covariance matrices, this result may be of particular interest in connection with discretized versions of continuous Gaussian processes.

Our last main result deals with approximations involving the Gaussian multiplier bootstrap procedure~\citep{CCK13}, which generates bootstrap samples of $S_n$ defined by $S_n^*=\frac{1}{\sqrt n}\sum_{i=1}^n \xi_i^* X_i$, where $\xi_1^*,\dots,\xi_n^*$ are i.i.d.~$\mathcal{N}(0,1)$ random variables that are independent of $X_1,\dots,X_n$. If the entries of $X_1$ are normalized to have unit variance and satisfy a suitable sub-Gaussian type condition, then  Theorem~\ref{thm:boot} shows that with high probability, the random distance $\sup_{A\in\mathscr R}\big|\PP^*(S_n^*\in A)-\PP(S_n\in A)\big|$ satisfies a bound that is similar to the aforementioned one for $\mathsf{D}_n$, where $\PP^*$ denotes probability conditional on $X_1,\dots,X_n$. Likewise, Theorem~\ref{thm:boot} improves upon previous bootstrap approximation results in ways that are similar to the improvements of Theorem~\ref{thm:main} in the Gaussian approximation context. (The bulk of comparable results on high-dimensional bootstrap approximation are analogous to the ones for Gaussian approximation, and are included in the works cited above.)\\

\noindent\textbf{Proof techniques.} The high-level framework of the proof is based on interpolation and induction, which have been foundational in previous work. If $S_n$ and $G$ are as in eqn.~\eqref{eqn:firstdistance}, and $N\sim \mathcal{N}(0,I_p)$ is independent of them, then our interpolation is built from the ``smoothed orthant indicator'' $g_{u}(y)=\PP\big(y+\sqrt{u}\,G+\sigma N\le\beta\big)$, where $y\in\R^p$ is viewed as the argument, $\beta\in\R^p$ is the corner point of the orthant, $\sigma>0$ is a small smoothing parameter, and $u=u(t)=1-t+\sigma^2$ for $t\in[0,1]$. Because $\sigma$ will be chosen to be negligible, it follows that the function $\Psi(t)=\E(g_{u}(\sqrt t S_n))$ satisfies $\Psi(1)\approx \PP(S_n\leq \beta)$ and $\Psi(0)\approx \PP(G\leq \beta)$. Thus, the primary substance of the proof consists in controlling the difference $\PP(S_n\leq \beta)-\PP(G\leq \beta)\approx \int_0^1 \Psi'(t)dt$ through the derivatives of $g_{u}(y)$.

The essential difficulty in quantifying these derivatives is due to the unrestricted correlation structure of $G$. To see this, note that $g_{u}(y)$ can be rewritten as the  probability that a standard Gaussian vector hits a generic polyhedron with $p$ facets. In other words, the complexity of $\E(X_1X_1^{\top})$ translates into geometric complexity, which then appears in the derivatives of $g_{u}(y)$. Also, it should be emphasized that in the case when $\lambda_{\min}\gtrsim 1$, much of this complexity is concealed. Roughly speaking, this is because the function $g_{u}(y)$ can be effectively replaced with the distribution function of $\mathcal{N}(y,\sigma^2 I_p)$, which has special structure in its higher-order derivatives that can be leveraged in a non-trivial way when $p$ is large~\citep{Ben90,AHT98,FK21,Lop22,CCK23}.

In light of the complexity of high-dimensional polyhedra, much prior work has pursued techniques that circumvent it. Two recent works that have stood out by dealing with polyhedra in a relatively direct way are~\citep{FKLZ23,GZ26}, and our work is similar in this regard. Accordingly, some of the key objects in our analysis are what we call ``facet functions'' $\psi_1(y),\dots,\psi_p(y)$ (see eqns.~\eqref{eqn:nujdef} and~\eqref{eqn:psidef}), which are proportional to the standard Gaussian surface measures of the $p$ facets of a polyhedron $P(y)$ parameterized by $y\in\R^p$. To a large extent, the problem of controlling $\Psi'(t)$ when the matrix $\E(X_1X_1^{\top})$ is degenerate amounts to controlling the first several derivatives of $\psi_1,\dots,\psi_p$. This connection is clarified by a quantity we call the ``facet sensitivity'' $\mathcal{F}(y)$ of $P(y)$, which is a global measure of the sensitivity of the functions $\log \psi_1,\dots,\log \psi_p$ to translations of $P(y)$ (see eqns.~\eqref{eqn:facetsensitivity} and~\eqref{eq:Py} and Proposition~\ref{prop:persummand}). On this point, it should be noted that the upper bounds on $\mathcal{F}(y)$ that we present in Proposition~\ref{thm:FNB} make use of Lemmas 26 and 30 in~\cite{GZ26}. However, our introduction of $\mathcal{F}(y)$ as a way to unify the quantities needed for upper bounding $|\Psi'(t)|$ is a novel component of our work, and our interpolation is quite different from the two separate interpolations used in~\cite{GZ26}.

Another novel part of the proof deals with bounding the curvature of the functions $\tau\mapsto -\log\psi_j(x+\tau w)$ for $\tau\in \R$, where $x,w\in\R^p$ are fixed. Because this curvature is difficult to control in a direct manner, we  reformulate it in terms of the negative second derivative $-A''(\tau)$ of the log-partition function $A(\tau)$ of a certain Gibbs measure (see Proposition~\ref{lem:curv} and Lemma~\ref{lem:sections}). 
The main benefit of this reformulation is that it becomes possible to express $-A''(\tau)$ in terms of a variance computed with respect to the Gibbs measure---and there are a variety of tools for analyzing such variances. In particular, eqn.~\eqref{eq:bartlett} and Lemma~\ref{lem:HSvar} lead to a variational representation of $-A''(\tau)$ in eqn.~\eqref{eqn:Aformula} that ultimately yields the desired curvature bound. This representation could be considered a variant of the classical Helffer-Sj\"ostrand representation of covariances computed with respect to a Gibbs measure~\citep[p.352]{HS94}, and
is derived in the proof of Lemma~\ref{lem:HSvar} using a technique from~\citet[Section~2]{BC13}.\\

\noindent\textbf{Outline.} Section~\ref{sec:main} states the main results. Section~\ref{sec:reductions} gives a high-level proof of the main result on Gaussian approximation (Theorem~\ref{thm:main}), based on lower-level results whose proofs occupy Sections~\ref{sec:onesummand}--\ref{sec:tools}. The main results on Gaussian comparison (Theorem~\ref{thm:bt}) and bootstrap approximation~(Theorem~\ref{thm:boot}) are proved in Sections~\ref{sec:bootstrap} and~\ref{sec:bootproof} respectively.  Section~\ref{sec:AI} is a statement on AI use.\\

\noindent\textbf{Notation and conventions.}
If $a$ and $b$ are real numbers, we write $a\lesssim b$ if $a\le Cb$ for a positive absolute constant $C$ whose value may change from line to line. Also, we write $a=\mathcal{O}(b)$ if $|a|\lesssim|b|$. 
The all-ones vector is $\1$, while the expressions $\1\{\cdot\}$ and $\1_A$ denote indicator functions.  If $x\in\R^p$, we write $\|x\|_2$ for the Euclidean norm and let $\|x\|_{\infty}=\max_{j\leq p}|x_j|$.   If $A$ is a matrix or tensor, then $\|A\|_{\infty}$ refers to the maximum absolute value of the entries, while we write $\|A\|_F$ for the Frobenius norm and $\|A\|_{\textup{op}}$ for the operator (spectral) norm. For a real-valued function $F$, let $\norm F_{\sup}=\sup_x\abs{F(x)}$. When $F$ takes values that are vectors, let $\|F\|_{\sup}=\sup_x \|F(x)\|_2$, and when $F$ takes values that are matrices or tensors, let $\|F\|_{\sup}=\sup_x\norm{F(x)}_{\textup{op}}$. If $x$ and $y$ are vectors in $\R^p$, then $x\leq y$ means $x_j\leq y_j$ for all $j\leq p$. The Loewner partial order on symmetric matrices is denoted by $\preceq$. The
standard Gaussian measure and density on $\R^p$ are denoted by $\gamma_p$ and $\phi_p$, while $\Phi$ and $\phi$ denote the distribution and density functions for $\mathcal{N}(0,1)$. For a linear subspace $E\subset \R^p$,
the standard Gaussian measure on $E$ is denoted as $\gamma_E$. We also write $\mathcal H^d$ for the $d$-dimensional Hausdorff measure.
The derivative in the direction of a vector $w$ is denoted as $\partial_w$, and we write $\partial_j=\partial_{e_j}$, where $e_j$ is the $j$th standard basis vector. If $\zeta$ is a real-valued random variable, its $\psi_2$-Orlicz norm  is given by $\norm{\zeta}_{\psi_2}=\inf\{t>0:\E e^{\zeta^2/t^2}\le2\}$. Finally, the collection of $p$-dimensional rectangles $\mathscr{R}$ consists of all sets of the form $A=\prod_{j=1}^p \mathcal{I}_j$, where $\mathcal{I}_1,\dots,\mathcal{I}_p$ are arbitrary intervals in $\R$ that may be bounded or unbounded, and contain zero, one, or two endpoints.

\section{Main results}\label{sec:main}

This section states our three main results on Gaussian approximation, Gaussian comparison, and bootstrap approximation.

\subsection{Gaussian approximation}
Our first result concerns the normalized sum of centered  i.i.d.\
random vectors whose coordinates have unit variance and sub-Gaussian tails. No assumption is made on the correlations between
the coordinates.

\begin{theorem}\label{thm:main}
There is an absolute constant $C$ such that the following holds. Let $n,p\ge2$, and let $X_1,\dots,X_n$ be i.i.d.\ centered
random vectors in $\R^p$ with $\Var(X_{1j})=1$ and $\norm{X_{1j}}_{\psi_2}\le B$ for all $j\le p$. Also, let $G$ be a centered Gaussian vector with covariance matrix $\E(X_1X_1^{\top})$. Then,
\[
\sup_{A\in\mathscr{R}}\Big|\PP\Big(\frac1{\sqrt n}\sum_{i=1}^nX_i\in A\Big)-\PP(G\in A)\Big|\;\le\;\frac{C\,B\,\log^{3}(pn)\log(n)}{\sqrt n}.
\]
\end{theorem}
This result is proven in Section~\ref{sec:reductions}, and it only depends on the distribution of
$X_1$ through the parameter $B$. In particular, the covariance matrix $\Sigma$ is allowed to be degenerate, and no
lower bound on its eigenvalues is required. If the variances $\Var(X_{1j})$ are positive but otherwise arbitrary, then the theorem can be applied by replacing $X_{1j}$ with $X_{1j}/\sqrt{\Var(X_{1j})}$ for $j\leq p$ because rescaling the coordinates maps rectangles to rectangles.

\subsection{Gaussian comparison}
Our next result compares $p$-dimensional Gaussian vectors $\tilde G\sim \mathcal{N}(0,\tilde\Sigma)$ and $G\sim \mathcal{N}(0,\Sigma)$ over rectangles. The comparison is made in terms of the following measure of
discrepancy between their covariance matrices

\begin{equation}\label{eq:bt-disc}
d(\tilde\Sigma,\Sigma)=\inf\Big\{t>0: |\ip{v}{(\tilde\Sigma-\Sigma)v}|\leq t\, \ip{v}{\Sigma v} \text{ for all } v\in\mathcal{V}\Big\},
\end{equation}
where $\mathcal V=\{e_1,\dots,e_p\}\cup\{e_j\pm e_l :j<l\}$.
 From this definition, it is straightforward to check that $d(\tilde\Sigma,\Sigma)=\infty$ can only arise if $G$ has an entry with 0 variance, or if it has a pair of (redundant) entries with correlation equal to $\pm 1$. If these exceptional cases do not occur, then the simpler formula  
 $$d(\tilde\Sigma,\Sigma)=\max_{v\in\mathcal{V}}\frac{|\langle v,(\tilde\Sigma-\Sigma)v\rangle|}{\langle v,\Sigma v\rangle }$$
 can be used instead. Importantly,
the following result only relies on the finiteness of $d(\tilde\Sigma,\Sigma)$ once the coordinates of $G$ have been normalized. Its proof is given in Section~\ref{sec:bootstrap}.

\begin{theorem}\label{thm:bt}
There is an absolute constant $C$ such that the following holds. Let $p\ge2$, let $\Sigma$ and $\tilde\Sigma$ be positive semidefinite
$p\times p$ matrices such that $\Sigma_{jj}=1$ for $j\leq p$. Also, let $ G\sim\mathcal N(0,\Sigma)$ and $\tilde G\sim\mathcal N(0,\tilde\Sigma)$.
If $d(\tilde\Sigma,\Sigma)\le\delta$ for some number $\delta>0$, then
\begin{equation*}
\sup_{A\in\mathscr R}\big|\PP(\tilde G\in A)-\PP(G\in A)\big|\le C\,\delta\,\log(p)\log(1+1/\delta).
\end{equation*}
\end{theorem}
If the diagonal entries of $\Sigma$ are positive but different from 1, then Theorem~\ref{thm:bt} can still be used by conjugating $\tilde\Sigma$ and $\Sigma$ with $\textup{Diag}(\Sigma)^{-1/2}$.
Apart from the generality of the covariance matrices, another valuable feature of Theorem~\ref{thm:bt} is that the bound has a \emph{near-linear} dependence on the parameter $\delta$. Previously, Gaussian comparison results over rectangles have been established with near-linear dependence on the distinct discrepancy measure $\|\tilde\Sigma-\Sigma\|_{\infty}$~\citep{FK21,Lop22}.
However, such bounds have the limitation of depending on  the reciprocal of the smallest eigenvalue $\lambda_{\min}$ of the correlation matrix of $G$. To show that $d(\tilde\Sigma,\Sigma)$  compares favorably, it can be checked that if $\Sigma_{jj}=1$ for $j\leq p$ and if $\lambda_{\min}>0$, then the following bound holds for any choice of $\tilde\Sigma$,
$$d(\tilde\Sigma,\Sigma)\leq \frac{2\|\tilde\Sigma-\Sigma\|_{\infty}}{\lambda_{\min}}.$$ 
Earlier still, Gaussian comparison results were established in~\citep{CCK15,CCKK22} without placing constraints on the correlation matrices of $\tilde  G$ and $G$, but with rates of order $\|\tilde\Sigma-\Sigma\|_{\infty}^{1/3}$ or $\|\tilde\Sigma-\Sigma\|_{\infty}^{1/2}$, up to logarithmic factors.

One further benefit of the discrepancy measure $d(\tilde\Sigma,\Sigma)$ is that it is well behaved when $\tilde\Sigma$ is replaced by a sample covariance matrix $\hat\Sigma$ of centered observations $X_1,\dots,X_n$ whose population covariance matrix is $\Sigma$. That is, even if $\ip{v}{\Sigma v}=0$ for some $v\in\mathcal{V}$, then it follows that $\ip{v}{X_i}=0$ almost surely for all $i$, and thus $d(\hat\Sigma,\Sigma)$ remains finite almost surely. This special property of $d$ will make Theorem~\ref{thm:bt} useful toward obtaining a near $1/\sqrt n$ rate of bootstrap approximation below.

\subsection{Bootstrap approximation}
The Gaussian comparison result in Theorem~\ref{thm:bt} serves as a bridge to extend our Gaussian approximation result in Theorem~\ref{thm:main} to a bootstrap approximation result. Let
$\xi_1^*,\dots,\xi_n^*$ be i.i.d.\ $\mathcal N(0,1)$ random variables that are independent of $X_1,\dots,X_n$, and consider the Gaussian multiplier bootstrap
version of $S_n$ given by
\[
S_n^*=\frac1{\sqrt n}\sum_{i=1}^n\xi_i^*X_i,
\]
which was popularized by~\cite{CCK13}. Conditionally on $X_1,\dots,X_n$, the vector $S_n^*$ is centered Gaussian with covariance matrix $\hat\Sigma=\frac1n\sum_{i=1}^nX_iX_i^\top$. Likewise, we write $\PP^*$ for
probability conditional on $X_1,\dots,X_n$.

\begin{theorem}\label{thm:boot}
There is an absolute constant $C$ such that the following holds. Let $n,p\ge2$, and let $X_1,\dots,X_n$ be i.i.d.\ centered random vectors
in $\R^p$ with covariance matrix $\Sigma=\E(X_1X_1^\top)$ such that $\Sigma_{jj}=1$ for all $j\le p$, and
$\norm{\ip{X_1}{v}}_{\psi_2}\le B\sqrt{\ip{v}{\Sigma v}}$ for all $v\in\mathcal V$. Then, for every $\alpha\in(0,1)$, the following event occurs
with probability at least $1-\alpha$,
\begin{equation*}
\sup_{A\in\mathscr R}\big|\PP^*(S_n^*\in A)-\PP(S_n\in A)\big|\le \frac{C B\log(n)}{\sqrt n}\Big(B\sqrt{\log(p/\alpha)}\,\log(p)+\log^3(pn)\Big).
\end{equation*}
\end{theorem}
The proof of Theorem~\ref{thm:boot} is given in Section~\ref{sec:bootproof}. Analogously to Theorem~\ref{thm:main}, the main novelty of this result is that it provides a near $n^{-1/2}$ rate of approximation even when $\Sigma$ may be far from invertible.
 The tail conditions on $X_1$ in this result are only modestly stronger than those in Theorem~\ref{thm:main}, to the extent that the set $\mathcal{V}$ expands $\{e_1,\dots,e_p\}$ by including the vectors in $\{e_j\pm e_l: j<l\}$. For context, it is common in the literature on high-dimensional probability and statistics to study settings where $\mathcal{V}$ is replaced with $\R^p$ i.e.~where there is some $B>0$ such that $\norm{\ip{X_1}{v}}_{\psi_2}\le B\sqrt{\ip{v}{\Sigma v}}$ holds for all $v\in\R^p$, \citep[e.g.][]{KL17,Ver18,Zhi24}.

\section{Proof of Theorem~\ref{thm:main}}\label{sec:reductions}

Standard arguments show that Theorem~\ref{thm:main} can be reduced to the case when  $X_1,\dots,X_n$ have bounded entries and the rectangles $A$ are  orthants, as recorded in Lemma~\ref{lem:reduction} of the background results section. Hence, Theorem~\ref{thm:main} is implied by the following result, whose proof will be covered in the rest of this section.

\begin{theorem}\label{thm:bounded}
There is an absolute constant $C$ such that the following holds. Let \smash{$p,n\ge2$,} let $\bar B>0$, and let $X_1,\dots,X_n$ be i.i.d.\ centered random vectors in $\R^p$ with
$\abs{X_{1j}}\le \bar B$ almost surely and $\Var(X_{1j})=1$ for all $j\leq p$. Also, let $G$ be a centered Gaussian vector with covariance matrix $\E(X_1X_1^\top)$, and for any $k=1,\dots,n$, let
$$ \bar{\mathsf D}_k=\sup_{\beta\in\R^p}\Big|\PP\Big(\frac1{\sqrt k}\sum_{i=1}^kX_i\le\beta\Big)-\PP(G\le\beta)\Big|.$$
Then, for any $k=1,\dots,n$, we have
\begin{equation}\label{eqn:thm:bound}
    \bar{\mathsf D}_k \leq \frac{C \bar B\,(\log(pn))^{5/2}\log(n)}{\sqrt k}.
\end{equation}
\end{theorem}

The proof of Theorem~\ref{thm:bounded} is organized as follows. After some preliminary reductions, the distributions of $\frac1{\sqrt k}\sum_{i=1}^kX_i$ and $G$ are smoothed by a small Gaussian noise, and then
connected by an interpolation (Lemma~\ref{lem:interp}). The derivatives appearing in the interpolation are controlled through a special quantity we call the facet sensitivity, which is introduced and analyzed in  Sections~\ref{subsec:sensitivity}-\ref{subsec:onesummand}. (The proofs of the lower-level results are deferred to Sections~\ref{sec:onesummand}--\ref{sec:tools}.) Finally, the remaining ingredients are assembled in Section~\ref{subsec:induction}, where the proof of Theorem~\ref{thm:bounded} is completed with an induction on $k$.\\

\noindent \textbf{Technical preliminaries and further reductions.} 
Throughout the proof of Theorem~\ref{thm:bounded}, we may restrict our attention to the values of $k=1,\dots,n$ that satisfy
\begin{equation}\label{eq:small}
\frac{\bar B(\log(pn))^{5/2}\log(n)}{\sqrt k}\le\frac14,
\end{equation}
for otherwise the conclusion~\eqref{eqn:thm:bound} is trivially true. Similarly, we may assume that
$pn$ is large enough so that $pn\geq c$ holds when $c>0$ is an absolute constant. To see this, recall that $p,n\geq 2$, and note that the almost-sure bound $|X_{ij}|\leq \bar B$ with $\Var(X_{ij})=1$ implies $\bar B\geq 1$. If the inequality $pn<c$ held, then we would have $k\leq n\leq c$, implying $\big[c^{1/2} (\log 4)^{-5/2}(\log 2)^{-1}\big]\bar B (\log(pn))^{5/2}\log(n)/\sqrt k\geq 1$, and~\eqref{eqn:thm:bound} would again be trivially true. For this reason, various intermediate results that depend on $pn$ exceeding an absolute constant will use phrases  such as ``when $pn$ is sufficiently large'' in order to avoid the clutter of an extra recurring absolute constant.

From now until the induction step at the end of the proof in Section~\ref{subsec:induction}, we fix $k=2,\dots,n$ such that~\eqref{eq:small} holds, and we put
\[
\sigma=\frac{C_0\,\bar B(\log(pn))^{3/2}}{\sqrt k},\qquad S=\frac1{\sqrt k}\sum_{i=1}^kX_i,\qquad S_{\setminus i}=S-\frac{X_i}{\sqrt k},
\]
where $C_0\ge1$ is an absolute constant that is chosen large enough for the few steps below that require it. (We do not put a subscript $k$ on $\sigma$ and $S$ above, because $k$ is
held fixed for almost the entire proof.) Unlike the generic absolute constant $C$, the value
of $C_0$ stays fixed throughout the proof and does not change from line to line. Also, when we use the notation $a\lesssim b$ to express the statement that $a\leq C b$ for some absolute constant $C>0$, the implicit constant $C$ will not depend on $C_0$.

Some frequently used relationships among the parameters are as follows. Since $\bar B\ge1$, we have $\sigma\ge \bar B/\sqrt k\ge n^{-1/2}$, and~\eqref{eq:small} gives
$\bar B\le\sqrt k$ and $\sigma\le C_0(\log(pn))^{-1}(\log(n))^{-1}/4\le1$, as $pn$ is large. Moreover, each coordinate of $X_i/\sqrt k$ is bounded
almost surely in magnitude by
\begin{equation}\label{eq:ratio}
\frac{\bar B}{\sqrt k}=\frac{\sigma}{C_0(\log(pn))^{3/2}},
\end{equation}
which is small compared to the smoothing parameter $\sigma$.

Next, let $\Sigma=\E(X_1X_1^{\top})$ so that $G\sim\mathcal{N}(0,\Sigma)$. Also, let $G'$ be an independent copy of $G$, and let $N\sim\mathcal N(0,I_p)$, with the collection $G,G',N,X_1,\dots,X_n$ being independent. Both $S$ and $G$ will be perturbed by the
noise vector $\sigma G'+\sigma N$. The importance of $\sigma N$ is that its covariance matrix is equal to $\sigma^2 I_p$, which provides a small degree of
smoothing in every direction, and thus ensures a basic level of differentiability even when $\Sigma$ is degenerate. The smoothing effect of $\sigma G'+\sigma N$ is quantified by a standard argument (Lemma~\ref{lem:smoothing}) that gives
\begin{equation}\label{eq:smoothing}
\bar{\mathsf D}_k\le\sup_{\beta\in\R^p}\bigg|\PP\Big(S+\sigma G'+\sigma N\le\beta\Big)-\PP\Big(G+\sigma G'+\sigma N\le\beta\Big)\bigg|+C\sigma\log(pn).
\end{equation}
Thus, it remains to bound the supremum
on the right side of \eqref{eq:smoothing}, and this will be done using the interpolation described below.\\

\noindent\textbf{Interpolation.} 
From now on $\beta\in\R^p$ is fixed. We interpolate between the two smoothed distributions in \eqref{eq:smoothing} by means of
the following functions. For $u>0$ and $y\in\R^p$, let
\[
g_u(y)=\PP\Big(y+\sqrt u\,G+\sigma N\le\beta\Big),
\]
which may be viewed as an orthant indicator function that has been smoothed by convolution with a Gaussian density. Also, for $t\in[0,1]$, define the interpolation function
\[
\Psi(t)=\E g_{u(t)}(\sqrt tS),
\]
where we put $u(t)=1-t+\sigma^2$.
This leads to
\begin{equation}\label{eq:interpFTC}
\PP\Big(S+\sigma G'+\sigma N\le\beta\Big)-\PP\Big(G+\sigma G'+\sigma N\le\beta\Big)=\int_0^1\Psi'(t)\,dt.
\end{equation}
Strictly speaking, the integral is understood as a limit with the lower endpoint tending to 0, which is justified by Proposition~\ref{prop:pertime}. The following lemma provides an explicit expression for the derivative $\Psi'(t)$ in terms of the operator $\LapS f=\tr(\Sigma\nabla^2f)$ acting on $f\in C^2(\R^p)$.

\begin{lemma}\label{lem:interp} If $t\in (0,1]$ and $u=u(t)$, then the derivative $\Psi'(t)$ satisfies
\begin{equation}\label{eqn:interpformula}
\Psi'(t)=\frac1{2\sqrt t}\sum_{i=1}^k\E\Ip{\nabla g_{u}\Big(\sqrt tS_{\setminus i}+\sqrt{\tfrac tk}\,X_i\Big)-\nabla g_{u}(\sqrt tS_{\setminus i})}{\frac{X_i}{\sqrt k}}
-\frac12\E\LapS g_u(\sqrt tS).
\end{equation}
\end{lemma}
\noindent This formula is a basic consequence of Stein's identity for Gaussian vectors and is proved in \smash{Section~\ref{sec:remaining}}. 

\subsection{Polyhedra and the facet sensitivity}\label{subsec:sensitivity}

The terms in the formula~\eqref{eqn:interpformula} will eventually be controlled by a weighted sum of
the Gaussian surface measures of the facets of certain polyhedra, which we now define. \\

\begin{figure}[ht]
\centering
\begin{tikzpicture}[scale=1.35,>=stealth]
%%% (a) 0 \in P %%%
% P = {x : <n_1,x> <= a_1, <n_2,x> <= a_2} is an unbounded wedge (p=2); only a window is shaded
\begin{scope}[shift={(0,0.93)}]
\fill[gray!15] (-2.1,1.3) -- (0.751,1.3) -- (2.368,-1.5) -- (-2.1,-1.5) -- cycle;
\draw[->,gray] (-2.1,0) -- (2.3,0) node[right,black] {$E_1=n_1^\perp$};
\draw[very thick] (-2.1,1.3) -- (0.751,1.3);
\draw[very thick] (0.751,1.3) -- (2.368,-1.5);
\draw[very thick,blue] (-2.1,0) -- (0.751,0);
\draw[dashed,blue] (0.751,0) -- (0.751,1.3);
\fill (0,0) circle (0.03) node[below left] {$0$};
\draw[->,thick] (0,0) -- (0,1.3) node[midway,left] {$a_1n_1$};
\draw[->,thick] (-0.35,1.3) -- (-0.35,1.85) node[left] {$n_1$};
\draw[->,thick] (1.125,0.650) -- (1.601,0.925) node[right] {$n_2$};
\node[above] at (-1.0,1.3) {$F_1$};
\node[right] at (1.95,-0.45) {$F_2$};
\node[below,blue] at (-0.6,0) {$K_1$};
\node at (-0.9,-0.8) {$P$};
\end{scope}
\end{tikzpicture}
\caption{An illustration of the notation describing a polyhedron $P$ in $\R^2$ with two linearly independent normal vectors. Note that $P$ is unbounded, and so only part of it is shaded.}
\label{fig:polyhedra}
\end{figure}
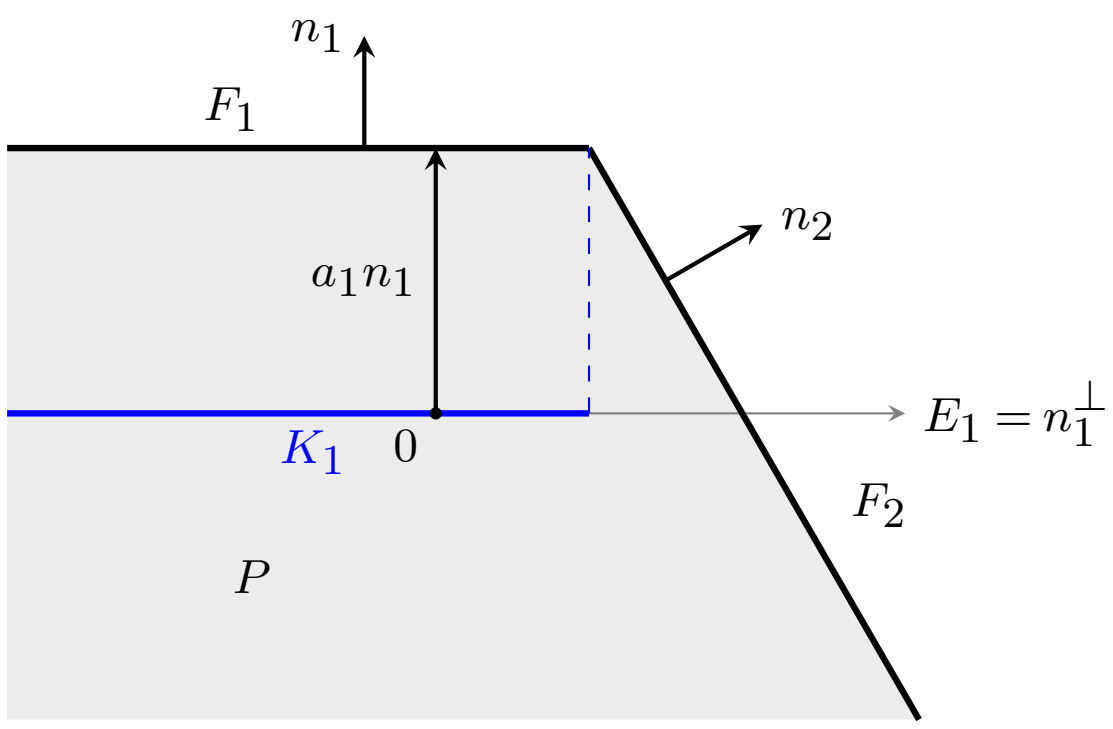

\noindent\textbf{Notation for polyhedra.}
A polyhedron in $\R^p$ with linearly independent normal vectors $n_1,\dots,n_p$ (with unit norm) and offset parameters $a=(a_1,\dots,a_p)$ is denoted as
\begin{equation}\label{eqn:Polygeneral}
P=\{z\in\R^p:\ip{n_j}{z}\le a_j\ \text{for } j=1,\dots,p\}.
\end{equation}
For each $j\leq p$, the $j$th facet of $P$ is
\begin{equation*}F_j=P\cap\{z\in\R^p:\ip{n_j}{z}=a_j\},
\end{equation*}
and its Gaussian surface measure is
\begin{equation}\label{eqn:nujdef}
\nu_j(P)=\int_{F_j}\phi_p\,d\mathcal H^{p-1},
\end{equation}
where $\mathcal H^{p-1}$ denotes the $(p-1)$-dimensional Hausdorff measure. (An illustration of the notation describing polyhedra is given in Figure~\ref{fig:polyhedra} later.)
When $P$ is fixed, we also write $\nu_j=\nu_j(P)$. For $x\in\R^p$, the translate $P-x$ is a polyhedron of the same form, with
the same normal vectors, while the $j$th offset and facet measure are respectively
$a_j-\ip{n_j}{x}$
and $\nu_j(P-x)$. For the sake of understanding various formulas given later, it is worth noting that the function $x\mapsto\nu_j(P-x)$ is positive and $C^\infty$ on $\R^p$ (Lemma~\ref{lem:facetlog}).

The following quantity plays a central role in our work. It quantifies the sensitivity of the facet measures $\nu_j$ to translations of the polyhedron $P$.

\begin{definition}[facet sensitivity]\label{def:sensitivity}
If $P\subset\R^p$ is a polyhedron as in~\eqref{eqn:Polygeneral} with $p$ facets having linearly independent normal vectors, then the facet sensitivity of $P$ is defined as
\begin{equation}\label{eqn:facetsensitivity}
\Fc(P)=\sum_{j=1}^p \nu_j(P) \bigg[\norm{\nabla_x\log\nu_j(P-x)}_2^2+\tr\big(-\!\nabla_x^2\log\nu_j(P-x)\big)\bigg]\bigg|_{x=0}\,.
\end{equation}
\end{definition}
\noindent Importantly, the facet sensitivity satisfies $\mathcal{F}(P)\geq 0$, because the functions $x\mapsto \nu_j(P-x)$ are log-concave (Lemma~\ref{lem:facetlog}). Some interpretation for $\mathcal{F}(P)$ is possible, since its terms are akin to the two classical forms of the trace of the Fisher information. An explicit formula for $\mathcal{F}(P)$ in terms of the characteristics of $P$ will be given in Lemma~\ref{lem:sensitivityformula}.\\

\noindent\textbf{Parameterized polyhedra.} Recall that we let $\Sigma=\E(X_1X_1^{\top})$, and for any $u>0$, let
 $$\Sigma_u=u\Sigma+\sigma^2I_p,$$
 which is a positive definite matrix with diagonal entries $u+\sigma^2$. For any $j\le p$ and $y\in\R^p$, let
\begin{equation*}n_j=\frac{\Sigma_u^{1/2}e_j}{\sqrt{u+\sigma^2}},\qquad a_j(y)=\frac{\beta_j-y_j}{\sqrt{u+\sigma^2}},
\end{equation*}
and denote the polyhedron with offsets $a_1(y),\dots,a_p(y)$ as
\begin{equation}\label{eq:Py}
\begin{split}
    P(y)&=\{z\in\R^p:\ip{n_j}{z}\le a_j(y)\ \text{for } j\le p\}\\
    &=\{z\in\R^p:\Sigma_u^{1/2}z\le\beta-y\}.
    \end{split}
\end{equation}
The Gram matrix of $n_1,\dots,n_p$ is $\Sigma_u/(u+\sigma^2)$, which is invertible with unit diagonal entries. Hence, the vectors $n_1,\dots,n_p$ are linearly independent unit vectors with $\ip{n_j}{n_l}=u\Sigma_{jl}/(u+\sigma^2)$ for $j\ne l$, and so $P(y)$ is an instance of the polyhedron $P$ discussed earlier. Likewise, we write $\Fc(y)=\Fc(P(y))$ for any $y\in\R^p$.

The polyhedron $P(y)$ is directly connected to the smoothed orthant indicator $g_u(y)$, because~\eqref{eq:Py} shows that
\begin{equation*}
g_u(y)=\gamma_p(P(y)),
\end{equation*}
where $\gamma_p$ is the standard Gaussian measure on $\R^p$. Furthermore, the derivatives of $g_u$ are linked to the facet measures $\nu_j$ through the formula
\begin{equation}\label{eqn:psidef}
-\partial_jg_u(y)=\psi_j(y):=\frac{\nu_j(P(y))}{\sqrt{u+\sigma^2}},
\end{equation}
which will be established in Lemma~\ref{lem:gu}(b), and provides some explanation as to why the facet sensitivity $\mathcal{F}(y)$ will allow us to control the interpolation.
We call $\psi_1,\dots,\psi_p$ the \emph{facet functions}, and they will be a recurring part of our analysis. Lastly, note that because $P(y)$ is defined in terms of $\Sigma_u$, the various objects associated with $P(y)$ depend on $u$ in general, but this will usually be suppressed for ease of notation.

\subsection{Bounds involving the facet sensitivity}\label{subsec:onesummand}
The following proposition is key, because it expresses the connection between the derivatives of the smoothed orthant indicator $g_u(y)$ and the facet sensitivity $\mathcal{F}(y)$ in quantitative terms. It is proved in Section~\ref{sec:onesummand}.

\begin{proposition}\label{prop:persummand}
Let $t\in(0,1]$, $u=u(t)$ and $x\in\R^p$. Then,
\[
\bigg|\E\Ip{\nabla g_u\Big(x+\sqrt{\tfrac tk}\,X_1\Big)-\nabla g_u(x)}{\frac {X_1}{\sqrt k}}-\frac{\sqrt t}k\LapS g_u(x)\bigg|
\lesssim\frac{t\bar B\,\Fc(x)}{(ku)^{3/2}}+\frac1{(pn)^4}.
\]
\end{proposition}

The next proposition provides upper bounds on the facet sensitivity. Notably, these bounds do not depend on the
angles between the normal vectors, and they have well controlled dependence on $p$. Part (a) follows from Lemmas~26 and~30 in~\cite{GZ26}. Part (b) is a new result, and part (c) combines parts (a) and (b) in the case of $P=P(y)$.

\begin{proposition}\label{thm:FNB}
If $P$ is a polyhedron in $\R^p$ of the form~\eqref{eqn:Polygeneral} having $p$ facets with linearly independent normal vectors, then the following statements hold.
\begin{itemize}
\item[(a)] The facet sensitivity satisfies $\Fc(P)\lesssim(\log(2p))^{3/2}$.
\item[(b)] If $r>0$ and  $\abs{\min_{j\le p}a_j}\ge r$, then $\Fc(P)\lesssim p^2e^{-r^2/16}$.
\item[(c)] If $P=P(y)$ with $u>0$ and $y\in\R^p$ as in~\eqref{eq:Py}, then
\[
\Fc(y)\lesssim(\log(pn))^{3/2}\,\1\Big\{\Big|\min_{j\le p}(\beta_j-y_j)\Big|<8\sqrt{(u+\sigma^2)\log(pn)}\Big\}+(pn)^{-2}.
\]
\end{itemize}
\end{proposition}
\noindent This result is proved in Section~\ref{sec:FNB}. For the polyhedron $P(y)$, the offset parameters are
$a_j(y)=(\beta_j-y_j)/\sqrt{u+\sigma^2}$, and so the assumption of part (b) holds when $y$ is  ``far'' from the boundary of the orthant
$\{y\le\beta\}$ with respect to the scale $\sqrt{u+\sigma^2}$. 

The facet sensitivity also gives pointwise bounds on the third derivatives $\partial_w\LapS g_u$, as shown in the following lemma that is proved in Section~\ref{sec:tools}.

\begin{lemma}\label{lem:lap3}
For all $u>0$ and $y,w\in\R^p$, we have
\[
\abs{\partial_w\LapS g_u(y)}\le \textstyle\frac{\norm w_\infty}{u^{3/2}}\Fc(y).
\]
\end{lemma}

We can now obtain a pointwise bound for the derivative $\Psi'$ in the interpolation.

\begin{proposition}\label{prop:pertime}
For $t\in(0,1]$ and $u(t)=1-t+\sigma^2$, let
\[
r_t=20\sqrt{u(t)\log(pn)}+\frac{\bar B}{\sqrt k}\qquad\text{and}\qquad\pi_t=\PP\Big(\Big|\min_{j\le p}(\beta-\sqrt tS_{\setminus 1})_j\Big|<r_t\Big).
\]
Then,
\[
\abs{\Psi'(t)}\lesssim\frac{\sqrt t\,\bar B(\log(pn))^{3/2}}{\sqrt k\,u(t)^{3/2}}\big(\pi_t+(pn)^{-2}\big)+\frac{1}{pn\sqrt t}.
\]
\end{proposition}
\begin{proof}
For each $y\in\R^p$, define the remainder function
\[
\mathcal R(y)=\E\Ip{\nabla g_{u(t)}\Big(y+\sqrt{\tfrac tk}\,X_1\Big)-\nabla g_{u(t)}(y)}{\frac {X_1}{\sqrt k}}-\frac{\sqrt t}k\LapS g_{u(t)}(y),
\]
which is bounded in Proposition~\ref{prop:persummand}.  Since $S_{\setminus i}$ is independent of $X_i$, and $X_i$ has the same distribution as $X_1$, we may express the formula in  Lemma~\ref{lem:interp} with this function and then rearrange the terms to obtain
\begin{align}\label{eq:s7-split}
\Psi'(t)&=\frac1{2\sqrt t}\sum_{i=1}^k\E\Big(\mathcal R(\sqrt tS_{\setminus i})+\frac{\sqrt t}k\LapS g_{u(t)}(\sqrt tS_{\setminus i})\Big)-\frac12\E\LapS g_{u(t)}(\sqrt tS)\\[2mm]
&=\frac1{2k}\sum_{i=1}^k\E\Big(\LapS g_{u(t)}(\sqrt tS_{\setminus i})-\LapS g_{u(t)}(\sqrt tS)\Big)+\frac1{2\sqrt t}\sum_{i=1}^k\E\,\mathcal R(\sqrt tS_{\setminus i}).\label{eq:s7-splitsecond}
\end{align}

We will control both of the sums in~\eqref{eq:s7-splitsecond} using the facet sensitivity near the points $\sqrt tS_{\setminus i}$. To this end, put $w_i=\sqrt{t/k}\,X_i$, so that
$\sqrt tS=\sqrt tS_{\setminus i}+w_i$ and $\norm{w_i}_\infty\le\sqrt t\,\bar B/\sqrt k$. As a preparatory step, we first show that
\begin{equation}\label{eq:s7-layer}
\E\sup_{0\le\theta\le1}\Fc(\sqrt tS_{\setminus i}+\theta w_i)\lesssim(\log(pn))^{3/2}\pi_t+(pn)^{-2},\qquad i\le k.
\end{equation}
Since $u(t)+\sigma^2\le2u(t)$, the threshold in Proposition~\ref{thm:FNB}(c) satisfies $8\sqrt{(u(t)+\sigma^2)\log(pn)}\le r_t-\bar B/\sqrt k.$
Due to this inequality and $\norm{\theta w_i}_\infty\le \bar B/\sqrt k$, if we put $y=\sqrt tS_{\setminus i}+\theta w_i$ inside the indicator of Proposition~\ref{thm:FNB}(c), then it is at most $\1\{|\min_{j\leq p}(\beta-\sqrt t  S_{\backslash i})_j|<r_t\}$. Furthermore, the latter indicator has expectation equal to $\pi_t$, because $S_{\setminus i}$ and $S_{\setminus 1}$ have the same distribution, and so \eqref{eq:s7-layer} follows.

For the sum of remainders in \eqref{eq:s7-splitsecond}, we apply Proposition~\ref{prop:persummand} with $x=\sqrt tS_{\setminus i}$, and then \eqref{eq:s7-layer} with $\theta=0$.
In turn, summing over $i\le k$ shows that
\[
\frac1{\sqrt t}\sum_{i=1}^k\abs{\E\,\mathcal R(\sqrt tS_{\setminus i})}\lesssim\frac{\sqrt t\,\bar B}{\sqrt k\,u(t)^{3/2}}\big((\log(pn))^{3/2}\pi_t+(pn)^{-2}\big)+\frac{1}{pn\sqrt t}.
\]

 Turning to the first sum in \eqref{eq:s7-splitsecond}, the mean value theorem and Lemma~\ref{lem:lap3} give
\begin{align*}
\abs{\LapS g_{u(t)}(\sqrt tS_{\setminus i}+w_i)-\LapS g_{u(t)}(\sqrt tS_{\setminus i})}&\le\sup_{0\le\theta\le1}\abs{\partial_{w_i}\LapS g_{u(t)}(\sqrt tS_{\setminus i}+\theta w_i)}\\
&\le\frac{\sqrt t\,\bar B}{\sqrt k\,u(t)^{3/2}}\sup_{0\le\theta\le1}\Fc(\sqrt tS_{\setminus i}+\theta w_i).
\end{align*}
Taking expectations on both sides and applying \eqref{eq:s7-layer} leads to the stated result.
\end{proof}

\subsection{Concluding the proof of Theorem~\ref{thm:bounded}}\label{subsec:induction}
\noindent The probability $\pi_t$ in Proposition~\ref{prop:pertime} involves only the $k-1$ vectors $X_2,\dots,X_k$. This allows for a bound in terms of $\bar{\mathsf D}_{k-1}$ and eventually leads to a finite induction on $k$.

\begin{lemma}\label{lem:induction} 
Let $t\in[1/2,1]$, and let $r_t$ and $\pi_t$ be as in Proposition~\ref{prop:pertime}. Then,
\[
\pi_t\lesssim r_t\sqrt{\log(pn)}+\bar{\mathsf D}_{k-1}.
\]
\end{lemma}
\begin{proof}
Let $S_{k-1}=(k-1)^{-1/2}\sum_{i=1}^{k-1}X_i$, and note that $\sqrt tS_{\setminus 1}$ is equal in distribution to $\sqrt{t(k-1)/k}\,S_{k-1}$. Also observe that if the condition $\abs{\min_j(\beta_j-y_j)}<r_t$ holds, then so does
$y\le\beta+r_t\1$ but not $y\le\beta-r_t\1$. Hence, $\pi_t$ is at most the difference of the probabilities that $S_{k-1}$ lies in two orthants, and replacing $S_{k-1}$ by $G$ in each of them incurs an error of at most $2\bar{\mathsf{D}}_{k-1}$. This leads to
\[
\pi_t\le\PP\Big(G\le\frac{\sqrt k\,(\beta+r_t\1)}{\sqrt{t(k-1)}}\Big)-\PP\Big(G\le\frac{\sqrt k\,(\beta-r_t\1)}{\sqrt{t(k-1)}}\Big)+2\bar{\mathsf D}_{k-1}.
\]
The corners of the orthants in these two probabilities  differ by $2\sqrt{k/(t(k-1))}\,r_t\1\leq 4r_t \1$. So, by Nazarov's inequality (Lemma~\ref{fact:nazarov}), the difference of the Gaussian
probabilities is at most of order $r_t\sqrt{\log(pn)}$.
\end{proof}

We now derive a recursion that upper bounds $\bar{\mathsf{D}}_k$ in terms of $\bar{\mathsf{D}}_{k-1}$.  Let $k=2,\dots,n$ continue to satisfy \eqref{eq:small}, so that the results of this section apply, and let  $\beta\in\R^p$ still be fixed. We will upper bound $\int_0^1 |\Psi'(t)|dt$ by integrating the derivative bound in
Proposition~\ref{prop:pertime}. Its last term contributes $\frac{1}{pn}\int_0^1 t^{-1/2} dt\le\frac{\bar B}{\sqrt k}$, which is crude but adequate.
Integrating the first term over $(0,1/2]$ leads to a contribution of order $(\bar B/\sqrt k)(\log(pn))^{3/2}$. When $t\in[1/2,1]$ we use
Lemma~\ref{lem:induction}, which gives
\[
u(t)^{-3/2}\big(\pi_t+(pn)^{-2}\big)\lesssim\frac{\log(pn)}{u(t)}+\frac{\bar B\sqrt{\log(pn)}}{\sqrt k\,u(t)^{3/2}}+\frac{\bar{\mathsf D}_{k-1}+(pn)^{-2}}{u(t)^{3/2}}.
\]
Combining the last bound with
$
\int_{1/2}^1\frac{dt}{u(t)}\le\log(n)$ 
and
$\int_{1/2}^1\frac{dt}{u(t)^{3/2}}\le\frac2\sigma$, we have
\begin{equation*}
\begin{split}
\int_0^1|\Psi'(t)|\,dt&\lesssim\frac{\bar B}{\sqrt k}(\log(pn))^{3/2}\Big(\log(pn)\log(n)+\frac{\bar B\sqrt{\log(pn)}}{\sqrt k\,\sigma}+\frac{\bar{\mathsf D}_{k-1}+(pn)^{-2}}{\sigma}\Big).
\end{split}
\end{equation*}
Furthermore, the smoothing error in~\eqref{eq:smoothing} satisfies the bound
\begin{equation}
    \begin{split}
        \sigma\log(pn)
        %=CC_0\bar B(\log(pn))^{5/2}/\sqrt k\\
        %
        &\leq \frac{C_0}{\log (2)}\,\frac{\bar B(\log(pn))^{5/2}\log(n)}{\sqrt k}, 
    \end{split}
\end{equation}
because $\log(n)\geq\log(2)$.
 Combining this with~\eqref{eq:interpFTC}, the previous display, and $C_0\geq 1$, we obtain an absolute
constant $C_1$, not depending on $C_0$, such that
\begin{equation}\label{eq:s7-rec}
\bar{\mathsf D}_k\le C_1C_0\frac{\bar B(\log(pn))^{5/2}\log(n)}{\sqrt k}+C_1C_0^{-1}\,\bar{\mathsf D}_{k-1},
\end{equation}
for every $k=2,\dots,n$ satisfying \eqref{eq:small}.

We conclude with a finite induction on $k=1,\dots,n$.
As $C_1$ does not depend on $C_0$, we may increase the value of $C_0$ if necessary  so that $C_1C_0^{-1}\le\frac14$. Letting $M=\max(4,2C_1C_0)$, we will show $\bar{\mathsf D}_k\le M\bar B(\log(pn))^{5/2}\log(n)/\sqrt k$ for
$k=1,\dots,n$, which will give the stated result. The base case of  $k=1$ is clear. Now, let $k\geq 2$ and assume the claim holds for $k-1$. If \eqref{eq:small} fails, then the
claimed bound is at least $M/4\ge1\ge \bar{\mathsf D}_k$. Otherwise, the following bounds hold by \eqref{eq:s7-rec}, $C_1C_0^{-1}\le\frac14$, and the induction hypothesis
\[
\begin{split}
\bar{\mathsf D}_k&\le \frac{C_1C_0\bar B(\log(pn))^{5/2}\log(n)}{\sqrt k}+\frac{\frac{1}{4}M\bar B(\log(pn))^{5/2}\log(n)}{\sqrt{k-1}}\\
&\le\Big(\frac12+\frac{\sqrt2}4\Big)\frac{M\bar B(\log(pn))^{5/2}\log(n)}{\sqrt k}\\
&\le\frac{M\bar B(\log(pn))^{5/2}\log(n)}{\sqrt k},
\end{split}
\]
where we have used $\sqrt{k/(k-1)}\leq \sqrt 2$ and $C_1C_0\leq M/2$ in the second step. This completes the proof of Theorem~\ref{thm:bounded}. \qed

\section{Proof of Proposition~\ref{prop:persummand}}\label{sec:onesummand}

Before proceeding directly to the proof of Proposition~\ref{prop:persummand} in Section~\ref{subsec:persumproof}, we will develop a collection of preparatory results in Sections~\ref{sec:properties} and~\ref{subsec:regular}.

\subsection{Properties of facet functions and the facet sensitivity}\label{sec:properties}

We begin by introducing some parameters associated with the sections of polyhedra.\\

\noindent\textbf{Sections of polyhedra and related parameters.} Let $P$ be a polyhedron in $\R^p$ with facets $F_1,\dots,F_p$ having linearly independent normal vectors $n_1,\dots,n_p$.
For each $j\leq p$, let the $(p-1)$-dimensional subspace orthogonal to $n_j$ be denoted as $E_j=n_j^\perp$, and let $\Pi_j$ be the orthogonal projection onto $E_j$. The $j$th section of $P$ and its Gaussian measure are denoted as
\[
K_j=\{z\in E_j:a_jn_j+z\in F_j\}\qquad\text{and}\qquad q_j=\gamma_{E_j}(K_j).
\]
It will be shown in Lemma~\ref{lem:structure} that $q_j>0$. Also, let $W$ be a random
vector in $E_j$ with distribution $\gamma_{E_j}$, and let
\[
m_j=\E(W\mid W\in K_j)\qquad\text{and}\qquad \mathfrak S_j=\Cov(W\mid W\in K_j).
\]
%
%

% preamble: \usepackage{tikz} \usepackage{graphicx}

The next lemma relates $\nu_j(P)$ to $q_j$ and is proved in Section~\ref{subsec:polyhedra}.

\begin{lemma}\label{lem:structure} Let $P\subset\R^p$ be a polyhedron as in~\eqref{eqn:Polygeneral} with $p$ facets having linearly independent normal vectors.
Then, for any $j\leq p$, the facet $F_j$ can be represented as $F_j=a_jn_j+K_j$, and $K_j$ is a closed convex
subset of $E_j$ with Gaussian measure $q_j>0$. Moreover,
\[
\nu_j(P)=\phi(a_j)\,q_j.
\]
\end{lemma}

Having introduced the parameters associated with sections, the facet sensitivity of Definition~\ref{def:sensitivity} can now be expressed in the following explicit form.

\begin{lemma}\label{lem:sensitivityformula}
If $P$ is as in the statement of Lemma~\ref{lem:structure}, then
\[
\Fc(P)=\sum_{j=1}^p\big(1+a_j^2+\norm{m_j}_2^2+\tr(\Pi_j-\mathfrak S_j)\big)\,\nu_j(P).
\]
\end{lemma}
\begin{proof}
Lemma~\ref{lem:facetlog} will show that the gradient and Hessian of the function $x\mapsto\log\nu_j(P-x)$ at $x=0$ are respectively given by $a_jn_j+m_j$ and 
$-(n_jn_j^\top+\Pi_j-\mathfrak S_j)$. As $m_j\in E_j$ is orthogonal to $n_j$, the squared norm of the gradient is $a_j^2+\norm{m_j}_2^2$, and the
trace of minus the Hessian is $1+\tr(\Pi_j-\mathfrak S_j)$. Applying these observations to the formula~\eqref{eqn:facetsensitivity} completes the proof.
\end{proof}

We now return to the polyhedra $P(y)$ of Section~\ref{subsec:sensitivity}. For the objects of Sections~\ref{subsec:sensitivity} and~\ref{sec:properties} attached to
$P(y)$, we write 
$$\nu_j(y)=\nu_j(P(y))$$
and similarly for $q_j(y)$, $m_j(y)$ and $\mathfrak{S}_j(y)$. Also recall that
$\psi_j(y)=\nu_j(y)/\sqrt{u+\sigma^2}$.
The importance of the next lemma is that it expresses the derivatives of $g_u$ in the interpolation 
through the polyhedron $P(y)$.

\begin{lemma}\label{lem:gu}
Let $u>0$ and $y,w\in\R^p$. Also, let the normal vectors $n_1,\dots,n_p$ and the polyhedron $P(y)$ be as in~\eqref{eq:Py}. Then, the following statements hold.
\begin{itemize}
\item[(a)] We have $P(y+w)=P(y)-\Sigma_u^{-1/2}w$.
\item[(b)] We have
\begin{align*}
\partial_wg_u(y)&=-\sum_{j=1}^p w_j\psi_j(y),\\[2mm]
\LapS g_u(y)&=-\sum_{j=1}^p(\Sigma\nabla\psi_j(y))_j=\frac1u\tr\big(\Gamma\,\nabla_x^2\gamma_p(P(y)-x)|_{x=0}\big),
\end{align*}
where $\Gamma=u\,\Sigma_u^{-1/2}\Sigma\,\Sigma_u^{-1/2}$ satisfies $0\preceq\Gamma\preceq I_p$.
\end{itemize}
\end{lemma}
\begin{proof}
Part (a) holds because $z\in P(y+w)$ if and only if $\Sigma_u^{1/2}(z+\Sigma_u^{-1/2}w)\leq \beta-y$.
Turning to part (b), fix any $y\in \R^p$. From part (a) above and Lemma~\ref{lem:transl}(a), we have
\[
\partial_wg_u(y)=-\sum_{j=1}^p\frac{w_j}{\sqrt{u+\sigma^2}}\,\nu_j(y)=-\sum_{j=1}^pw_j\psi_j(y).
\]
Hence
$\LapS g_u=\sum_{j,l}\Sigma_{jl}\partial_l\partial_jg_u=-\sum_j(\Sigma\nabla\psi_j)_j$. In addition, we have the Hessian formula
$$\nabla^2g_u(y)=\Sigma_u^{-1/2}\,\nabla_x^2\gamma_p(P(y)-x)|_{x=0}\,\Sigma_u^{-1/2},$$
which gives the second formula for
$\LapS g_u(y)=\tr(\Sigma\nabla^2g_u(y))$. Finally, $0\preceq\Gamma\preceq I_p$ follows from $0\preceq u\Sigma\preceq\Sigma_u$ by conjugation with
$\Sigma_u^{-1/2}$.
\end{proof}

The following lemma records further properties of the facet functions and is proved in Section~\ref{subsec:polyhedra}.
\begin{lemma}\label{lem:facetfn}
If $u>0$ and $j\le p$, then for all $y\in\R^p$, we have
\[
\psi_j(y)=\frac{\phi(a_j(y))\,q_j(y)}{\sqrt{u+\sigma^2}},
\]
and the functions $\psi_j$ and $q_j$ are positive, $C^\infty$ and log-concave on $\R^p$. In particular, for all $y,w\in\R^p$,
\begin{equation}\label{eq:s7-Dlogpsi}
\partial_w\log\psi_j(y)=\frac{w_j\,a_j(y)}{\sqrt{u+\sigma^2}}+\partial_w\log q_j(y).
\end{equation}
Moreover, for all $y\in\R^p$, we have
\begin{align}
\ip{\nabla\log\psi_j(y)}{\Sigma\,\nabla\log\psi_j(y)}&\le\frac{a_j(y)^2+\norm{m_j(y)}_2^2}u,\label{eq:facetfn1}\\
\tr\big(\Sigma(-\nabla^2\log\psi_j(y))\big)&\le\frac{1+\tr(\Pi_j-\mathfrak S_j(y))}u,\label{eq:facetfn2}
\end{align}
as well as,
\begin{align}
\sum_{j=1}^p\psi_j(y)&\le u^{-1/2}\Fc(y),\label{eq:facetfn3}\\
\sum_{j=1}^p\psi_j(y)\,\ip{\nabla\log\psi_j(y)}{\Sigma\,\nabla\log\psi_j(y)}&\le u^{-3/2}\Fc(y),\label{eq:facetfn4}\\
\sum_{j=1}^p\psi_j(y)\tr\big(\Sigma(-\nabla^2\log\psi_j(y))\big)&\le u^{-3/2}\Fc(y).\label{eq:facetfn5}
\end{align}
\end{lemma}

The second term in~\eqref{eq:s7-Dlogpsi} is controlled by the following lemma, which is proved in Section~\ref{subsec:sectionsproof}.

\begin{lemma}\label{lem:hazard}
For all $u\in(0,2]$, $y,w\in\R^p$ and $j\le p$,
\[
\abs{\partial_w\log q_j(y)}\lesssim\frac{\norm w_\infty}\sigma\Big(\sqrt{\log(pn)}\,\log\frac1{q_j(y)}+1\Big).
\]
\end{lemma}

\subsection{Regular facets}\label{subsec:regular}

This subsection distinguishes facets into two types that we call \emph{regular} and \emph{non-regular} (Definition~\ref{def:regular}). Roughly speaking, a facet $F_j$ is regular when its offset parameter $a_j$ is not very large, and its Gaussian section measure $q_j$ is non-negligible.  
One reason why this matters is that if $F_j$ is a regular facet, then it will turn out that the curvature of $\log\psi_j$ can be controlled (Proposition~\ref{lem:curv}), and this will be
essential for the proof of Proposition~\ref{prop:persummand}.

\begin{definition}[regular facets]\label{def:regular}
For $u>0$, $x\in\R^p$, $j\leq p$, we say the facet $F_j$ (or the index $j$) is regular at $x$ if
\[
\abs{a_j(x)}\le\sqrt{20\log(pn)}\qquad\text{and}\qquad q_j(x)\ge(pn)^{-10}.
\]
\end{definition}

If the $j$th facet is regular at $x$, then the proof of Proposition~\ref{prop:persummand} relies on the curvature of $\log\psi_j$ being well controlled at points near $x$, which is ensured by the following result.

\begin{proposition}\label{lem:curv}
Let $u\in(0,2]$, $j\le p$ and $x,w\in\R^p$ be fixed with $\norm w_\infty\le \bar B/\sqrt k$. If $j$ is regular at $x$, then the function  $f_j(\tau)=\log\psi_j(x+\tau w)$ satisfies
\[
-f_j''(\tau)\lesssim-f_j''(0)+\frac{\bar B^2}{kpn}\qquad\text{for all }\tau\in[0,1].
\]
\end{proposition}
 \noindent This result is crucial, and it will be explained in Section~\ref{sec:Gibbs} that it is a consequence of Proposition~\ref{prop:kappacomp}.

\subsection{Concluding the proof  of Proposition~\ref{prop:persummand}}\label{subsec:persumproof}

Before starting the proof, we need the following lemma, which records some properties about regular facets and is proved in Section~\ref{sec:remaining}. 

\begin{lemma}\label{lem:segment}
Let $j\leq p$, $u\in(0,2]$ and $x,w\in\R^p$ with $\norm w_\infty\le \bar B/\sqrt k$.
\begin{itemize}
\item[(a)] If $j$ is regular at $x$, then $\abs{\partial_w\log\psi_j(x)}\le\frac1{100}$.
\item[(b)] If $j$ is not regular at $x$, then for all $\tau\in[0,1]$ we have
$$\big(\abs{a_j(x)}+\log\frac1{q_j(x)}+1\big)\,\nu_j(x)\lesssim(pn)^{-5} \qquad \text{ and } \qquad \nu_j(x+\tau w)\lesssim(pn)^{-5}.$$
\end{itemize}
\end{lemma}

For the purposes of proving Proposition~\ref{prop:persummand}, the objects $t\in(0,1]$, $u=u(t)$ and $x\in\R^p$ will remain fixed throughout this subsection. For each $j\le p$ and $w\in\R^p$, define the linear function $h_j(w)=\ip{\nabla\log\psi_j(x)}{w}$ and the remainder function
\begin{equation}\label{eqn:remj}
R_j(w)=\psi_j(x+\sqrt tw)-\psi_j(x)-\sqrt t\,\psi_j(x)h_j(w).
\end{equation}
 By Lemma~\ref{lem:gu}(b), we have 
\begin{equation}\label{eq:mainrem}
\ip{\nabla g_u(x+\sqrt tw)-\nabla g_u(x)}{w}=-\sqrt t\sum_{j=1}^pw_j\psi_j(x)h_j(w)-\sum_{j=1}^pw_jR_j(w).
\end{equation}
We will apply this with $w=X_1/\sqrt k$ and take expectations. Since $h_j$ is linear, we have $\E\big(X_{1j}h_j(X_1)\big)=(\Sigma\nabla\log\psi_j(x))_j$. Together with $\psi_j\nabla\log\psi_j=\nabla\psi_j$ and
Lemma~\ref{lem:gu}(b), this shows that the first term on the right of~\eqref{eq:mainrem} has expectation
\begin{equation}\label{eqn:mainterm}
\begin{split}
-\frac{\sqrt t}k\sum_{j=1}^p\psi_j(x)\,\E\big(X_{1j}h_j(X_1)\big)
&=-\frac{\sqrt t}k\sum_{j=1}^p(\Sigma\nabla\psi_j(x))_j = \frac{\sqrt t}k\LapS g_u(x).
\end{split}
\end{equation}
It therefore remains to bound the sum
\[
\frac1{\sqrt k}\sum_{j=1}^p\E\abs{X_{1j}R_j(X_1/\sqrt k)},
\]
and this will be done by handling the regular and non-regular indices $j$ separately.
\\

\noindent\textit{Remainder bound for regular facets.} Let $j$ be regular at $x$, and first let $w\in\R^p$ be deterministic with $\norm w_\infty\le\bar B/\sqrt k$. By
Lemma~\ref{lem:facetfn}, the function $f_j(\tau)=\log\psi_j(x+\tau w)$ is $C^\infty$ and concave with $f_j'(0)=h_j(w)$, so Taylor's formula with
integral remainder leads to
\[
\psi_j(x+\sqrt tw)=\psi_j(x)\,e^{\sqrt th_j(w)-\varepsilon_j(w)}\qquad\text{where }\qquad \varepsilon_j(w)=\int_0^{\sqrt t}(\sqrt t-r)\,(-f_j''(r))\,dr\ge0 ,
\]
where we recall that $t\in(0,1]$.
Consequently, we may express the remainder $R_j(w)$ defined in~\eqref{eqn:remj} as
\begin{equation*}
R_j(w)=\psi_j(x)\Big[e^{\sqrt t h_j(w)-\varepsilon_j(w)}-1-\sqrt t h_j(w)\Big].
\end{equation*}
By Lemma~\ref{lem:segment}(a), we have $\abs{h_j(w)}\le\frac1{100}$, and then applying the inequality $\abs{e^{a-b}-1-a}\le a^2+2b$ for $\abs a\le\frac1{100}$
and $b\ge0$ gives 
\begin{equation}\label{eqn:midrem}
\abs{R_j(w)}\le\psi_j(x)\Big[th_j(w)^2+2\varepsilon_j(w)\Big].
\end{equation}
Next, a step that deserves special emphasis is that $\varepsilon_j(w)$ can be bounded pointwise in $w$ by
Proposition~\ref{lem:curv}. Combining that result with $f_j''(0)=\ip{w}{\nabla^2\log\psi_j(x)\,w}\leq 0$  gives
$$ \varepsilon_j(w)\lesssim -t\ip{w}{\nabla^2\log\psi_j(x)\,w}+ \frac{t\bar B^2}{kpn}.$$
%$$-f_j''(r)\lesssim-f_j''(0)+\bar B^2/(kpn)$$
%
After applying the last bound to~\eqref{eqn:midrem} and substituting $w=X_1/\sqrt k$, we take the expectation over $X_1$ to obtain
\begin{equation*}
\begin{split}
    \frac{1}{\sqrt k}\E|X_{1j} R_j(X_1/\sqrt k)| &\lesssim \frac{t \psi_j(x)\bar B}{\sqrt k}\E\bigg( h_j(X_1/\sqrt k)^2-\frac{1}{k} \ip{X_1}{\nabla^2 \log \psi_j(x)\, X_1} +\frac{\bar B^2}{kpn}\bigg)\\[0.2cm]
    &= \frac{t \psi_j(x)\bar B}{ k^{3/2}}\bigg(\ip{\nabla\log\psi_j(x)}{\Sigma\,\nabla\log\psi_j(x)}+ \tr\big(\Sigma(-\nabla^2\log\psi_j(x))\big)
    +\frac{\bar B^2}{pn}\bigg),
\end{split}
\end{equation*}
where we have used $\|X_1\|_{\infty}\leq \bar B$ in the first step. By summing the terms on the right over regular $j$ and then applying the facet function bounds~\eqref{eq:facetfn3}-\eqref{eq:facetfn5} of Lemma~\ref{lem:facetfn} (which hold for all $j\leq p$), we conclude that 
\begin{equation}\label{eqn:regularsum}
\begin{split}
\frac1{\sqrt k}\sum_{j\text{ regular}}\E\abs{X_{1j}R_j(X_1/\sqrt k)}
&\lesssim\frac{t\bar B\,\Fc(x)}{(ku)^{3/2}},
\end{split}
\end{equation}
where we have simplified the last bound by using $\bar B^2/(pn)\le1$ and $u^{-1/2}\leq 2u^{-3/2}$ (as $u\leq 2$).\\

\noindent\textit{Remainder bound for non-regular facets.} Let $j$ be non-regular at $x$, and again let $w\in\R^p$ be deterministic with $\norm w_\infty\le\bar B/\sqrt k$.
Since the definition of $\psi_j$ implies the pointwise upper bound $\psi_j=\nu_j/\sqrt{\sigma^2+u}\leq \nu_j/\sigma$, we may apply
Lemma~\ref{lem:segment}(b) to the first two terms in the formula~\eqref{eqn:remj}. Together with $t\le1$, this gives
\begin{equation*}|R_j(w)| \lesssim \frac{1}{\sigma (pn)^5} + \psi_j(x)|h_j(w)|.
\end{equation*}
The second term on the right can be bounded by combining the formula~\eqref{eq:s7-Dlogpsi} for $\partial_w\log\psi_j(x)=h_j(w)$ with
Lemma~\ref{lem:hazard}, and using $\sqrt{\log(pn)}\norm w_\infty/\sigma\le 1$ from~\eqref{eq:ratio}. Along with
the pointwise bound $\psi_j\le\nu_j/\sigma$ and Lemma~\ref{lem:segment}(b), this gives
\begin{align*}
\psi_j(x)\abs{h_j(w)}
&\lesssim\frac{\nu_j(x)}{\sigma}\Big(\abs{a_j(x)}+\textstyle \log\frac1{q_j(x)}+1\Big)
\lesssim \, \frac{1}{\sigma(pn)^5}.
\end{align*}
In turn, using the last two displays with $\norm w_\infty\leq \bar B/\sqrt k\le\sigma$ leads to
$$\abs{w_jR_j(w)}\lesssim \frac{1}{(pn)^5}.$$
Since this holds for every such $w$, it holds in particular for $w=X_1/\sqrt k$. Taking the expectation over $X_1$, and summing over the set of 
(at most $p$) non-regular indices $j$, we have
\[
\frac1{\sqrt k}\sum_{j\text{ not regular}}\E\abs{X_{1j}R_j(X_1/\sqrt k)}\lesssim \frac1{(pn)^4}.
\]
Combining this with the bound~\eqref{eqn:regularsum} for the regular facets and the expression~\eqref{eqn:mainterm} for the main term, the proof of
Proposition~\ref{prop:persummand} is complete.\qed

\section{Proof of Proposition~\ref{lem:curv}: control of curvature}\label{sec:HS}
The curvature of the log-facet functions $\log \psi_j$ appearing in the statement of Proposition~\ref{lem:curv} is difficult to bound in a direct manner. This section reframes the curvature in terms of the log-partition function of a particular Gibbs measure, which is introduced in Section~\ref{sec:Gibbs}. The benefit of this reframing is that the curvature can be expressed in terms of a variance with respect to the Gibbs measure, and well-developed tools are available for analyzing such variances. Throughout this section, the quantities $u$, $k$, and $\sigma=C_0\bar B(\log(pn))^{3/2}/\sqrt k$ will remain fixed as in Section~\ref{sec:reductions}, and so they will be suppressed in much of the notation.

\subsection{Reframing in terms of log-partition functions of Gibbs measures}\label{sec:Gibbs}
Let $\fa=(\fa_1,\dots,\fa_{p-1})$ be a tuple of real-valued affine functions on $\R^{p-1}\times \R$. That is, for each $l=1,\dots,p-1$, the function $\fa_l(z,\tau)$ is jointly affine in $z\in\R^{p-1}$ and $\tau\in \R$. For a given tuple $\fa$, consider the associated  function
\begin{equation}\label{eq:elldef}
\ell_\tau(z)=-\sum_{l=1}^{p-1}\log\Phi\Big(\frac{\fa_l(z,\tau)}\sigma\Big),
\end{equation}
which is non-negative and convex. As such, this function induces a Gibbs measure $\mu_{\tau}$ whose density with respect to $\gamma_{p-1}$ on $\R^{p-1}$ is defined by
\begin{equation}\label{eqn:Gibbs}
d\mu_\tau(z)=e^{-\ell_\tau(z)-A(\tau)}\,d\gamma_{p-1}(z),
\end{equation}
where $A(\tau)$ is the \emph{log-partition function}. More explicitly, the function $A(\tau)$ is given by
\begin{equation}\label{eq:Adef}
A(\tau)=\log\int_{\R^{p-1}}\prod_{l=1}^{p-1}\Phi\Big(\frac{\fa_l(z,\tau)}\sigma\Big)\,d\gamma_{p-1}(z).
\end{equation}
Below, we will usually suppress the fact that  $\mu_{\tau}$ and $A$ depend on $\fa$, and will simply refer to them as the Gibbs measure and log-partition function associated with $\fa$.

The close connection between log-facet functions and log-partition functions is established in Lemma~\ref{lem:sections} below, which is proved in Section~\ref{subsec:sectionsproof}. As a clarification about notation, the gradient $\nabla \fa_l$ below is understood \emph{with respect to the $z$ argument} of $\fa_l(z,\tau)$, and it should also be noted that the vectors $\nabla\fa_l$ and  $\partial_{\tau}\fa=(\partial_{\tau} \fa_1,\dots,\partial_{\tau}\fa_{p-1})$ are constant with respect to $z$ and $\tau$.

\begin{lemma}\label{lem:sections}
For any fixed $u\in(0,2]$, $j\le p$ and $y,w\in\R^p$, there is a tuple of affine functions $\fa=(\fa_1,\dots,\fa_{p-1})$ on $\R^{p-1}\times \R$ with an associated log-partition function $A$
such that 
\begin{equation}\label{eqn:curvatureidentity}
\partial_{\tau}^2\log\psi_j(y+\tau w)=-\frac{w_j^2}{u+\sigma^2}+A''(\tau)\qquad\text{for all }\tau\in\R.
\end{equation}
Furthermore, the tuple $\fa$ and function $A$ satisfy
$\norm{\partial_\tau\fa}_\infty\le2\norm w_\infty$, $\max_{l\le p-1}\norm{\nabla\fa_l}_2\le2$, and $A(0)=\log q_j(y)$.
\end{lemma}

Based on this lemma, it is straightforward to check that Proposition~\ref{lem:curv} is implied by the following analogous statement involving the log-partition function, which is proven in Section~\ref{subsec:HScomp}.

\begin{proposition}\label{prop:kappacomp}
Let $\fa=(\fa_1,\dots,\fa_{p-1})$ be a tuple of affine functions on $\R^{p-1}\times \R$ with associated log-partition function $A$. Also, let $n\ge2$ be an integer with $pn$ sufficiently large. If $\norm{\partial_\tau\fa}_\infty\le2\bar B/\sqrt k$,
$\max_{l\le p-1}\norm{\nabla\fa_l}_2\le2$,  and $A(0)\ge-10\log(pn)$, then
\[
-A''(\tau)\lesssim-A''(0)+\frac{\bar B^2}{kpn}\qquad\text{for all }\tau\in[0,1].
\]
\end{proposition}

\subsection{Variational representation of $-A''(\tau)$}\label{subsec:HSstrategy}

The main goal of this subsection is to develop a formula showing that for a fixed $\tau\in \R$, the negative curvature $-A''(\tau)$ can be represented as the minimum value of a functional of the form $h\mapsto \E_{\mu_{\tau}} \mathsf J_{\tau} h$, where $\mathsf{J}_{\tau}$ is a certain differential operator and $h$ ranges over a certain set of test functions~(see \eqref{eqn:Aformula} and~\eqref{eq:Jgamma}). This representation of $-A''(\tau)$ is useful, because a minimum value lends itself to deriving an upper bound on $-A''(\tau)$ that is needed in Proposition~\ref{prop:kappacomp}, and because the functional will have a well-behaved dependence on $\tau$ (which facilitates the comparison of $A''(\tau)$ and $A''(0)$).
Throughout the rest of this section, $\fa=(\fa_1,\dots,\fa_{p-1})$ is a tuple of affine functions on $\R^{p-1}\times\R$, with associated function
$\ell_\tau$, Gibbs measure $\mu_\tau$ and log-partition function $A$ as in Section~\ref{sec:Gibbs}.

The curvature of $A(\tau)$ can be calculated by differentiating the identity $e^{A(\tau)}=\int e^{-\ell_\tau}\,d\gamma_{p-1}$ twice, leading to
\begin{equation}\label{eq:bartlett}
-A''(\tau)=\E_{\mu_\tau}\partial_\tau^2\ell_\tau-\Var_{\mu_\tau}(\dot\ell_\tau),
\end{equation}
where we use the notation $\dot\ell_{\tau}=\partial_{\tau}\ell_{\tau}$.
 When the density of $\mu_{\tau}$ is computed with respect to Lebesgue measure on $\R^{p-1}$ (rather than $\gamma_{p-1}$), it is proportional to $e^{-W_{\tau}}$, where the function
\[
W_\tau(z)=\tfrac12\norm z_2^2+\ell_\tau(z)
\]
is jointly convex in $(z,\tau)$, by the convexity of $\ell_\tau$. Another related classical fact is that $\mu_{\tau}$ is the invariant measure of the Langevin diffusion (see, e.g., \citealp{Pav14,BGL14}),
whose generator is given by the following operator $\mathsf{L}_{\tau}$ defined in terms of $W_{\tau}$,
\[
\Lop_\tau h=\Delta h-\ip{\nabla W_\tau}{\nabla h},
\]
where $h$ is a suitable test function on $\R^{p-1}$. This operator and the function $W_{\tau}$ are needed to obtain the variational representation for $\Var_{\mu_{\tau}}(\dot\ell_{\tau})$ in the lemma below, which is proved in Section~\ref{subsec:HSvar}. As a matter of notation, we use $\mathscr H$ to denote the class of functions
$h\in C^3(\R^{p-1})$ with bounded first and second derivatives.

\begin{lemma}\label{lem:HSvar}
For any $\tau\in\R$ and any $h\in\mathscr H$, we have
\begin{equation}\label{eqn:Helffer}
\begin{split}
\Var_{\mu_\tau}(\dot\ell_\tau)&=\E_{\mu_\tau}\Big(2\ip{\nabla\dot\ell_\tau}{\nabla h}-\ip{\nabla h}{\nabla^2W_\tau\nabla h}-\norm{\nabla^2h}_F^2\Big)\\
&\quad+\E_{\mu_\tau}\big[(\dot\ell_\tau-\E_{\mu_\tau}\dot\ell_\tau+\Lop_\tau h)^2\big].
\end{split}
\end{equation}
In particular, $\Var_{\mu_\tau}(\dot\ell_\tau)$ is lower bounded by the first term on the right, with equality holding when $h$ solves the equation
$\Lop_\tau h=\E_{\mu_\tau}\dot\ell_\tau-\dot\ell_\tau$.
\end{lemma}
The formula~\eqref{eqn:Helffer} may be regarded as a variant of the Helffer--Sj\"ostrand representation for covariances computed with respect to a Gibbs
measure~\citep[p.352]{HS94}. This formula and~\eqref{eq:bartlett} motivate defining the operator
\begin{equation}\label{eq:Jgamma}
\J_\tau h=\partial_\tau^2\ell_\tau-2\ip{\nabla\dot\ell_\tau}{\nabla h}+\ip{\nabla h}{\nabla^2W_\tau\nabla h}+\norm{\nabla^2h}_F^2,
\end{equation}
which allows $-A''(\tau)$ to be expressed as
\begin{equation}\label{eqn:Aformula}
    -A''(\tau)=\E_{\mu_{\tau}} \mathsf{J}_{\tau} h -\E_{\mu_\tau}\big[(\dot\ell_\tau-\E_{\mu_\tau}\dot\ell_\tau+\Lop_\tau h)^2\big]
\end{equation}
for every $\tau\in\R$ and every $h\in\mathscr H$. This formula for $-A''(\tau)$ now puts us in position to prove Proposition~\ref{prop:kappacomp}.

\subsection{Proof of Proposition~\ref{prop:kappacomp}}\label{subsec:HScomp}
The variational representation~\eqref{eqn:Aformula} shows that the upper bound
\begin{equation}\label{eq:Aupper}
    -A''(\tau)\leq \E_{\mu_{\tau}} \mathsf{J}_{\tau} h
\end{equation}
holds for every $h\in\mathscr H$ and $\tau\in[0,1]$. It also shows that if $h^{\star}\in\mathscr H$ solves the equation $\Lop_{0} h^{\star}=\E_{\mu_{0}}\dot\ell_{0}-\dot\ell_{0}$, then we have the equality
\begin{equation*}-A''(0)=\E_{\mu_{0}} \mathsf{J}_{0} h^{\star}.
\end{equation*}
(The fact that such a solution $h^{\star}$ exists is ensured by Lemma~\ref{lem:varrep} later.)
To proceed, we take the test function in~\eqref{eq:Aupper} to be $h^{\star}$, which gives
\begin{equation*}-A''(\tau)\leq \E_{\mu_{\tau}} \mathsf{J}_{\tau} h^{\star},
\end{equation*}
for all $\tau\in[0,1]$.
Thus, to prove Proposition~\ref{prop:kappacomp}, it remains to relate $\E_{\mu_\tau}\J_\tau h^\star$ to $\E_{\mu_0}\J_0h^\star$. We do this in two steps, by replacing $\J_{\tau}$ with $\J_0$, and then
replacing $\mu_{\tau}$ with $\mu_0$, while keeping track of the errors introduced by these substitutions. First, by choosing the test function in Lemma~\ref{lem:Jcomp} to be $h^{\star}$, it follows that
\begin{equation}\label{eq:stepi}
\E_{\mu_\tau}\J_\tau h^\star\le2\,\E_{\mu_\tau}\J_0h^\star+\frac{2\norm{\partial_\tau\fa}_\infty^2}{(pn)^4}.
\end{equation}
Next, it is shown in Lemma~\ref{lem:stab} that
\begin{equation}\label{eq:stepii}
\E_{\mu_\tau}\J_0h^\star\lesssim\E_{\mu_0}\J_0h^\star+\frac{\norm{\J_0h^\star}_{\sup}}{(pn)^{10}}.
\end{equation}
Putting the last few pieces together gives
\begin{equation*}-A''(\tau)\lesssim -A''(0)+\frac{\norm{\J_0h^\star}_{\sup}}{(pn)^{10}}+\frac{\norm{\partial_\tau\fa}_\infty^2}{(pn)^4}.
\end{equation*}
In Lemma~\ref{lem:varrep}, it will be shown that the function $h^{\star}$ can be taken to satisfy $\norm{\J_0h^\star}_{\sup}\lesssim(pn)^5\norm{\partial_\tau\fa}_\infty^2$. Finally, applying the assumption $\norm{\partial_\tau\fa}_\infty\le2\bar B/\sqrt k$ completes the proof of
Proposition~\ref{prop:kappacomp}.\qed

\subsection{The dependence of $\J_\tau$ on the parameter $\tau$}\label{subsec:derivs}

This subsection will establish the upper bound~\eqref{eq:stepi} in Lemma~\ref{lem:Jcomp}, which relies on the explicit formula for $\mathsf{J}_{\tau}$ below.  In particular, the formula shows that the dependence of $\J_\tau h$ on $\tau$ is essentially governed by the function $\rho:=-(\log\Phi)''$, which has well understood properties.

\begin{lemma}\label{lem:Jformula}
For every tuple $\fa$ of affine functions on $\R^{p-1}\times\R$, every $\tau\in\R$, and every function $h\in C^2(\R^{p-1})$, we have
\begin{equation}\label{eq:Jdef}
\J_\tau h \, = \, \norm{\nabla h}_2^2+\norm{\nabla^2h}_F^2+ \sum_{l=1}^{p-1}\rho\Big(\frac{\fa_l(\cdot,\tau)}\sigma\Big)\frac{(\partial_\tau\fa_l-\ip{ \nabla\fa_l}{\nabla h})^2}{\sigma^2}.
\end{equation}
Furthermore, $\J_\tau h\ge0$.
\end{lemma}

\begin{proof} From looking at the definition of $\mathsf J_{\tau}$ in~\eqref{eq:Jgamma} and noting that $\nabla^2 W_{\tau}=I+\nabla^2 \ell_{\tau}$, we need to calculate $\partial_{\tau}^2\ell_{\tau}$, $\nabla \dot \ell_{\tau}$ and $\nabla^2 \ell_{\tau}$.
 Direct calculations with $\ell_{\tau}$ give
\begin{equation}\label{eq:derivs}
\begin{split}
\partial_\tau^2\ell_\tau&=\sum_{l=1}^{p-1}\rho\Big(\frac{\fa_l(\cdot,\tau)}\sigma\Big)\frac{(\partial_\tau\fa_l)^2}{\sigma^2},\\[2mm]
\nabla\dot\ell_\tau&=\sum_{l=1}^{p-1}\rho\Big(\frac{\fa_l(\cdot,\tau)}\sigma\Big)\frac{(\partial_\tau\fa_l)\,\nabla\fa_l}{\sigma^2},\\[2mm]
\nabla^2\ell_\tau&=\sum_{l=1}^{p-1}\rho\Big(\frac{\fa_l(\cdot,\tau)}\sigma\Big)\frac{\nabla\fa_l\,\nabla\fa_l^\top}{\sigma^2}.
\end{split}
\end{equation}
Substituting these formulas into~\eqref{eq:Jgamma} and completing the square gives~\eqref{eq:Jdef}. The non-negativity of $\mathsf{J}_{\tau} h$ follows from that of $\rho$, i.e.~log-concavity of $\Phi$.
\end{proof}

The next lemma compares $\J_\tau$ with $\J_0$.

\begin{lemma}\label{lem:Jcomp}
Let $\fa$ be a tuple of affine functions on $\R^{p-1}\times\R$, and let $pn$ be sufficiently large.
If $\max_{l\le p-1}\norm{\nabla\fa_l}_2\le2$ and $\norm{\partial_\tau\fa}_\infty\le2\bar B/\sqrt k$, then for every $\tau\in[0,1]$, and every
$h\in C^2(\R^{p-1})$, the following inequality holds pointwise on $\R^{p-1}$,
\[
\J_\tau h\le2\J_0h+2\norm{\partial_\tau\fa}_\infty^2(pn)^{-4}.
\]
\end{lemma}

\begin{proof}
Because the functions~$\fa_l$ in the formula~\eqref{eq:Jdef} are affine, the derivatives $\partial_{\tau}\fa_l$ and $\nabla \fa_l$ do not depend on $\tau$. Thus, a comparison of $\J_{\tau}h$ and $\J_0 h$ reduces to comparing $\rho(\fa_l(z,\tau)/\sigma)$ and $\rho(\fa_l(z,0)/\sigma)$ at any fixed $z\in\R^{p-1}$.
In Lemma~\ref{lem:hazardfacts}(c) it is shown that
$$\rho(\tilde s)\le2\rho(s)+(pn)^{-5}$$
when
 $s,\tilde s\in\R$ satisfy $\abs{\tilde s-s}\le2(\log(pn))^{-3/2}$ and $pn$ is sufficiently large. Here, we will take $\tilde s=\fa_l(z,\tau)/\sigma$ and $ s=\fa_l(z,0)/\sigma$, and it can be checked that $|\tilde s-s|\leq \tau\|\partial_{\tau}\fa\|_{\infty}/\sigma\leq 2(\log(pn))^{-3/2}$, due to the assumption on $\fa$, and the basic relation~\eqref{eq:ratio}.
To proceed, we insert the bound on $\rho(\tilde s)$ into~\eqref{eq:Jdef} and compare with $2\J_0h(z)$, to obtain
\begin{equation}\label{eqn:Jtaualmost}
\J_\tau h(z)\le2\J_0h(z)-\norm{\nabla h(z)}_2^2+(pn)^{-5}\textstyle\sum_{l=1}^{p-1}v_l,
\end{equation}
where we put $v_l=(\partial_\tau\fa_l-\ip{\nabla\fa_l}{\nabla h(z)})^2/\sigma^2$. Furthermore, we have $v_l\leq (2\norm{\partial_\tau\fa}_\infty^2+8\norm{\nabla h(z)}_2^2)\,n$, because $\norm{\nabla\fa_l}_2\le2$ and $\sigma^{-2}\le n$, and so the bound~\eqref{eqn:Jtaualmost} leads to the stated result.
\end{proof}

\subsection{Changing the measure $\mu_{\tau}$ to $\mu_0$}\label{subsec:HSfamily}

This subsection establishes the bound~\eqref{eq:stepii} in the following lemma.

\begin{lemma}\label{lem:stab}
Let $\fa$ be a tuple of affine functions on $\R^{p-1}\times\R$, and let $pn$ be sufficiently large. Suppose that
$\max_{l\le p-1}\norm{\nabla\fa_l}_2\le2$, $\norm{\partial_\tau\fa}_\infty\le2\bar B/\sqrt k$ and $A(0)\ge-10\log(pn)$, and let $h^\star\in\mathscr H$ be
the function given by Lemma~\ref{lem:varrep}. Then, for every $\tau\in[0,1]$,
\[
\E_{\mu_\tau}\J_0h^\star\lesssim\E_{\mu_0}\J_0h^\star+\frac{\norm{\J_0h^\star}_{\sup}}{(pn)^{10}}.
\]
\end{lemma}

\begin{proof}

The proof is based on the formula
\begin{equation}\label{eq:densformula}
\frac{d\mu_\tau}{d\mu_0}=e^{\ell_0-\ell_\tau+A(0)-A(\tau)},
\end{equation}
which follows from the definition $d\mu_r=e^{-\ell_r-A(r)}d\gamma_{p-1}$ for $r\in\R$. We will show that the exponent in~\eqref{eq:densformula} is
at most a constant multiple of $\eps\ell_0+1$, where we put $\eps=1/\log(pn)$ throughout the current proof. Consequently, the density $d\mu_{\tau}/d\mu_0$ is bounded by a constant, except on a set where $\ell_0$ is
large, and the assumption  $A(0)\geq-10\log(pn)$ will turn out to make the $\mu_0$-probability of that set negligible. To bound the exponent in~\eqref{eq:densformula}, we will proceed in two stages by comparing $\ell_\tau$ with
$\ell_0$, and then comparing $A(\tau)$ with $A(0)$.

\medskip\noindent\emph{Comparison of $\ell_\tau$ and $\ell_0$.} Our aim in this stage is to show that
\begin{equation}\label{eq:lcomp}
\abs{\ell_\tau-\ell_0}\lesssim\eps\ell_0+1,
\end{equation}
and we will do this by developing a suitable bound on the derivative $\dot \ell_{r}=\partial_r \ell_r$ that holds for any affine tuple $\fa$. Differentiating the formula for $\ell_r$ in~\eqref{eq:elldef} gives
\begin{equation*}\dot\ell_r=-\sum_{l=1}^{p-1}\lambda(\fa_l(\cdot,r)/\sigma)\textstyle\frac{\partial_\tau\fa_l}{\sigma},
\end{equation*}
where we let
$\lambda=(\log\Phi)'=\phi/\Phi$. The function $\lambda$ is bounded pointwise in Lemma~\ref{lem:hazardfacts}(b) according to
\[
\lambda(s)\le\begin{cases}
-2\big(1+s_+\big)\log\Phi(s) & \text{for all } s\in\R,\\[1mm]
 \ \ 2\phi(s) & \text{for all } s\ge0,
\end{cases}
\]
where we write $s_+=\max(s,0)$.
To lighten notation, we will also put $s_l=\fa_l(z,r)/\sigma$ for fixed $z\in \R^{p-1}$ and $r\in\R$. Applying the first bound when $s_l$ is at most $\sqrt{2\log(pn)}$ and applying the second bound when $s_l$ exceeds $\sqrt{2\log(pn)}$, it follows that
\begin{equation*}\begin{split}
\sum_{l=1}^{p-1}\lambda(s_l)&\, \lesssim \, 1+ \sqrt{\log(pn)}\,\sum_{l=1}^{p-1}(-\log\Phi(s_l)) .
\end{split}
\end{equation*}
Noticing that the sum on the right is simply $\ell_r(z)$,
combining the last few steps shows that for all $r\in\R$ and $z\in\R^{p-1}$ we have
\begin{equation}\label{eq:slope}
\abs{\dot\ell_r(z)}\lesssim\frac{\norm{\partial_\tau\fa}_\infty}\sigma\Big(\sqrt{\log(pn)}\,\ell_r(z)+1\Big).
\end{equation}
We now divide both sides by $\ell_r+1$ to get $\abs{\partial_r\log(\ell_r+1)}\lesssim
\sqrt{\log(pn)}\,\norm{\partial_\tau\fa}_\infty/\sigma$ and then integrate the latter bound over
$r\in[0,\tau]$ to obtain
\begin{equation}\label{eqn:almostratio}
\begin{split}
\Big|\log\frac{\ell_\tau+1}{\ell_0+1}\Big|&\lesssim
\sqrt{\log(pn)}\,\frac{\tau \norm{\partial_\tau\fa}_\infty}{\sigma}\lesssim \eps
\end{split}
\end{equation}
where we have used the upper bound on $\|\partial_{\tau}\fa\|_{\infty}$ in the assumption, the basic relation~\eqref{eq:ratio}, and $\tau\leq 1$. Lastly, the inequality~\eqref{eqn:almostratio} leads to~\eqref{eq:lcomp} by letting $\log((\ell_{\tau}+1)/(\ell_0+1))$ play the role of $x$ in the bound $\abs{e^x-1}\lesssim \varepsilon e^{\mathcal{O}(\varepsilon)}$  that holds for all $|x|\lesssim \varepsilon$. This finishes the comparison of $\ell_{\tau}$ and $\ell_0$.

\medskip\noindent\emph{Comparison of $A(\tau)$ and $A(0)$.} 
Integrating the identity $e^{-\ell_{\tau}-A(0)}=e^{\ell_0-\ell_{\tau}}\frac{d\mu_0}{d\gamma_{p-1}}$ against
$\gamma_{p-1}$ leads to 
\begin{equation*}
    \begin{split}
    A(0)-A(\tau)=&-\log\E_{\mu_0}(e^{\ell_0-\ell_\tau})\\
    \leq & \ \E_{\mu_0}(\ell_\tau-\ell_0)\\
      \lesssim &   \ \varepsilon\E_{\mu_0}(\ell_0)+1
    \end{split}
\end{equation*}
where the second step uses Jensen's inequality and the third step uses~\eqref{eq:lcomp}. Since the KL-divergence between $\mu_0$ and $\gamma_{p-1}$ satisfies $0\leq \int \log(d\mu_0/d\gamma_{p-1})d\mu_0=\E_{\mu_0}(-\ell_0-A(0))$, we get
\begin{equation*}
    \begin{split}
    A(0)-A(\tau)&\lesssim -\varepsilon A(0)+1\\
    &\lesssim 1,
    \end{split}
    \end{equation*}
   which uses the assumed lower bound on $A(0)$ and the choice of $\varepsilon=1/\log(pn)$.
Applying this and~\eqref{eq:lcomp} to~\eqref{eq:densformula} shows that 
\begin{equation}\label{eq:dens}
\frac{d\mu_\tau}{d\mu_0}\;\lesssim \; e^{\mathcal{O}(\eps\ell_0)}.
\end{equation}

\medskip\noindent\emph{Combining.} As a temporary shorthand, let $\varphi=\J_0h^\star$. This function is nonnegative and has finite sup-norm due to Lemma~\ref{lem:Jformula} (as $\rho<1$ and $h^\star\in\mathscr H$). To conclude the proof, note that $\E_{\mu_\tau}\varphi = \E_{\mu_0}(\varphi \,d\mu_{\tau}/d\mu_0)$ and consider the set $\Omega=\{\ell_0\le20\log(pn)\}$. On the set $\Omega$, the density  bound~\eqref{eq:dens} implies $d\mu_{\tau}/d\mu_0\lesssim 1$.
Meanwhile, on the set $\Omega^c$, we will use the same density bound to write $(d\mu_{\tau}/d\mu_0)d\mu_0=\exp(\mathcal{O}(\varepsilon\ell_0)-\ell_0-A(0))d\gamma_{p-1}$. Hence,
\begin{equation}\label{eqn:keyexponent}
    \E_{\mu_\tau}\varphi \lesssim \E_{\mu_0}\varphi \, +\, \|\varphi\|_{\sup}\int_{\Omega^c} e^{\mathcal{O}(\varepsilon \ell_0)-\ell_0-A(0)}d\gamma_{p-1}.
\end{equation}
From the definition of $\Omega^c$ and the assumed bound $-A(0)\leq 10\log(pn)$, it follows that the exponent in~\eqref{eqn:keyexponent} is at most $-10\log(pn)+\mathcal{O}(1)$, which gives the stated result.
\end{proof}

\subsection{The function $h^{\star}$ and the proof of Lemma~\ref{lem:HSvar}}\label{subsec:HSvar}

Before proving Lemma~\ref{lem:HSvar}, we first prove Lemma~\ref{lem:varrep} below, which establishes the existence and regularity of the function $h^{\star}$ used in Section~\ref{subsec:HScomp}. Also, recall that $\mathscr{H}$ denotes the class of functions $h\in C^3(\R^{p-1})$ with bounded first and second derivatives.

\begin{lemma}\label{lem:varrep}
Let $\fa$ be a tuple of affine functions on $\R^{p-1}\times\R$ with $\max_{l\le p-1}\norm{\nabla\fa_l}_2\le2$. Then, there exists a function $h^\star\in\mathscr H$ such that $\Lop_0h^\star=\E_{\mu_0}\dot\ell_0-\dot\ell_0$ and
\begin{equation}\label{eq:J0bound}
\norm{\J_0h^\star}_{\sup}\lesssim(pn)^5\norm{\partial_\tau\fa}_\infty^2.
\end{equation}
\end{lemma}

\begin{proof}
It is shown in~\citet[Theorem 2.1]{MG16} that a solution $h^{\star}\in C^2(\R^{p-1})$ to the equation $\Lop_0h=\E_{\mu_0}\dot\ell_0-\dot\ell_0$ exists. Once a $C^2$ solution is known to exist, the differentiability can be upgraded to $C^{\infty}$, because the operator $\Lop_0=\Delta-\ip{\nabla W_0}{\nabla\cdot}$ has $C^{\infty}$ coefficients and the right-hand side $\E_{\mu_0}\dot\ell_0-\dot\ell_0$ is $C^{\infty}$~\citep[Theorem~6.17]{GT01}.

It remains to show that the first and second derivatives of $h^{\star}$ are bounded and that~\eqref{eq:J0bound} holds. After adapting the notation of~\citet[Theorem 2.1]{MG16} to our setting, it provides the following quantitative derivative bounds 
\begin{equation*}\norm{\nabla h^{\star}}_{\sup}\lesssim \norm{\nabla\dot\ell_0}_{\sup} \qquad 
    \text{ and } \qquad \norm{\nabla^2 h^{\star}}_{\sup}\lesssim  \norm{\nabla^3 \ell_0 }_{\sup}\norm{\nabla\dot\ell_0}_{\sup}+\norm{\nabla^2\dot\ell_0}_{\sup}.
\end{equation*}
From the formula~\eqref{eq:elldef}, the derivatives of $\dot \ell_0$ and $\ell_0$ can be computed explicitly, and since the functions $\rho,\rho', \rho''$ are bounded (Lemma~\ref{lem:hazardfacts}(a)), it can be checked that
\begin{equation*}\norm{\nabla\dot\ell_0}_{\sup}\le\frac{2p\norm{\partial_\tau\fa}_\infty}{\sigma^2},\qquad\norm{\nabla^2\dot\ell_0}_{\sup}\lesssim\frac{p\norm{\partial_\tau\fa}_\infty}{\sigma^3},\qquad
\norm{\nabla^3 \ell_0}_{\sup}\lesssim\frac p{\sigma^3}.
\end{equation*}
(Note that parts of~\eqref{eq:derivs} can be re-used here.) In particular, this shows that $\|\nabla h^{\star}\|_{\sup}$ and $\|\nabla^2 h^{\star}\|_{\sup}$ are finite, and hence $h^{\star}\in\mathscr{H}$. Finally, for the bound on $\|\J_0 h^{\star}\|_{\sup}$, applying the last two displays to~\eqref{eq:Jdef}, and using the basic inequalities $\rho<1$ and
$\sigma^{-2}\le n$, 
leads to~\eqref{eq:J0bound}.
\end{proof}

\begin{proof}[Proof of Lemma~\ref{lem:HSvar}] 
The proof proceeds in a manner that is similar to~\citet[Section~2]{BC13}.
We begin by considering $h\in C^3(\R^{p-1})$ with compact support, and then we will extend the argument later to handle $h\in\mathscr{H}$. As notation that is specific to this proof, let $\varphi\in C^1(\R^{p-1})$ be a test function, and let $g=\dot\ell_\tau-\E_{\mu_\tau}\dot\ell_\tau$ so that $\Var_{\mu_{\tau}}(\dot\ell_{\tau})=\E_{\mu_{\tau}} g^2.$
The basic idea is to expand the square $ \E_{\mu_{\tau}}(g+\Lop_{\tau} h)^2$ to obtain $\E_{\mu_{\tau}}g^2= -2\E_{\mu_{\tau}}(g \Lop_{\tau} h)-\E_{\mu_{\tau}}(\Lop_{\tau} h)^2+\E_{\mu_\tau}(g+\Lop_{\tau}h)^2$ and then show that the right side has the form of the stated result in~\eqref{eqn:Helffer}.

By the divergence theorem and the compact support of $h$, the divergence of $\varphi e^{-W_{\tau}}\nabla h$ integrates to 0 over $\R^{p-1}$. Rearranging that equation and noting that the density of $\mu_{\tau}$ with respect to Lebesgue measure on $\R^{p-1}$ is proportional to $e^{-W_{\tau}}$, it can be checked that
\begin{equation}\label{eq:HSibp}
\E_{\mu_\tau}\ip{\nabla\varphi}{\nabla h}=-\E_{\mu_\tau}(\varphi\,\Lop_\tau h).
\end{equation}
We will use this identity three times with certain choices of $\varphi$ and $h$. First, we take $\varphi=\Lop_\tau h$, which lies in $C^1(\R^{p-1})$ because $h\in C^3(\R^{p-1})$, and
obtain $\E_{\mu_\tau}(\Lop_\tau h)^2=-\E_{\mu_\tau}\ip{\nabla\Lop_\tau h}{\nabla h}$. Using $\partial_j\Lop_\tau h=\Lop_\tau\partial_jh-(\nabla^2W_\tau\nabla h)_j$, this becomes
\[
\E_{\mu_\tau}(\Lop_\tau h)^2=-\sum_{j=1}^{p-1}\E_{\mu_\tau}\big(\partial_jh\,\Lop_\tau\partial_jh\big)+\E_{\mu_\tau}\ip{\nabla h}{\nabla^2W_\tau\nabla h}.
\]
Second, we apply~\eqref{eq:HSibp} with $(\varphi,h)$ replaced by $(\partial_jh,\partial_jh)$. This gives $-\E_{\mu_\tau}(\partial_jh\,\Lop_\tau\partial_jh)=\E_{\mu_\tau}\norm{\nabla\partial_jh}_2^2$, and summing over $j$ yields
\[
\E_{\mu_\tau}(\Lop_\tau h)^2=\E_{\mu_\tau}\Big(\|\nabla^2 h\|_F^2+\langle \nabla h,\nabla^2 W_{\tau} \nabla h\rangle\Big).
\]
Third, taking $\varphi=g$ in~\eqref{eq:HSibp} gives $-\E_{\mu_\tau}(g\,\Lop_\tau h)=\E_{\mu_\tau}\ip{\nabla g}{\nabla h}=\E_{\mu_\tau}\ip{\nabla\dot\ell_\tau}{\nabla h}$.
Inserting the last two facts into the expression for $\E_{\mu_\tau}g^2$ at the beginning shows that it has the form of~\eqref{eqn:Helffer}.

Thus, it remains to extend to the case when $h\in\mathscr{H}$. For such $h$, let $f_r=\chi(\cdot/r)h$, where $r>0$ is a cutoff parameter and $\chi$ is a smooth radial cutoff function on $\R^{p-1}$ with
$\1_{\mathbb B(0,1)}\le\chi\le\1_{\mathbb B(0,2)}$, and $\mathbb B(0,r)=\{z\in\R^{p-1}:\|z\|_2\leq r\}$.
Since $f_r$ is $C^3$ with compact support, the first part of the proof shows that~\eqref{eqn:Helffer} holds with $f_r$ in place of $h$. We now
let $r\to\infty$. Since $f_r=h$ on $\mathbb B(0,r)$, the functions $f_r$, $\nabla f_r$, $\nabla^2f_r$ and $\Lop_\tau f_r$ converge pointwise to $h$,
$\nabla h$, $\nabla^2h$ and $\Lop_\tau h$ as $r\to\infty$. To justify passing the limit under the expectation in~\eqref{eqn:Helffer}, we use dominated convergence,
based on the following observations.
\begin{itemize}
\item[(i)] The first two derivatives of $f_r$ are bounded uniformly in $r\ge1$. Indeed, $\nabla h$ and $\nabla^2h$ are bounded, $h$ grows at most
linearly, and the $k$th derivative of $\chi(\cdot/r)$ is bounded by a constant multiple of $r^{-k}$ and vanishes outside $\mathbb B(0,2r)$.
\item[(ii)] The functions $\|\nabla W_\tau\|_2$, $g$ and $\Lop_\tau f_r$ grow at most linearly in $\norm z_2$, uniformly in $r\ge1$. For $\nabla W_\tau$, this
follows from $\nabla W_\tau=z-\sum_{l\leq p-1}\lambda(\fa_l(\cdot,\tau)/\sigma)\,\nabla\fa_l/\sigma$ and $\lambda(s)\le1+\abs s$ (Lemma~\ref{lem:hazardfacts}(a)). For $g$,
it follows from the boundedness of $\nabla\dot\ell_\tau$, which holds by~\eqref{eq:derivs} and $\rho<1$. For $\Lop_\tau f_r=\Delta f_r-\ip{\nabla W_\tau}{\nabla f_r}$,
it follows from~(i).
\item[(iii)] The Hessian of $W_{\tau}$ satisfies $\|\nabla^2 W_{\tau}\|_{\sup}\leq 1+\|\nabla^2\ell_{\tau}\|_{\sup}<\infty$, due to~\eqref{eq:derivs} and $\rho<1$.
\item[(iv)] Since $d\mu_\tau\le e^{-A(\tau)}d\gamma_{p-1}$, every function of $z$ that grows at most polynomially is $\mu_\tau$-integrable.
\end{itemize}
By (i)--(iii), each integrand in~\eqref{eqn:Helffer} for $f_r$ is bounded by a polynomial in $\norm z_2$ that does not depend on $r\ge1$, and so
(iv) allows us to let $r\to\infty$. This proves~\eqref{eqn:Helffer} for all $h\in\mathscr H$.
\end{proof}

\subsection{Proofs of Lemmas~\ref{lem:sections} and~\ref{lem:hazard}}\label{subsec:sectionsproof}

\begin{proof}[Proof of Lemma~\ref{lem:sections}]
Let $s=\sqrt{u+\sigma^2}$ and observe that Lemma~\ref{lem:facetfn} gives
\[
\log\psi_j(y+\tau w)=\log q_j(y+\tau w)-\tfrac12a_j(y+\tau w)^2-\tfrac12\log(2\pi s^2),
\]
where $a_j(y+\tau w)$ is affine in $\tau$ with slope $-w_j/s$. It therefore suffices to find a tuple of affine functions satisfying the two stated
bounds such that
\begin{equation}\label{eq:Aq}
\log q_j(y+\tau w)=A(\tau)\qquad\text{for all }\tau\in\R,
\end{equation}
since differentiating twice gives the curvature identity~\eqref{eqn:curvatureidentity}, and setting $\tau=0$ gives $A(0)=\log q_j(y)$.

To proceed, let $Z\sim\gamma_p$ and let $\Pi_j$ denote the orthogonal projection onto $E_j$ so that  $\Pi_jZ\sim\gamma_{E_j}$. Consequently, the definitions of $q_j$ and $K_j$ give 
\begin{equation}\label{eq:chain}
    \begin{split}
        q_j(y)&=\PP(\Pi_j Z\in K_j)\\[0.2cm]
        &=\PP\Big(a_j(y)n_j+\Pi_jZ\in P(y)\Big)\\[0.2cm]
&=\PP\Big(\Sigma_u^{1/2}\big(a_j(y)n_j+\Pi_jZ\big)\leq \beta-y\Big),
    \end{split}
\end{equation}
where we recall $n_j=\frac{1}{s}\Sigma_u^{1/2}e_j$. Next, observe that the orthogonal projection $\Pi_j$ can be expressed as $\Pi_j=I_p-n_jn_j^{\top}$, and define the deterministic vector $v_j=\frac{1}{s}\Sigma_u^{1/2}n_j$, as well as the Gaussian vector $\breve G=\Sigma_u^{1/2}Z$. Thus, the last expression in~\eqref{eq:chain} becomes
\begin{equation}\label{eq:afterchain}
    q_j(y)=\PP\Big( \big[s a_j(y)  -\breve G_j\big] v_j+ \breve G\leq \beta-y\Big).
\end{equation}
Writing $v_{jl}$ for the $l$th entry of $v_j$, it can be checked that $v_{jj}=1$ and $\abs{v_{jl}}\le1$. Since $s\,a_j(y)=\beta_j-y_j$, the $j$th inequality
in~\eqref{eq:afterchain} holds trivially, and so
\[
q_j(y)=\PP\Big(\breve G_l-v_{jl}\breve G_j\le b_l\ \text{for all }l\ne j\Big),\qquad \text{where} \ b_l=\beta_l-y_l-v_{jl}(\beta_j-y_j).
\]
Next, recall that $G\sim \mathcal{N}(0,\Sigma)$ and $N\sim\mathcal{N}(0,I_p)$ denote independent Gaussian vectors, which implies $\breve G\overset d=\sqrt uG+\sigma N$. Based on this, each $\breve G_l-v_{jl}\breve G_j$ may be written in the form $\zeta_l+\sigma N_l$, where we put $\zeta_l=\sqrt uG_l-v_{jl}(\sqrt uG_j+\sigma N_j)$.
Importantly, the vector $\zeta=(\zeta_l)_{l\ne j}$ only depends on $N$ through its $j$th entry, and so $\zeta$ is independent of $(N_l)_{l\ne j}$. Taking advantage of this independence, as well as the fact that the entries of $N$ are i.i.d.~standard Gaussian, we have
\begin{equation}\label{eqn:Gintegral}
\begin{split}
    q_j(y)&=\PP\Big(N_l\leq \textstyle \frac{b_l-\zeta_l}{\sigma} \text{ for } l\neq j\Big)\\
    &=\E\prod_{l:l\neq j}^p\Phi\Big( \frac{b_l-\zeta_l}{\sigma}\Big)\\
    &=\int_{\R^{p-1}}\prod_{l:l\neq j}^p\Phi\Big(\frac{b_l-\ip{\eta_l}{z}}{\sigma}\Big)\,d\gamma_{p-1}(z)
    \end{split}
\end{equation}
where we have used $\zeta_l \overset d= \langle \eta_l, N'\rangle$, with $N'\sim \gamma_{p-1}$ and the vectors $(\eta_l)_{l\neq j}$ corresponding to the rows of 
$(\E(\zeta\zeta^\top))^{1/2}$.
To conclude the proof, we replace $y$ with $y+\tau w$ and define the affine functions
\[
\fa_l(z,\tau)=b_l-\ip{\eta_l}z-\tau(w_l-v_{jl}w_j),\qquad l\ne j.
\]
In this notation, the integral in~\eqref{eqn:Gintegral} matches the form of the one that defines $A(\tau)$ in equation~\eqref{eq:Adef}, provided that the affine functions are relabelled as $\fa_1,\dots,\fa_{p-1}$.
Finally, to verify that the affine functions satisfy the required conditions, observe that $\abs{\partial_\tau\fa_l}=\abs{w_l-v_{jl}w_j}\le2\norm w_\infty$. Furthermore, we have
$\norm{\nabla\fa_l}_2^2=\norm{\eta_l}_2^2=\Var(\zeta_l)$, and a direct computation using $v_{jl}=\frac{u}{s^2}\Sigma_{jl}$ and $u\leq 2$ shows that $\Var(\zeta_l)=u-s^2v_{jl}^2\le2$.
\end{proof}

\begin{proof}[Proof of Lemma~\ref{lem:hazard}]
Let $\fa$ be the tuple of affine functions given by Lemma~\ref{lem:sections} with associated log-partition function $A$. Also, let $u$, $y$, $w$ and $j$ be fixed as in the stated result.
By~\eqref{eq:Aq}, we have $\log q_j(y+\tau w)=A(\tau)$ for all $\tau\in\R$, and differentiating at $\tau=0$ gives
\[
\partial_w\log q_j(y)=A'(0).
\]
So, it suffices to bound $\abs{A'(0)}$. Differentiating the identity $e^{A(\tau)}=\int e^{-\ell_\tau}\,d\gamma_{p-1}$ gives $A'(0)=-\E_{\mu_0}\dot\ell_0$.
The derivative bound~\eqref{eq:slope} in the proof of Lemma~\ref{lem:stab} holds for an arbitrary tuple $\fa$, and so
\[
\abs{A'(0)}\le\E_{\mu_0}\abs{\dot\ell_0}\lesssim\frac{\norm{\partial_\tau\fa}_\infty}\sigma\Big(\sqrt{\log(pn)}\,\E_{\mu_0}\ell_0+1\Big).
\]
We will now express the right side in terms of $w$ and $q_j(y)$. By Lemma~\ref{lem:sections}, the tuple $\fa$ satisfies
$\norm{\partial_\tau\fa}_\infty\le2\norm w_\infty$. Next, a simple calculation shows that the KL-divergence between $\mu_0$ and $\gamma_{p-1}$ is equal to $\E_{\mu_0}(-\ell_0-A(0))$, and so its non-negativity gives 
$\E_{\mu_0}\ell_0\le-A(0)$. Lastly, Lemma~\ref{lem:sections} gives $-A(0)=\log(1/q_j(y))$.
Substituting these bounds into the last display completes the proof.
\end{proof}

\section{Proof of Proposition~\ref{thm:FNB}: facet sensitivity bounds}\label{sec:FNB}
For proving parts (a) and (b) of Proposition~\ref{thm:FNB}, recall that $P$ is a polyhedron in $\R^p$ with $p$ linearly independent normal vectors as in~\eqref{eqn:Polygeneral}, and that Lemma~\ref{lem:sensitivityformula} gives the following explicit formula for the facet sensitivity,
\begin{equation}\label{eq:Lambdaexplicit}
\Fc(P)=\sum_{j=1}^p\big(1+a_j^2+\norm{m_j}_2^2+\tr(\Pi_j-\mathfrak S_j)\big)\nu_j(P).
\end{equation}

\medskip\noindent\emph{Part (a).} By~\eqref{eq:Lambdaexplicit}, it suffices to show that
\begin{equation*}\sum_{j=1}^p\big(1+a_j^2+\norm{m_j}_2^2\big)\nu_j(P)\lesssim(\log(2p))^{3/2}\ \ \ \text{and}\ \ \ \sum_{j=1}^p\tr(\Pi_j-\mathfrak S_j)\,\nu_j(P)\lesssim(\log(2p))^{3/2}.
\end{equation*}
These two bounds have been established by \citet[Lemmas~26 and~30]{GZ26} for polyhedra that satisfy a few basic conditions, and so it is enough to check that they are satisfied by $P$. (For the purpose of referencing the cited results, the notation in~\cite{GZ26} corresponds to ours via $m_P=p$, $s_j\leftrightarrow\nu_j(P)$, $K_j\leftrightarrow \mathfrak{S}_j$, $\mu_j\leftrightarrow m_j$, $\Delta_j\leftrightarrow\tr(\Pi_j-\mathfrak S_j)$, and $Q_j\leftrightarrow K_j$.)

The conditions for $P$ that we must check are (1) that it is a proper subset of $\R^p$, (2) that it is non-empty, (3) that it has non-empty interior, and (4) that none of the constraint inequalities in~\eqref{eqn:Polygeneral} are redundant (i.e.~none may be removed without changing $P$). All four conditions can be checked using the observation that for any $(b_1,\dots,b_p)\in\R^p$, there is some $z\in\R^p$ such that
\begin{equation}\label{eqn:linearlyindep}
    \langle n_j,z\rangle = b_j \qquad \text{  for all $j\leq p$},
\end{equation}
because the unit normal vectors $n_1,\dots,n_p$ are linearly independent. The first two conditions are trivial and so we only give details for the latter two. For condition (3), let $z$ correspond to the choice $b_j=a_j-1$ for all $j\leq p$, and notice that for any $u\in\R^p$ with $\|u\|_2\leq 1$, we have $\langle n_j,z+u\rangle \leq a_j-1+\|n_j\|_2\|u\|_2\leq a_j$ for all $j$. In other words, the unit $\ell_2$-ball centered at $z$ is contained in $P$, and so $P$ has non-empty interior. For condition (4), fix any $j_0\in\{1,\dots,p\}$ and let $z'$ be the vector that satisfies $\langle n_j,z'\rangle=a_j$ for all $j\neq j_0$, as well as $\langle n_{j_0},z'\rangle=a_{j_0}+1$. Hence, we have $z'\not\in P$, but if the constraint inequality corresponding to $j_0$ were removed from the definition of $P$, then $z'$ would be included in $P$, and thus condition (4) is satisfied.

\medskip\noindent\emph{Part (b).} We first bound $\Fc(P)$ by two simpler sums. Recall that $\nu_j(P)=\phi(a_j)q_j$ with $0<q_j\le1$
(Lemma~\ref{lem:structure}). Also note that $\tr(\Pi_j-\mathfrak S_j)\leq \tr(\Pi_j)= p-1$ and $\norm{m_j}_2^2\le2\log(1/q_j)$, where the latter bound follows from letting the subspace $E_j$ play the role of Euclidean space in Lemma~\ref{lem:condgauss}(a).
Applying these facts to~\eqref{eq:Lambdaexplicit} leads to
\begin{equation}\label{eq:Bsplit}
\Fc(P)\lesssim\sum_{j=1}^p\big(p+a_j^2\big)\nu_j(P)+\sum_{j=1}^p\phi(a_j)\,q_j\log(1/q_j).
\end{equation}
The main step is to show that every facet has small Gaussian surface measure, in the sense that
\begin{equation}\label{eq:nusmall}
\nu_j(P)\lesssim e^{-r^2/8}\qquad\text{for all }j\le p.
\end{equation}
To prove this, we will consider an arbitrary fixed $j$ and then separately address the two cases when $|a_j|\geq r/2$ and $|a_j|<r/2$. In the first case, the bound~\eqref{eq:nusmall} holds because 
$$\nu_j(P)=\phi(a_j)q_j\leq \phi(a_j)\le\phi(r/2)\lesssim e^{-r^2/8}.$$
The second case  $\abs{a_j}<r/2$ is more involved. Let $l^{\star}$ denote the minimizing index satisfying $a_{l^{\star}}=\min_{l\leq p} a_l$ and let $W\sim\gamma_{E_j}$.  Observe that the event $W\in K_j$ is the same as $(a_j n_j+W) \in F_j$, which implies that the point $(a_j n_j+W)$ is in $P$, and hence satisfies the $l^{\star}$th inequality in $P$. Thus, $\{W\in K_j\}\subset \{\langle n_{l^{\star}}, a_j n_j +W\rangle \leq a_{l^{\star}}\}$. Furthermore, because we are in the case that $\abs{a_j}<r/2$, the assumption $|\min_{l\leq p}a_l|\geq r$ implies $a_{l^{\star}}\leq -r$, and so the last few observations give
\[
q_j=\gamma_{E_j}(K_j)=\PP(W\in K_j)\le\PP\big(\ip{n_{l^{\star}}}{a_jn_j+W}\le-r\big).
\]
Furthermore, note that $\ip{n_{l^{\star}}}{a_jn_j+W}$ is Gaussian with mean $a_j\ip{n_{l^{\star}}}{n_j}\ge-r/2$ and variance
$\norm{\Pi_jn_{l^{\star}}}_2^2\le1$, which means that the last probability is at most $1-\Phi(r/2)\leq e^{-r^2/8}$. Combining this with $\nu_j(P)=\phi(a_j)q_j\lesssim q_j$ proves~\eqref{eq:nusmall}.

It remains to bound the two sums in~\eqref{eq:Bsplit}. For the terms in the first sum, we take the geometric mean of the upper bounds~\eqref{eq:nusmall} and $\nu_j(P)\le\phi(a_j)$ to get
\begin{equation}\label{eqn:firstterms}
    \begin{split}
        (p+a_j^2)\nu_j(P)&\lesssim (p+a_j^2)\sqrt{\phi(a_j)e^{-r^2/8}}\\
        &\lesssim \Big(p +\sup_{a\in \R}(a^2 e^{-a^2/4})\Big)e^{-r^2/16}\\
        &\lesssim p e^{-r^2/16} .
    \end{split}
\end{equation}
For the terms in the second sum in~\eqref{eq:Bsplit}, the elementary inequality $q\log(1/q)\le\sqrt{q}$ for $q\in(0,1]$, and the formula $\sqrt{q_j}=\sqrt{\nu_j(P)/\phi(a_j)}$ give
\begin{equation}\label{eqn:secondterms}
\phi(a_j)\,q_j\log(1/q_j)\le\sqrt{\phi(a_j)\,\nu_j(P)}\lesssim e^{-r^2/16},
\end{equation}
where we have used~\eqref{eq:nusmall} in the second step.
Summing the bounds in~\eqref{eqn:firstterms} and~\eqref{eqn:secondterms} over $j=1,\dots,p$ produces the stated result for part (b).

\medskip\noindent\emph{Part (c).} Let $r>0$ and consider multiplying $\Fc(y)$ by the indicator  $\1\{|\min_{j\leq p}a_j(y)|\geq r\}$ and its complementary version. Since the bound in part (a) holds uniformly with respect to all choices of $P(y)$, and since we can apply the bound in part (b) when the condition $|\min_{j\leq p}(a_j(y))|\geq r$ holds, it follows that
\begin{equation*}
\begin{split}
    \Fc(y) &\lesssim (\log(pn))^{3/2}\1\Big\{\Big|\min_{j\leq p}a_j(y)\Big|< r\Big\} \, + \, p^2 e^{-r^2/16}.
\end{split}
\end{equation*}
Taking $r=8\sqrt{\log(pn)}$ and noting that
$\min_{j\leq p}a_j(y)=\frac{1}{\sqrt{u+\sigma^2}}\min_{j\leq p}(\beta_j-y_j)$ leads to the stated result for part (c).\qed

\section{Proofs of Lemmas~\ref{lem:interp} and~\ref{lem:segment}}\label{sec:remaining}

\begin{proof}[Proof of Lemma~\ref{lem:interp}]
As a temporary shorthand, let $f(y)=\PP(y+\sigma N\le\beta)=\prod_j\Phi((\beta_j-y_j)/\sigma)$, so that the smoothed orthant indicator is
$g_u(y)=\E f(y+\sqrt uG)$.
Differentiating in $u$ and applying Stein's identity~\citep[eqn.~(A.17)]{Tal11}   
\begin{equation}\label{eqn:steinidentity}
\E(G_j\varphi(G))=\sum_{l=1}^p\Sigma_{jl}\E\,\partial_l\varphi(G)
\end{equation}
to the function
$\varphi(G)=\partial_jf(y+\sqrt u G)$, it follows that
\begin{equation*}
\partial_ug_u(y)
=\textstyle\frac12\LapS g_u(y).
\end{equation*}
In turn, by the chain rule and recalling $u(t)=1-t+\sigma^2$, we have
\begin{equation}\label{eqn:almost}
\begin{split}
\frac d{dt}\,g_{u(t)}(\sqrt tS)&=\frac1{2\sqrt t}\ip{\nabla g_{u(t)}(\sqrt tS)}{S}-\frac12\LapS g_{u(t)}(\sqrt tS).
\end{split}
\end{equation}
Finally, note that the first term on the right side of~\eqref{eqn:almost} matches the sum on the right side of Lemma~\ref{lem:interp} because
$\E\ip{\nabla g_u(\sqrt tS_{\setminus i})}{X_i/\sqrt k}=0$ for each $i$. 
\end{proof}

\begin{proof}[Proof of Lemma~\ref{lem:segment}]
 Both parts rely on a
bound for $\partial_w\log q_j(x)$, which we will derive first for any fixed $j\leq p$. By~\eqref{eq:ratio}, we have $\norm w_\infty/\sigma\le C_0^{-1}(\log(pn))^{-3/2}$, and so
Lemma~\ref{lem:hazard} gives
\begin{equation}\label{eq:s7-hazard}
\abs{\partial_w\log q_j(x)}\le\frac C{C_0}\Big(\frac{\log(1/q_j(x))}{\log(pn)}+(\log(pn))^{-3/2}\Big),
\end{equation}
where $C$ is the absolute constant implicit in Lemma~\ref{lem:hazard}.

\medskip\noindent\emph{Part (a).} Here, $j$ is assumed to be regular at $x$, and so if we let $r=\sqrt{20\log(pn)}$, then $\abs{a_j(x)}\le r$ and $q_j(x)\ge(pn)^{-10}$.
 We bound the two terms of $\partial_w\log\psi_j(x)$ in~\eqref{eq:s7-Dlogpsi}
separately, which gives
\[
\abs{\partial_w\log\psi_j(x)}\le\frac{\norm w_\infty\abs{a_j(x)}}{\sqrt{u+\sigma^2}}+\abs{\partial_w\log q_j(x)}\le\frac{\sqrt{20}}{C_0\log(pn)}+\frac{11C}{C_0},
\]
where we used 
$\log(1/q_j(x))\le10\log(pn)$ from~\eqref{eq:s7-hazard}. Since $C_0$ may be taken to be large, the right side is at most $\frac1{100}$.

\medskip\noindent\emph{Part (b).} Now suppose that $j$ is not regular at $x$, and let $\tau\in[0,1]$. By the non-regularity, we may assume that either $\abs{a_j(x)}>r$ holds, or that the two conditions
$\abs{a_j(x)}\le r$ and $q_j(x)<(pn)^{-10}$ hold. We will treat these two cases separately, but in both cases we will use $\nu_j=\phi(a_j)q_j$ with $q_j\le1$
(Lemma~\ref{lem:structure}).

In the first case $\abs{a_j(x)}>r$, the definition of the offset $a_j(y)=(\beta_j-y_j)/\sqrt{u+\sigma^2}$ implies
\[
\abs{a_j(x+\tau w)-a_j(x)}=\frac{\tau\abs{w_j}}{\sqrt{u+\sigma^2}}\le\frac{\norm w_\infty}\sigma\le1.
\]
Hence $\abs{a_j(x+\tau w)}\ge r-1$, and so $\nu_j(x+\tau w)\le\phi(r-1)\le(pn)^{-5}$. For the
other bound in~(b), we use $q\log(1/q)\le1$ for $q\in(0,1]$, together with the fact that $(a+2)\phi(a)$ is decreasing on $[1,\infty)$, to obtain
\begin{equation}
    \begin{split}
        \Big(\abs{a_j(x)}+\textstyle\log\frac1{q_j(x)}+1\Big)\nu_j(x)&\le\big(\abs{a_j(x)}+2\big)\,\phi(a_j(x))\\
        &\le(r+2)\,\phi(r)\\
        &\lesssim(pn)^{-5}.
    \end{split}
\end{equation}

Now we turn to the second case. Consider the function
$\vartheta(\tau)=\log(1/q_j(x+\tau w))$, which is convex because $q_j$ is log-concave (Lemma~\ref{lem:facetfn}). Its slope at the origin is
$\vartheta'(0)=-\partial_w\log q_j(x)$, and so~\eqref{eq:s7-hazard} gives $\vartheta'(0)\ge-\frac14\vartheta(0)-1$, as $C_0$ is large. (Here, the constants $-\frac{1}{4}$ and $-1$ have been used for simplicity.)  Thus,
\begin{equation}
    \begin{split}
        \vartheta(\tau)
        %&
        %\geq \vartheta(0)+\vartheta'(0) \tau\\
        &\ge\vartheta(0)-\tau\big(\textstyle\frac14\vartheta(0)+1\big)\\
        &\ge\frac34\vartheta(0)-1\\
        &\ge5\log(pn),
    \end{split}
\end{equation}
where $\vartheta(0)>10\log(pn)$ holds because we are considering the second case.
Rewriting the previous display as an upper bound on $q_j(x+\tau w)$ then leads to $\nu_j(x+\tau w)\le q_j(x+\tau w)\le(pn)^{-5}$. For the other bound in~(b), we use $\phi\le1$,
and the fact that $q\mapsto (r+1+\log(1/q))q$ is increasing on $(0,1]$ to obtain
\[
\begin{split}
\textstyle \Big(\abs{a_j(x)}+\log\frac1{q_j(x)}+1\Big)\nu_j(x)&\le\Big(\textstyle r+1+\log\frac1{q_j(x)}\Big)q_j(x)\\[1mm]
&\le\big(r+1+10\log(pn)\big)(pn)^{-10}\\[1mm]
&\lesssim(pn)^{-5}.
\end{split}
\]
\end{proof}

\section{Background results}\label{sec:tools}

This section collects the reduction lemma used in Section~\ref{sec:reductions}, results from the literature, a standard smoothing lemma, a lemma on the function $\phi/\Phi$, and elementary facts on the Gaussian
measures of convex sets and of the polyhedra of Section~\ref{subsec:sensitivity}.

\begin{lemma}\label{lem:reduction}
Theorem~\ref{thm:bounded} implies Theorem~\ref{thm:main}.
\end{lemma}
\begin{proof}
First, the supremum in Theorem~\ref{thm:main} can be replaced with a supremum over compact rectangles $A=\prod_j[\alpha_j,\beta_j]$ with
$\alpha_j\le\beta_j$. This is due to continuity of measure and the fact that taking a limit of compact rectangles approaching a general rectangle
will not affect a bound that is uniform over compact rectangles.

For a compact rectangle $A=\prod_j[\alpha_j,\beta_j]$ and a random vector $U$, observe that the event $\{U\in A\}$ is equivalent to
$\{(U,-U)\le(\beta,-\alpha)\}$, where $\alpha=(\alpha_1,\dots,\alpha_p)$ and $\beta=(\beta_1,\dots,\beta_p)$. Furthermore, the vectors
$(X_1,-X_1),\dots,(X_n,-X_n)$ are i.i.d.\ and centered in $\R^{2p}$, and their coordinates have unit variances and $\psi_2$-norms at most $B$. Also,
the Gaussian counterpart of these vectors is $(G,-G)$, and so it suffices to work with orthants in dimension $2p$. (However, for simplicity we
continue to write the dimension as $p$ in the rest of this proof.)

It remains to show there is an absolute constant $C>0$ and centered i.i.d.~random vectors $Y_1,\dots,Y_n\in\R^p$ with the following properties. These vectors must satisfy $\Var(Y_{1j})=1$ and  $|Y_{1j}|\leq \bar B\lesssim B\sqrt{\log(pn)}$  almost surely for all $j\leq p$, and also, if $H$ is a centered Gaussian vector with covariance matrix $\E(Y_1 Y_1^{\top})$, then the bound below must hold for any $\beta\in\R^p$, 
\begin{equation*}
\Big|\PP\Big(\frac1{\sqrt n}\sum_{i=1}^nX_i\le\beta\Big)-\PP(G\le\beta)\Big|\le\sup_{\beta'\in\R^p}\Big|\PP\Big(\frac1{\sqrt n}\sum_{i=1}^nY_i\le\beta'\Big)-\PP(H\le\beta')\Big|+\frac Cn.
\end{equation*}

We now construct each $Y_i$ by truncation. We may assume without loss of generality that $B^2\le n$, because otherwise Theorem~\ref{thm:main} holds trivially with $C\ge1$. 
Also note that  $1=(\E X_{1j}^2)^{1/2}\leq (\E X_{1j}^4)^{1/4}\le B$, and that if we define $\Omega_i=\{\norm{X_i}_\infty\le  B\sqrt{10\log(pn)}\}$, then each $\Omega_i$ fails with probability at most $2(pn)^{-9}$.
To proceed, define the mean and variance parameters  $\tilde\mu=\E(X_1\1_{\Omega_1})$ and $\tilde\sigma_j^2=\Var(X_{1j}\1_{\Omega_1})$, and put
$$Y_{ij}=\textstyle\frac{1}{\tilde\sigma_j}(X_{ij}\1_{\Omega_i}-\tilde\mu_j),$$ 
so that $Y_1,\dots,Y_n$ are centered and i.i.d.~with standardized entries.

Next, observe that the matrix
\[
\Delta=\Sigma-\Cov(X_1\1_{\Omega_1})=\E\big(X_1X_1^\top\1_{\Omega_1^c}\big)+\tilde\mu\tilde\mu^\top
\]
is positive semidefinite, and the Cauchy--Schwarz inequality implies $\norm{\tilde\mu}_\infty\le(pn)^{-4}$, as well as
$\Delta_{jj}\leq 3(pn)^{-3}$. Since the diagonal entries of $\E(X_1X_1^{\top})$ are all equal to $1$, it follows that $\tilde\sigma_j^2=1-\Delta_{jj}\ge\frac14$ for all $j\leq p$. Consequently, if we take $\bar B=10B\sqrt{\log(pn)}$, then we have the almost-sure bound $\|Y_1\|_{\infty}\le\bar B$. Moreover, on the event $\bigcap_{i\leq n}\Omega_i$, we have
\begin{equation}\label{eqn:nongaussian}
\Big\{\frac1{\sqrt n}\sum_{i=1}^nX_i\le\beta\Big\}=\Big\{\frac1{\sqrt n}\sum_{i=1}^nY_i\le\beta'\Big\},
\end{equation}
where we put $\beta'_j=\frac{1}{\tilde\sigma_j}(\beta_j-\sqrt n\,\tilde\mu_j)$ for $j\leq p$.

Due to~\eqref{eqn:nongaussian}, the proof will be complete if we can show $|\PP(H\leq \beta')-\PP(G\leq \beta)|\lesssim 1/n$. Since $\Delta\succeq0$, we may represent $G$ in distribution as
$G_j=\tilde\sigma_jH_j+D_j$, where $D\sim\mathcal N(0,\Delta)$ is independent of $H$. Thus, from the definition of $\beta'$,  it is enough to show
\begin{equation}\label{eqn:almostnazarov}
\Big|\PP\Big(H\leq \textstyle\frac{1}{\tilde\sigma}(\beta-\sqrt n \tilde\mu)\Big)-\PP\Big(H\leq \textstyle\frac{1}{\tilde\sigma}(\beta-D)\Big)\Big|\lesssim \frac{1}{n},
\end{equation}
where division by $\tilde\sigma$ is understood in an entrywise sense. By construction, the coordinates of $H$ satisfy $\Var(H_j)=1$ for $j\leq p$, and it can be checked that $\|\frac{1}{\tilde\sigma}(D-\sqrt n\tilde\mu)\|_{\infty}> 4(pn)^{-1}$ occurs with probability at most $\mathcal{O}(1/n)$.
Consequently, the proof is completed by Nazarov's inequality (Lemma~\ref{fact:nazarov}). 
\end{proof}

The following result is known as Nazarov's inequality~\citep{Naz03,CCKnaz}.

\begin{lemma}\label{fact:nazarov}
If $G$ is a Gaussian vector in $\R^d$ with $\min_{l\leq d}\Var(G_l)\ge1$, then for all $y\in\R^d$ and $\delta>0$ we have
$$\PP(G\le y+\delta\1)-\PP(G\le y)\le\delta(\sqrt{2\log d}+2).$$ 
\end{lemma}

The next result is a variant of a classical smoothing inequality. The proof is similar to that of~\citet[Lemma 11.2]{CCK23} and hence omitted.

\begin{lemma}\label{lem:smoothing}
Let $p,n\ge2$ and $n^{-1/2}\le\sigma\le1$. Let $S$ and $G$ be random vectors in $\R^p$, where $G$ is centered Gaussian with
$\min_{j\leq p}\Var(G_j)\ge1$, and let $W$ be a random vector in $\R^p$, independent of $S$ and of $G$, with $W_j\sim\mathcal N(0,2\sigma^2)$
for all $j\leq p$. Then,
\[
\sup_{\beta\in\R^p}\Big|\PP(S\le\beta)-\PP(G\le\beta)\Big|\, \le \, 
\sup_{\beta\in\R^p}\Big|\PP(S+W\le\beta)-\PP(G+W\le\beta)\Big|+C\sigma\log(pn).
\]
\end{lemma}
The following elementary facts are used in Section~\ref{sec:HS}. Below, we use $s_+=\max(s,0)$ to denote the non-negative part of $s\in\R$.

\begin{lemma}\label{lem:hazardfacts}
If $\lambda=(\log \Phi)'=\phi/\Phi$ and $\rho=-\lambda'$, then the following statements hold.
\begin{itemize}
\item[(a)] The function $\rho$ is decreasing, with $0<\rho<1$. Also, the derivatives $\rho'$ and $\rho''$ are bounded, and
$$\abs{\rho'(s)}\le(1+s_+)\,\rho(s) \qquad \text{and}\qquad \lambda(s)\le1+\abs s$$
for all real $s$.
\item[(b)] We have $\lambda(s)\le 2(1+s_+)\log(1/\Phi(s))$ for all real $s$, as well as $$\lambda(s)\le2\phi(s) \qquad \text{and} \qquad
\rho(s)\le2(1+s)\phi(s)$$
for all $s\ge0$.
\item[(c)] If $pn$ is sufficiently large, and $s,\tilde s\in\R$ satisfy $\abs{\tilde s-s}\le2(\log(pn))^{-3/2}$, then
$$\rho(\tilde s)\le2\rho(s)+(pn)^{-5}.$$
\end{itemize}
\end{lemma}

\begin{proof} 
All three parts are consequences of the following standard facts, where $\xi\sim\mathcal N(0,1)$ and $m(s)=s+\lambda(s)$. First, direct
differentiation gives $\lambda'=-\lambda m$, and so $\rho=\lambda m$. Second, integration by parts gives
\[
1-\rho(s)=\Var(\xi\mid \xi\le s).
\]
Third, the Mills ratio satisfies $\phi(x)/(1+x)\le1-\Phi(x)<\phi(x)/x$ for $x>0$, where the lower bound is Birnbaum's inequality~\citep{Bir42}.

\medskip\noindent\emph{Part (a).} 
The variance identity gives $\rho<1$, and $\rho$ is decreasing because $\Var(\xi\mid\xi\leq s)$ is increasing in
$s$~\citep[Proposition 1]{HH90}. The Mills ratio bounds give $0<m\le1+s_+$ and $\lambda(s)\le1+\abs s$, and so $\rho=\lambda m>0$. For the derivatives, note that $m'=1-\rho$. Differentiating $\rho=\lambda m$ then gives
\[
\rho'=-\lambda\big(m^2-(1-\rho)\big)\qquad\text{and}\qquad \rho''=-\rho'(m+\lambda)-2\rho(1-\rho).
\]
Since $\rho$ is decreasing, the first formula gives $\abs{\rho'}\le\lambda m^2=\rho m\le(1+s_+)\rho$. Using this in the second formula, together with $\lambda m=\rho$ and $\rho<1$, gives
\[
\abs{\rho''}\le\rho m^2+\rho^2+2\rho(1-\rho)\le(1+s_+)^2\rho(s)+1.
\]
Both bounds are bounded functions of $s$, because $\rho<1$, and because $\rho(s)\le2(1+s)\phi(s)$ for $s\ge0$ by part (b), whose proof uses only $m\le1+s_+$.

\medskip\noindent\emph{Part (b).} We treat the cases $s\ge0$ and $s<0$ separately. First let $s\ge0$. Then $\Phi(s)\ge\frac12$ gives
$\lambda(s)\le2\phi(s)$, and together with $m(s)\le1+s$ from (a), this gives $\rho(s)=\lambda(s)m(s)\le2(1+s)\phi(s)$. In addition, the lower
Mills ratio bound and the inequality $1-\Phi(s)\le\log(1/\Phi(s))$ give $\phi(s)\le(1+s)\log(1/\Phi(s))$, which proves the first bound in (b)
when $s\ge0$.

Next, let $s<0$. Since $m>0$, we have $\lambda(r)\ge-r$ for all $r\in\R$, and so
\[
\log\Phi(0)-\log \Phi(s)=\int_s^0\lambda(r)\,dr\ge \frac{s^2}2.
\]
Combining this with the bound $\lambda(s)\le1-s$ for negative $s$ from (a), and the elementary inequality $1-s\le2\log2+s^2$, shows that
$\lambda(s)\le2\log(1/\Phi(s))$, which proves the first bound in (b) when $s<0$.

\medskip\noindent\emph{Part (c).} Since $\rho$ is decreasing, we may assume that $\tilde s=s-\delta$ with $0\le\delta\le2(\log(pn))^{-3/2}$. We
treat the cases $s\le\log(pn)$ and $s>\log(pn)$ separately. In the first case, the bound $\abs{\rho'(t)}\le(1+\log(pn))\rho(t)$ holds for all
$t\le s$ by (a), and integrating it over $[s-\delta,s]$ gives
\[
\rho(s-\delta)\le\rho(s)\,e^{\delta(1+\log(pn))}\le2\rho(s).
\]
In the second case, we have $s-\delta\ge\log(pn)-1\ge1$, and so the bound on $\rho$ in (b) gives
$\rho(s-\delta)\lesssim\log(pn)\,\phi(\log(pn)-1)\le(pn)^{-5}$. In both cases, the last inequality holds because $pn$ is large.
\end{proof}

\begin{lemma}\label{lem:condgauss}
Let $Z\sim\gamma_d$, let $A\subset\R^d$ be a closed convex set with $\gamma_d(A)>0$, and let $m_A=\E(Z\mid Z\in A)$ and
$\mathfrak S_A=\Cov(Z\mid Z\in A)$.
\begin{itemize}
\item[(a)] The conditional mean satisfies $\norm{m_A}_2^2\le2\log(1/\gamma_d(A))$.
\item[(b)] The function $x\mapsto\gamma_d(A-x)$ is $C^\infty$, positive and log-concave on $\R^d$. Also, the function $x\mapsto \log\gamma_d(A-x)$ has gradient $m_A$ and Hessian $\mathfrak S_A-I_d\preceq0$ at $x=0$.
\end{itemize}
The same statements hold for the standard Gaussian measure $\gamma_E$ of a linear subspace $E$, with $I_d$ replaced by the
orthogonal projection onto $E$.
\end{lemma}
\begin{proof}
\emph{Part (a).} Let $\nu$ be the conditional distribution of $Z$ given $Z\in A$, so that $d\nu/d\gamma_d=\1_A/\gamma_d(A)$ and $d_{\textup{KL}}(\nu\,\|\,\gamma_d)=\log(1/\gamma_d(A))$.
 Applying the transportation inequality~\citep[Lemma~4.18]{BLM13} to the function
$z\mapsto\ip{m_A}{z}/\norm{m_A}_2$ (which is a standard Gaussian variable under $\gamma_d$) it follows that $\norm{m_A}_2\le\sqrt{2\log(1/\gamma_d(A))}$.\\

\noindent\emph{Part (b).} We write $\gamma_d(A-x)=\int_A\phi_d(z-x)\,dz$. The positivity and smoothness are clear. To address the log-concavity, observe that the convexity of $A$ implies that the function $(x,z)\mapsto\1_A(z)\phi_d(z-x)$ is log-concave. This means  $x\mapsto\gamma_d(A-x)$ is obtained from marginalizing a log-concave function, and is hence log-concave~\citep{Pre73}.  For the gradient and Hessian expressions, first note that when $\nabla_x\phi_d(z-x)$ and $\nabla_x^2\phi_d(z-x) $ are evaluated at $x=0$, they become $z\,\phi_d(z)$ and
$(zz^\top-I_d)\,\phi_d(z)$. Second, let $\nu$ be as in the proof of part (a) and recall that $d\nu/d\gamma_d=\1_A/\gamma_d(A)$.  A direct calculation of $\nabla_x \log\gamma_d(A-x)|_{x=0}$ then gives
$$ \nabla_x \log \gamma_d(A-x)|_{x=0}=\int_{\R^d} z (d\nu/d\gamma_d)\,d\gamma_d(z)= \int_{\R^d} z\,d\nu(z)=m_A,$$
and the calculation for the Hessian is similar. Meanwhile, the condition $\mathfrak{S}_A-I_d\preceq 0$ is due to the log-concavity of $x\mapsto \gamma_d(A-x)$.
Lastly, for a subspace $E$, the same arguments apply after identifying $E$ with a Euclidean space of the same dimension.
\end{proof}

\subsection{Polyhedra with linearly independent normal vectors}\label{subsec:polyhedra}

In this subsection, we prove Lemmas~\ref{lem:structure} and~\ref{lem:facetfn}, and two further
lemmas, on the derivatives of the Gaussian measures of $P$ and of its facets under translations.

\begin{proof}[Proof of Lemma~\ref{lem:structure}] The  representation of the $j$th facet as $F_j=a_jn_j+K_j$ is straightforward to check, as is its closedness and convexity. To show that $q_j$ is positive, first let $H_j=\{z\in\R^p:\ip{n_j}z=a_j\}$ denote the hyperplane containing the $j$th facet, so
that $F_j=P\cap H_j$. By~\eqref{eqn:linearlyindep} there is a point $z^{\circ}\in\R^p$ with $\ip{n_j}{z^{\circ}}=a_j$ and $\ip{n_l}{z^{\circ}}=a_l-1$ for
$l\ne j$.  Consequently, if a point $z\in H_j$ satisfies $\norm{z-z^{\circ}}_2\le1$, then,
\[
\ip{n_l}z\le\ip{n_l}{z^{\circ}}+\norm{z-z^{\circ}}_2\leq a_l\qquad\text{for all }l\ne j.
\]
This means that the unit $\ell_2$-ball of $H_j$ centered at $z^{\circ}$ is contained in $F_j$, and hence, that $F_j$ has non-empty interior relative to $H_j$. Lastly, because $K_j$ and $E_j$ are obtained from a common shift of $F_j$ and $H_j$, this implies that $K_j$ has non-empty interior relative to $E_j$, and thus $q_j=\gamma_{E_j}(K_j)$ is positive.

Finally, we verify the formula $\nu_j(P)=\phi(a_j)\,q_j$. For any $z'\in E_j$ we have $\norm{a_jn_j+z'}_2^2=a_j^2+\norm{z'}_2^2$, and so $\phi_p(a_jn_j+z')=\phi(a_j)\phi_{p-1}(z')$. Since the subspace $E_j$ is isometric to the hyperplane $H_j$ mentioned a moment ago, it follows that
\[
\nu_j(P)=\int_{F_j}\phi_p\,d\mathcal H^{p-1}=\int_{K_j}\phi_p(a_j n_j+z')d\mathcal{H}^{p-1}(z')=\phi(a_j)\int_{K_j}\phi_{p-1}\,d\mathcal H^{p-1}=\phi(a_j)\,q_j.
\]
\end{proof}

The content of the following result might be regarded as classical, but we provide a proof since the statements do not seem to be easily referenced in the form below.

\begin{lemma}\label{lem:transl}
Let $P\subset\R^p$ be a polyhedron as in~\eqref{eqn:Polygeneral} with $p$ facets having linearly independent normal vectors, and let $h=\pi\phi_p$ with $\pi$ being a polynomial on $\R^p$. Then, the function $x\mapsto\int_{P-x}h\,dz$ is $C^\infty$
on $\R^p$, and for all $x,x'\in\R^p$ we have
\begin{equation}\label{eq:transl}
\partial_{x'}\int_{P-x}h\,dz=-\sum_{j=1}^p\ip{n_j}{x'}\int_{F_j-x}h\,d\mathcal H^{p-1},
\end{equation}
where $F_j-x$ is the $j$th facet of $P-x$. Furthermore, the following hold.
\begin{itemize}
\item[(a)] The first derivatives are given by $\partial_{x'}\gamma_p(P-x)=-\sum_j\ip{n_j}{x'}\,\nu_j(P-x)$.
\item[(b)] The second and third derivatives are given by
\[
\begin{split}
\nabla_x^2\gamma_p(P-x)&=\int_{P-x}(zz^\top-I_p)\,\phi_p\,dz,\\[1mm]
\partial_{x'}\nabla_x^2\gamma_p(P-x)&=-\sum_{j=1}^p\ip{n_j}{x'}\,M_j(P-x),
\end{split}
\]
where we put $M_j(P)=\int_{F_j}(zz^\top-I_p)\,\phi_p\,d\mathcal H^{p-1}$.
\item[(c)] For every symmetric matrix $\Gamma$ with $0\preceq\Gamma\preceq I_p$ and every $j$, we have
$$\abs{\tr(\Gamma M_j(P))}\le\Big(1+a_j^2+\norm{m_j}_2^2+\tr(\Pi_j-\mathfrak S_j)\Big)\,\nu_j(P).$$
\end{itemize}
\end{lemma}

\begin{proof}
We first prove~\eqref{eq:transl}. Since $\int_{P-x}h\,dz=\int_Ph(z-x)\,dz$, differentiation under the integral shows that this
function is $C^\infty$, with
\[
\partial_{x'}\int_{P-x}h\,dz=-\int_{P-x}\ip{\nabla h}{x'}\,dz.
\]
The boundary of $P-x$ is the union of its facets $F_j-x$, and the outer unit normal on $F_j-x$ is $n_j$. So the divergence theorem, applied to
the vector field $hx'$ on $P-x$, gives $\int_{P-x}\ip{\nabla h}{x'}\,dz=\sum_{j=1}^p\ip{n_j}{x'}\int_{F_j-x}h\,d\mathcal H^{p-1}$, which
is~\eqref{eq:transl}. (Although $P-x$ is unbounded, the divergence theorem can still be used by intersecting $(P-x)$ with the $\ell_2$ ball of radius $r$ and then letting  $r\to\infty$.)

\medskip\noindent\emph{Part (a).} Taking $\pi=1$ in~\eqref{eq:transl} proves (a).

\medskip\noindent\emph{Part (b).} Since $\nabla^2\phi_p(z)=(zz^\top-I_p)\phi_p(z)$, differentiating $\gamma_p(P-x)=\int_P\phi_p(z-x)\,dz$ twice
under the integral gives the first statement of (b). The second formula follows by applying~\eqref{eq:transl} entrywise to the formula for $\nabla_x^2\gamma(P-x)$, while taking $h(z)=\pi(z)\phi_p(z)$ with the polynomial
$\pi(z)=z_iz_l-\1\{i=l\}$.\\

\noindent\emph{Part (c).} We first express $M_j(P)$ in terms of the parameters of the
$j$th section $K_j$. By Lemma~\ref{lem:structure}, every point of $F_j$ can be written as
$z=a_jn_j+z'$ with $z'\in K_j$, and for such points we have
$\phi_p(z)=\phi(a_j)\phi_{p-1}(z')$. Changing variables from $z$ to $z'$ in the
definition of $M_j(P)$ gives
\begin{equation*}
 M_j(P)\;=\;\phi(a_j)\int_{K_j}\Big((a_jn_j+z')(a_jn_j+z')^{\top}-I_p\Big)
 \,\phi_{p-1}(z')\,d\mathcal{H}^{p-1}(z').
\end{equation*}
Next, let $W\sim\gamma_{E_j}$ and $U=a_jn_j+W$. The integral above is the
expectation of $(UU^{\top}-I_p)1\{W\in K_j\}$, and so dividing and multiplying
by $q_j=\PP(W\in K_j)$ and using $\nu_j(P)=\phi(a_j)q_j$ leads to
\begin{equation*}
 M_j(P)\;=\;\nu_j(P)\,\E\big(UU^{\top}-I_p\,\big|\,W\in K_j\big).
\end{equation*}
To evaluate the conditional expectation, recall that $m_j$ and $\mathfrak{S}_j$
are the mean and covariance of $W$ given $W\in K_j$. Hence, the conditional
mean and covariance of $U$ are $a_jn_j+m_j$ and $\mathfrak{S}_j$, which gives
$\E(UU^{\top}\,|\,W\in K_j)=(a_jn_j+m_j)(a_jn_j+m_j)^{\top}+\mathfrak{S}_j$.
Combining this with the algebraic identity $I_p=n_jn_j^{\top}+\Pi_j$, we obtain
\begin{equation*}
 M_j(P)\;=\;\nu_j(P)\Big[(a_jn_j+m_j)(a_jn_j+m_j)^{\top}-n_jn_j^{\top}
 -(\Pi_j-\mathfrak{S}_j)\Big].
\end{equation*}
It remains to bound $\operatorname{tr}(\Gamma M_j(P))$. The matrix in the brackets
above has the form $A_1-A_2$, where $A_1=(a_jn_j+m_j)(a_jn_j+m_j)^{\top}$ and
$A_2=n_jn_j^{\top}+(\Pi_j-\mathfrak{S}_j)$ are positive semidefinite, with
$\Pi_j-\mathfrak{S}_j\succeq 0$ (by Lemma~\ref{lem:condgauss}(b) on $E_j$). Since
$0\preceq\Gamma\preceq I_p$, every positive semidefinite matrix $A$ satisfies
$0\le\operatorname{tr}(\Gamma A)\le\operatorname{tr}(A)$,
and so
\begin{equation*}
 \frac{|\operatorname{tr}(\Gamma M_j(P))|}{\nu_j(P)}
 \;\le\;\operatorname{tr}(A_1)+\operatorname{tr}(A_2)
 \;=\;1+a_j^2+\|m_j\|_2^2+\operatorname{tr}(\Pi_j-\mathfrak{S}_j),
\end{equation*}
where the last step uses the orthogonality of $m_j$ and $n_j$.
\end{proof}

\begin{lemma}\label{lem:facetlog}
Let $P\subset\R^p$ be a polyhedron as in~\eqref{eqn:Polygeneral} with $p$ facets having linearly independent normal vectors. Then, for all $j\le p$ and $x\in\R^p$, the $j$th section of $P-x$ can be represented as $K_j-\Pi_jx$, and 
\[
\nu_j(P-x)=\phi\big(a_j-\ip{n_j}{x}\big)\,q_j(P-x).
\]
Furthermore, both of the functions $x\mapsto q_j(P-x)$ and $x\mapsto \nu_j(P-x)$ are positive, $C^\infty$, and log-concave on $\R^p$. Lastly, at $x=0$, the function $x\mapsto\log\nu_j(P-x)$ has gradient
$a_jn_j+m_j$ and Hessian $\mathfrak S_j-n_jn_j^\top-\Pi_j \preceq0$.
\end{lemma}

\begin{proof}
We first identify the $j$th section of $P-x$. Note that the $j$th projection matrices of $P-x$ and $P$ are the same, because they have the same normal vectors. Also, it is simple to check that the $j$th section $K_j$ and facet $F_j$ of $P$ are related by $\Pi_j F_j=K_j$. So, because the $j$th facet of $P-x$ is given by $F_j-x$, it follows that the $j$th section of $P-x$ is given by $\Pi_j(F_j-x)=K_j-\Pi_j x$.  Taking Gaussian measures gives $q_j(P-x)=\gamma_{E_j}(K_j-\Pi_jx)$, and applying
Lemma~\ref{lem:structure} to $P-x$ gives the formula for $\nu_j(P-x)$.

It remains to establish the regularity properties and the derivatives at $x=0$. By the formula for $\nu_j(P-x)$, the function
$\log\nu_j(P-x)$ is the sum of $\log\phi(a_j-\ip{n_j}x)$ and $\log\gamma_{E_j}(K_j-\Pi_jx)$. The first term $\log\phi(a_j-\ip{n_j}x)$  is clearly concave, with gradient $a_jn_j$ and Hessian $-n_jn_j^\top$ at $x=0$. For the
second term $\log\gamma_{E_j}(K_j-\Pi_jx)$, recall from Lemma~\ref{lem:structure} that $K_j$ is closed and convex with $q_j>0$. So Lemma~\ref{lem:condgauss}(b) on $E_j$ shows
that $x'\mapsto\gamma_{E_j}(K_j-x')$ is positive, $C^\infty$ and log-concave on $E_j$, and that its logarithm has gradient $m_j$ and Hessian
$\mathfrak S_j-\Pi_j$ at $x'=0$. Composing with the orthogonal projection $\Pi_j$ of $\R^p$ onto $E_j$ preserves these properties, with the same
gradient and Hessian at $x=0$. Adding the two contributions gives the claims.
\end{proof}

\begin{proof}[Proof of Lemma~\ref{lem:facetfn}]
The formula for $\psi_j(y)$ follows from its definition combined with $\nu_j(P)=\phi(a_j)q_j$ (Lemma~\ref{lem:structure}).
Taking logarithms in this formula gives
\[
\log\psi_j(y)=\log q_j(y)-\tfrac12a_j(y)^2-\tfrac12\log\big(2\pi(u+\sigma^2)\big),
\]
and applying $\partial_w$ gives~\eqref{eq:s7-Dlogpsi}. Next, since we have
$P(y+w)=P(y)-\Sigma_u^{-1/2}w$ by Lemma~\ref{lem:gu}(a), the positivity, smoothness, and log-concavity of $\psi_j$ and $q_j$ follow from
Lemma~\ref{lem:facetlog}.

Before proceeding to the upper bounds in the statement of Lemma~\ref{lem:facetfn}, we first compute formulas for the derivatives of $\log\psi_j$. The chain rule and the last statement of Lemma~\ref{lem:facetlog} (replacing $x$ with $\Sigma_u^{-1/2}w$) give
\[
\nabla\log\psi_j(y)=\Sigma_u^{-1/2}\big(a_j(y)n_j+m_j(y)\big)\qquad\text{and}\qquad-\nabla^2\log\psi_j(y)=\Sigma_u^{-1/2}\Lambda_j(y)\,\Sigma_u^{-1/2},
\]
where we put $\Lambda_j(y)=n_jn_j^\top+\Pi_j-\mathfrak S_j(y)\succeq0$. To develop upper bounds, we use $u\Sigma\preceq\Sigma_u$, which is the same as
$\Sigma_u^{-1/2}\Sigma\Sigma_u^{-1/2}\preceq u^{-1}I_p$. Applying this matrix inequality to the gradient formula gives
$$\ip{\nabla\log\psi_j(y)}{\Sigma\,\nabla\log\psi_j(y)}\,\le\, \textstyle \frac{1}{u}\norm{a_j(y)n_j+m_j(y)}_2^2,$$
which is~\eqref{eq:facetfn1} because $m_j(y)$ and $n_j$ are orthogonal.
Similarly, combining the matrix inequality and the Hessian formula gives
\begin{equation*}
\begin{split}
\tr\big(\Sigma(-\nabla^2\log\psi_j(y))\big)
&\le\textstyle\frac{1}{u}\tr\Lambda_j(y),
\end{split}
\end{equation*}
and this reduces to~\eqref{eq:facetfn2} because $\tr\Lambda_j(y)=1+\tr(\Pi_j-\mathfrak S_j(y))$.

Now we conclude the proof by applying the gradient and Hessian inequalities to prove~\eqref{eq:facetfn3}--\eqref{eq:facetfn5}. Since $\tr(\Pi_j-\mathfrak S_j(y))\ge0$ (Lemma~\ref{lem:condgauss}(b)), the quantities $\nu_j(y)$, $(a_j(y)^2+\norm{m_j(y)}_2^2)\nu_j(y)$ and
$\tr\Lambda_j(y)\nu_j(y)$ are all at most the $j$th term of $\Fc(y)$ in Lemma~\ref{lem:sensitivityformula}. Finally, due to the basic inequality 
$\psi_j(y)\le\nu_j(y)/\sqrt u$, the bounds~\eqref{eq:facetfn3}--\eqref{eq:facetfn5} follow by multiplying~\eqref{eq:facetfn1}
and~\eqref{eq:facetfn2} by $\psi_j(y)$ and summing over $j$.
\end{proof}

\begin{proof}[Proof of Lemma~\ref{lem:lap3}]
Fix $y,w$ and put $x'=\Sigma_u^{-1/2}w$. By Lemma~\ref{lem:gu}, $P(y+\tau w)-x=P(y)-(x+\tau x')$ and hence
$\LapS g_u(y+\tau w)=\frac1u\tr\big(\Gamma\,\nabla_x^2\gamma_p(P(y)-x)|_{x=\tau x'}\big)$. Differentiating at $\tau=0$ with
Lemma~\ref{lem:transl}(b), applied to $P(y)$ at $x=0$,
\[
\partial_w\LapS g_u(y)=-\frac1u\sum_{j=1}^p\ip{n_j}{x'}\tr\big(\Gamma M_j(P(y))\big),
\]
where $M_j(P(y))$ is as defined in that lemma.
Here $\ip{n_j}{x'}=w_j/\sqrt{u+\sigma^2}$ with $u+\sigma^2\ge u$, and Lemma~\ref{lem:transl}(c) gives
$\abs{\tr(\Gamma M_j(P(y)))}\le\big(1+a_j(y)^2+\norm{m_j(y)}_2^2+\tr(\Pi_j-\mathfrak S_j(y))\big)\,\nu_j(y)$. Summing over $j$ and using Lemma~\ref{lem:sensitivityformula} gives the claim.
\end{proof}

\section{Proof of Theorem~\ref{thm:bt}}\label{sec:bootstrap}
It is sufficient to prove the result when the class of rectangles is replaced with the class of orthants, using essentially the same reasoning as in the proof of Lemma~\ref{lem:reduction}. 
Also, it can be checked that the discrepancy measure $d(\tilde\Sigma,\Sigma)$ between the covariance matrices of $\tilde G\sim\mathcal{N}(0,\tilde \Sigma)$ and $G\sim\mathcal{N}(0,\Sigma)$ does not change when $\tilde\Sigma$ and $\Sigma$ are replaced by the $2p\times 2p$ covariance matrices of $(\tilde G,-\tilde G)$ and $(G,-G)$. 

Thus, we will show that if $d(\tilde\Sigma,\Sigma)\leq \delta$ for some $\delta>0$, and if $\Sigma$ has ones along the diagonal, then (with $p$ in place of $2p$)
\begin{equation}\label{eqn:Gcompstart}
    \sup_{\beta\in\R^p}\Big|\PP(\tilde G\leq \beta)-\PP(G\leq \beta)\Big|\lesssim \delta \log(p)\log(1+1/\delta).
\end{equation}
Also, we may assume going forward that $\delta\leq 1/2$ holds without loss of generality, for otherwise the statement of Theorem~\ref{thm:bt} is trivially true.

To proceed, we will re-use an initial portion of the proof of Theorem~\ref{thm:bounded} by taking $X_1,\dots,X_n$ to be i.i.d.~copies of $\tilde G$ with $n=\lceil 1/\delta^2\rceil$, so that the Gaussianity of $\tilde G$ implies $S=n^{-1/2}\sum_{i=1}^n X_i\overset{d}{=} \tilde G$.\footnote{Note that even though the statement of Theorem~\ref{thm:bounded} deals with bounded random vectors $X_1,\dots,X_n$ whose covariance matrix matches that of $G$, the arguments used to reach~\eqref{eqn:interpagain} do not rely on boundedness or on the matching covariance matrices.} Furthermore, if we let $\sigma=\delta$, then the arguments at the beginning of Section~\ref{sec:reductions} show that the supremum on the left side of~\eqref{eqn:Gcompstart} can be bounded as
\begin{equation}\label{eqn:interpagain}
      \sup_{\beta\in\R^p}\Big|\PP(\tilde G\leq \beta)-\PP(G\leq \beta)\Big| \lesssim \int_0^1 |\Psi'(t)|dt + \delta\log(p/\delta),
\end{equation}
where $\Psi(t)=\E g_{u(t)}(\sqrt t\,\tilde G)$ and $u(t)=1-t+\delta^2$. 
As in the proof of Lemma~\ref{lem:interp}, the derivative of $\Psi$ can be calculated as $\Psi'(t)=\frac1{2\sqrt t}\E\ip{\nabla g_u(\sqrt t\tilde G)}{\tilde G}-\frac12\E\LapS g_u(\sqrt t\tilde G)$ for $t\in(0,1]$, and applying Stein's identity~\eqref{eqn:steinidentity} to the first term gives
\[
\Psi'(t)=\tfrac12\,\E\tr\big((\tilde\Sigma-\Sigma)\nabla^2g_{u(t)}(\sqrt t\,\tilde G)\big).
\]
In the later part of this section, Lemmas~\ref{lem:bt-pointwise} and~\ref{lem:bt-wedge} will collectively show that the trace in this expression can be bounded as
\begin{equation}\label{eq:bt-key}
\big|\tr\big((\tilde\Sigma-\Sigma)\nabla^2g_u(y)\big)\big| \, \lesssim \, \frac{\delta \log(2p)}{(u+\delta^2)},
\end{equation}
for all $y\in\R^p$ and $u>0$. Combining this with $\int_0^1(u(t)+\delta^2)^{-1}dt\le2\log(1/\delta)$, we obtain
\begin{equation*}
          \sup_{\beta\in\R^p}\Big|\PP(\tilde G\leq \beta)-\PP(G\leq \beta)\Big| \lesssim \delta\log(2p)\log(1/\delta) + \delta\log(p/\delta),
\end{equation*}
which implies the stated result.\qed

\subsection{Supporting results for the bound~\eqref{eq:bt-key}}\label{subsec:bt-hessian}
This subsection collects three lemmas that establish the bound~\eqref{eq:bt-key}. Throughout the subsection, $u>0$ is fixed, and we write $a(y)=(a_1(y),\dots,a_p(y))=(\beta-y)/\sqrt{u+\sigma^2}$ for the vector of offsets
of $P(y)$. Also, let $R$ denote the correlation matrix $R=\frac{1}{u+\sigma^2}\Sigma_u$, and let $Y\sim\mathcal{N}(0,R)$. Likewise, let $\phi_R$ and $\Phi_R$ denote the density and distribution function of $Y$ so that $\Phi_R(a)=\PP(Y\leq a)$ for all $a\in\R^p$.

The following lemma gives integral formulas for the first and second partial derivatives of $\Phi_R$, together with an identity for the
diagonal entries of $R\nabla^2\Phi_R$.

\begin{lemma}\label{lem:bt-stein}
For all $a\in\R^p$ and all $1\leq j\ne l\leq p$, we have
\[
\partial_j\Phi_R(a)=\int_{\mathrlap{\raisebox{-0.6ex}{$\scriptstyle\{z\le a,\ z_j=a_j\}$}}}\mkern12mu\phi_R(z)\,d\mathcal H^{p-1}(z)
\qquad\text{and}\qquad
\partial_j\partial_l\Phi_R(a)=\int_{\mathrlap{\raisebox{-0.6ex}{$\scriptstyle\{z\le a,\ z_j=a_j,\ z_l=a_l\}$}}}\mkern12mu\phi_R(z)\,d\mathcal H^{p-2}(z).
\]
Moreover, for all $a\in\R^p$ and $j\le p$,
\begin{equation}\label{eq:bt-stein}
\sum_{l=1}^pR_{jl}\,\partial_j\partial_l\Phi_R(a)=-a_j\,\partial_j\Phi_R(a).
\end{equation}
\end{lemma}

\begin{proof} The first two formulas are straightforward.
To prove the third formula~\eqref{eq:bt-stein}, we write $\Phi_R(a)=\int_{\{z\le0\}}\phi_R(z+a)\,dz$ and differentiate under the integral sign,
which gives $\nabla\Phi_R(a)=\int_{\{z\le a\}}\nabla\phi_R(z)\,dz$. The left side of~\eqref{eq:bt-stein} can be rewritten as
\begin{equation*}
    \begin{split}
         \sum_{l=1}^pR_{jl}\,\partial_j\partial_l\Phi_R(a) = \partial_j \ip{e_j}{R\nabla \Phi_R(a)}
         =\int_{\mathrlap{\raisebox{-0.6ex}{$\scriptstyle\{z\le a\}$}}}\mkern12mu\partial_j\big(R\nabla\phi_R\big)_j(z)\,dz.
    \end{split}
\end{equation*}
Furthermore, the basic formula $ R\nabla \phi_R(z)=-\phi_R(z) z$ allows the last integral to be calculated as
$$ 
-\int_{\mathrlap{\raisebox{-0.6ex}{$\scriptstyle\{z\le a\}$}}}\mkern12mu\partial_j\big(z_j\phi_R(z)\big)\,dz\\
=-a_j\int_{\mathrlap{\raisebox{-0.6ex}{$\scriptstyle\{z\le a,\ z_j=a_j\}$}}}\mkern12mu\phi_R(z)\,d\mathcal H^{p-1}(z).$$
By the earlier formula for $\partial_j\Phi_R(a)$, the right side is $-a_j\,\partial_j\Phi_R(a)$.
\end{proof}

As notation for the next lemma, we define for each $1\leq j\neq l\leq p$ and $y\in \R^p$ the function
$$\nu_{jl}(y)=\partial_j\partial_l\Phi_R(a(y)).$$
This notation is motivated by the formula 
$\nu_j(y)=\partial_j\Phi_R(a(y))$, which will be shown below in~\eqref{eq:nuQ}.

\begin{lemma}\label{lem:bt-pointwise}
Let $\Sigma$ and $\tilde\Sigma$ be positive semidefinite $p\times p$ matrices such that $\Sigma_{jj}=1$ for all $j\le p$ and $d(\tilde\Sigma,\Sigma)\le\delta$ for some $\delta>0$. Then, for all $y\in\R^p$,
\[
\big|\tr\big((\tilde\Sigma-\Sigma)\nabla^2g_u(y)\big)\big|
\le \frac{\delta}{u+\sigma^2}\bigg(\sum_{j=1}^p\abs{a_j(y)}\,\nu_j(y)+2\sum_{j<l}(1-R_{jl}^2)\,\nu_{jl}(y)\bigg).
\]
\end{lemma}

\begin{proof}
To ease notation, we fix $y$ and write $a=a(y)$, $\nu_j=\nu_j(y)$ and $\nu_{jl}=\nu_{jl}(y)$. Noting that the covariance matrix of
$(\sqrt uG+\sigma N)/\sqrt{u+\sigma^2}$ is equal to $R$, it follows from the definition of $g_u$ that $g_u(y)=\Phi_R(a(y))$, and so the chain rule leads to
$$\nabla^2g_u(y)=\textstyle\frac{1}{(u+\sigma^2)}\nabla^2\Phi_R(a).$$
Furthermore, since the off-diagonal entries of
$\nabla^2\Phi_R(a)$ are the numbers $\nu_{jl}$, and since $\Sigma_{jj}=1$, we have
\begin{equation}\label{eq:bt-tracestart}
\tr\big((\tilde\Sigma-\Sigma)\nabla^2\Phi_R(a)\big)=\sum_{j=1}^p(\tilde\Sigma_{jj}-1)\,\partial_j^2\Phi_R(a)+2\sum_{j<l}(\tilde\Sigma_{jl}-\Sigma_{jl})\,\nu_{jl}.
\end{equation}
Towards computing a formula for $\partial_j^2\Phi_R(a)$, we show that the first derivative $\partial_j\Phi_R(a)$ is related to $\nu_j(y)$. Indeed, the chain rule gives $\partial_jg_u(y)=-\partial_j\Phi_R(a)/\sqrt{u+\sigma^2}$, and
comparing this with~\eqref{eqn:psidef} gives
\begin{equation}\label{eq:nuQ}
\nu_j(y)=\partial_j\Phi_R(a(y)),
\end{equation}
for all $j\leq p$.
Now, to get $\partial_j^2\Phi_R(a)$, note that this quantity appears in the $l=j$ term of the sum on the left side
of~\eqref{eq:bt-stein}, and that $R_{jj}=1$. Based on this observation, rearranging~\eqref{eq:bt-stein} and recalling the definition of $\nu_{jl}$ leads to
\[
\partial_j^2\Phi_R(a)=-a_j\nu_j-\sum_{l\,:\,l\ne j}R_{jl}\nu_{jl}.
\]
Substituting this into~\eqref{eq:bt-tracestart} gives
\[
\tr\big((\tilde\Sigma-\Sigma)\nabla^2\Phi_R(a)\big)
=-\sum_{j=1}^p(\tilde\Sigma_{jj}-1)\,a_j\nu_j-\sum_{j\ne l}(\tilde\Sigma_{jj}-1)R_{jl}\,\nu_{jl}
+2\sum_{j<l}(\tilde\Sigma_{jl}-\Sigma_{jl})\,\nu_{jl}.
\]
To convert the sum over $j\ne l$ into a sum over $j<l$, each term $(\tilde\Sigma_{jj}-1)R_{jl}\nu_{jl}$ will become
$R_{jl}(\tilde\Sigma_{jj}+\tilde\Sigma_{ll}-2)\,\nu_{jl}$. Thus,  we obtain
\begin{equation}\label{eq:bt-trace}
\tr\big((\tilde\Sigma-\Sigma)\nabla^2\Phi_R(a)\big)=-\sum_{j=1}^p(\tilde\Sigma_{jj}-1)\,a_j\nu_j+\sum_{j<l}c_{jl}\,\nu_{jl},
\end{equation}
where we define the coefficients
$c_{jl}=2(\tilde\Sigma_{jl}-\Sigma_{jl})-R_{jl}(\tilde\Sigma_{jj}+\tilde\Sigma_{ll}-2)$.

Because the quantities $\nu_j$ and $\nu_{jl}$ are all non-negative by
Lemma~\ref{lem:bt-stein}, it remains to show that $\abs{\tilde\Sigma_{jj}-1}\le\delta$ and $\abs{c_{jl}}\le2\delta(1-R_{jl}^2)$.
By the assumption $d(\tilde\Sigma,\Sigma)\leq \delta$, we have
\begin{equation}\label{eqn:dbound}
\abs{\ip{v}{(\tilde\Sigma-\Sigma)v}}\le\delta\,\ip{v}{\Sigma v}
\end{equation}
for every $v\in\mathcal V$.
Taking $v=e_j$ gives $\abs{\tilde\Sigma_{jj}-1}\le\delta$, which handles the first sum in~\eqref{eq:bt-trace}. For the second sum, it will be convenient to introduce the following notation with the indices $j<l$ held fixed,
\[
\omega_\pm=\ip{e_j\pm e_l}{(\tilde\Sigma-\Sigma)(e_j\pm e_l)}.
\]
Taking $v=e_j\pm e_l$ in the bound~\eqref{eqn:dbound} gives $\abs{\omega_\pm}\le2\delta(1\pm\Sigma_{jl})$.
To relate $\omega_\pm$ to $c_{jl}$, we expand 
$$\omega_\pm=(\tilde\Sigma_{jj}-1)+(\tilde\Sigma_{ll}-1)\pm2(\tilde\Sigma_{jl}-\Sigma_{jl}).$$
Comparing this with the definition of $c_{jl}$, it can be checked that
\[
c_{jl}=\tfrac12(1-R_{jl})\,\omega_+-\tfrac12(1+R_{jl})\,\omega_-.
\]
Since $\abs{R_{jl}}\le1$, the coefficients $\frac12(1\mp R_{jl})$ are nonnegative, and so the bounds on $\abs{\omega_\pm}$ give
\[
\begin{split}
\abs{c_{jl}}&\le\tfrac12(1-R_{jl})\abs{\omega_+}+\tfrac12(1+R_{jl})\abs{\omega_-}\\[0.2cm]
&\le\delta(1-R_{jl})(1+\Sigma_{jl})+\delta(1+R_{jl})(1-\Sigma_{jl})\\[0.2cm]
&=2\delta\,(1-R_{jl}\Sigma_{jl})\ \le\ 2\delta\,(1-R_{jl}^2),
\end{split}
\]
where the last step holds because $R_{jl}\Sigma_{jl}=u\Sigma_{jl}^2/(u+\sigma^2)\ge R_{jl}^2$.
\end{proof}

\begin{lemma}\label{lem:bt-wedge}
For all $y\in\R^p$,
\begin{align}
\sum_{j=1}^p\abs{a_j(y)}\,\nu_j(y)&\lesssim\log(2p),\label{eq:bt-facetsum}\\[0.1cm]
\sum_{j<l}(1-R_{jl}^2)\,\nu_{jl}(y)&\lesssim\log(2p).\label{eq:bt-edgesum}
\end{align}
\end{lemma}

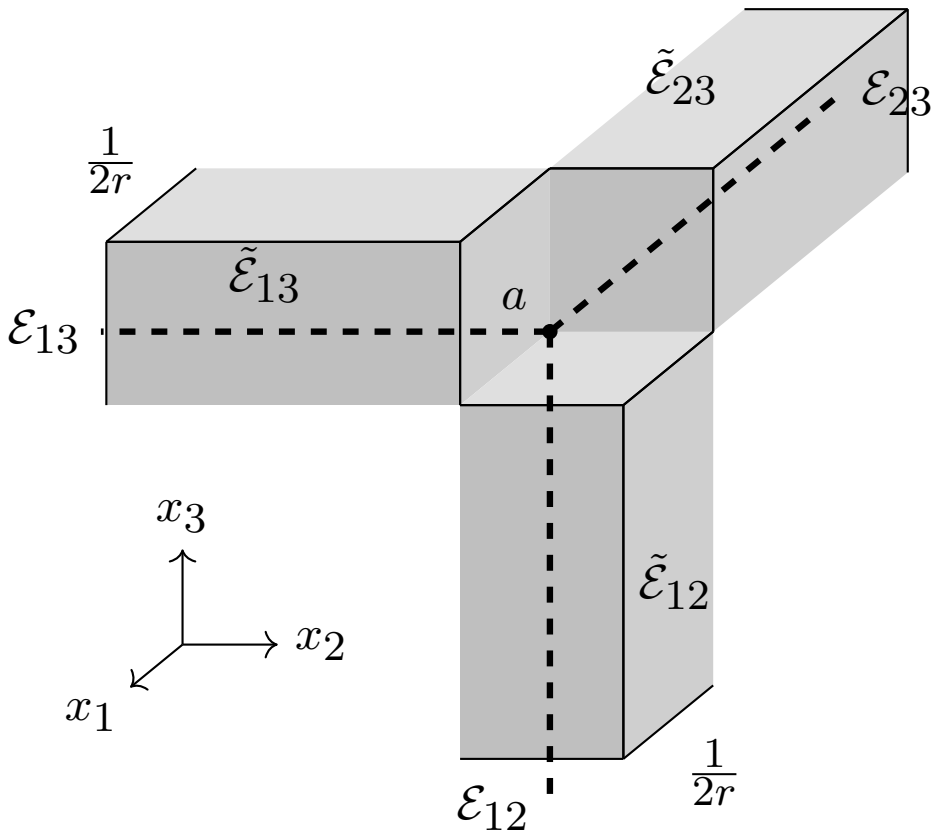
\begin{figure}[t]
\centering
\begin{tikzpicture}[x={(-0.55cm,-0.45cm)},y={(1cm,0cm)},z={(0cm,1cm)},scale=0.9,line join=round,font=\footnotesize]
\def\wh{1.2}\def\wL{2.6}
\fill[gray!50] (0,0,0) -- (0,\wh,0) -- (0,\wh,\wh) -- (0,0,\wh) -- cycle;
\fill[gray!38] (-\wL,\wh,0) -- (0,\wh,0) -- (0,\wh,\wh) -- (-\wL,\wh,\wh) -- cycle;
\fill[gray!25] (-\wL,0,\wh) -- (0,0,\wh) -- (0,\wh,\wh) -- (-\wL,\wh,\wh) -- cycle;
\draw (0,0,\wh) -- (0,\wh,\wh);
\draw (0,\wh,0) -- (0,\wh,\wh);
\draw (-\wL,\wh,\wh) -- (0,\wh,\wh);
\draw (-\wL,\wh,0) -- (-\wL,\wh,\wh);
\draw (0,\wh,0) -- (0,\wh,\wh);
\draw (-\wL,\wh,\wh) -- (0,\wh,\wh);
\draw (-\wL,0,\wh) -- (-\wL,\wh,\wh);
\draw (0,0,\wh) -- (0,\wh,\wh);
\fill[gray!50] (\wh,-\wL,0) -- (\wh,0,0) -- (\wh,0,\wh) -- (\wh,-\wL,\wh) -- cycle;
\fill[gray!38] (0,0,0) -- (\wh,0,0) -- (\wh,0,\wh) -- (0,0,\wh) -- cycle;
\fill[gray!25] (0,-\wL,\wh) -- (\wh,-\wL,\wh) -- (\wh,0,\wh) -- (0,0,\wh) -- cycle;
\draw (\wh,-\wL,\wh) -- (\wh,0,\wh);
\draw (\wh,-\wL,0) -- (\wh,-\wL,\wh);
\draw (\wh,0,0) -- (\wh,0,\wh);
\draw (0,0,\wh) -- (\wh,0,\wh);
\draw (\wh,0,0) -- (\wh,0,\wh);
\draw (0,-\wL,\wh) -- (\wh,-\wL,\wh);
\draw (0,0,\wh) -- (\wh,0,\wh);
\draw (\wh,-\wL,\wh) -- (\wh,0,\wh);
\fill[gray!50] (\wh,0,-\wL) -- (\wh,\wh,-\wL) -- (\wh,\wh,0) -- (\wh,0,0) -- cycle;
\fill[gray!38] (0,\wh,-\wL) -- (\wh,\wh,-\wL) -- (\wh,\wh,0) -- (0,\wh,0) -- cycle;
\fill[gray!25] (0,0,0) -- (\wh,0,0) -- (\wh,\wh,0) -- (0,\wh,0) -- cycle;
\draw (\wh,0,-\wL) -- (\wh,\wh,-\wL);
\draw (\wh,0,0) -- (\wh,\wh,0);
\draw (\wh,\wh,-\wL) -- (\wh,\wh,0);
\draw (0,\wh,-\wL) -- (\wh,\wh,-\wL);
\draw (0,\wh,0) -- (\wh,\wh,0);
\draw (\wh,\wh,-\wL) -- (\wh,\wh,0);
\draw (0,\wh,0) -- (\wh,\wh,0);
\draw (\wh,0,0) -- (\wh,\wh,0);
% edges E_jl (hidden behind the expanded edges W_jl, and excluded from them)
\draw[very thick,dashed] (0,0,0) -- (-\wL-1.3,0,0) node[right] {$\mathcal E_{23}$};
\draw[very thick,dashed] (0,0,0) -- (0,-\wL-0.7,0) node[left] {$\mathcal E_{13}$};
\draw[very thick,dashed] (0,0,0) -- (0,0,-\wL-0.9) node[left] {$\mathcal E_{12}$};
\fill (0,0,0) circle (1.6pt) node[above left=1pt] {$a$};
\node[inner sep=0pt] at (0.42*\wh,\wh,-0.62*\wL) {$\tilde{\mathcal E}_{12}$};
\node at (\wh,-0.55*\wL,0.8*\wh) {$\tilde{\mathcal E}_{13}$};
\node at (-0.58*\wL,0.12*\wh,\wh) {$\tilde{\mathcal E}_{23}$};
% widths
\node[above left=-1pt] at (0.5*\wh,-\wL,\wh) {$\frac{1}{2r}$};
\node[below right] at (0.5*\wh,\wh,-\wL) {$\frac{1}{2r}$};
% coordinate directions
\begin{scope}[shift={(0,-2.7,-2.3)}]
\draw[->] (0,0,0) -- (0.7,0,0) node[below left=-1pt] {$x_1$};
\draw[->] (0,0,0) -- (0,0.7,0) node[right] {$x_2$};
\draw[->] (0,0,0) -- (0,0,0.7) node[above] {$x_3$};
\end{scope}
\end{tikzpicture}
\caption{The expanded edges $\tilde{\mathcal E}_{12}$, $\tilde{\mathcal E}_{13}$ and $\tilde{\mathcal E}_{23}$ for the orthant $\{x\le a\}$ are denoted by the solid gray regions in the case when $p=3$ and $R=I_3$. The edges $\mathcal E_{12}$, $\mathcal E_{13}$ and $\mathcal E_{23}$ are denoted by the dashed lines.}
\label{fig:wedges}
\end{figure}

\begin{proof}
To ease notation as in the previous proof, we fix $y$, and write $a=a(y)$, $\nu_j=\nu_j(y)$ and $\nu_{jl}=\nu_{jl}(y)$. In both~\eqref{eq:bt-facetsum} and~\eqref{eq:bt-edgesum}, the indices $j$ corresponding to large values of $|a_j|$ will be handled in part by the inequalities
\begin{equation}\label{eq:bt-density}
\nu_j\le\phi(a_j)\qquad\text{and}\qquad\nu_{jl}\le(2\pi(1-R_{jl}^2))^{-1/2}\phi(a_j)\qquad\text{for }j\ne l.
\end{equation}
The first of these follows from the formula $\nu_j=\phi(a_j)q_j$ in Lemma~\ref{lem:structure}, and the second follows from the formula for $\partial_j\partial_l \Phi_R(a)$ in Lemma~\ref{lem:bt-stein} by expanding the domain of integration from $\{z\leq a, z_j=a_j,z_l=a_l\}$ to  $\{ z_j=a_j,z_l=a_l\}$. To see this, note that the latter domain corresponds to marginalizing out all the coordinates except for $j$ and $l$, and then the integral is simply the bivariate normal density of $(Y_j,Y_l)$ at $(a_j,a_l)$. In turn, that density can be upper bounded by factoring it as  $\phi(a_j)$ times the conditional density of $Y_l|Y_j$ at $a_l$, which is at most $(2\pi(1-R_{jl}^2))^{-1/2}$.

Now we turn to the proof of~\eqref{eq:bt-facetsum}. Since $Y$ has unit variances, Nazarov's inequality (Lemma~\ref{fact:nazarov}) gives
$\Phi_R(a+\eps\1)-\Phi_R(a)\le\eps(\sqrt{2\log p}+2)$ for all $\eps>0$. Dividing by $\eps$ and letting $\eps\downarrow0$, we get
$\sum_{j=1}^p\nu_j\le\sqrt{2\log p}+2$. Next, we let $r=2\sqrt{\log(2p)}$ and we split the indices $j$ according to whether or not $\abs{a_j}\le r$ holds. Using~\eqref{eq:bt-density} and the fact that $t\phi(t)$ is decreasing on $[1,\infty)$, we obtain
\[
\sum_{j=1}^p\abs{a_j}\nu_j \ \le 
\ \sum_{j:|a_j|\leq r}^p r\nu_j+\sum_{j:\,\abs{a_j}>r}\abs{a_j}\phi(a_j)\le r\big(\sqrt{2\log p}+2\big)+p\,r\phi(r)\lesssim\log(2p),
\]
which proves~\eqref{eq:bt-facetsum}.

The argument for~\eqref{eq:bt-edgesum} is more involved. We split the pairs $j<l$ according to whether $\max(a_j,a_l)\le r$. If $\max(a_j,a_l)>r$, say
$a_j>r$, then~\eqref{eq:bt-density} gives $(1-R_{jl}^2)\nu_{jl}\le\phi(a_j)\le\phi(r)\le(2p)^{-2}$, and since there are fewer than $p^2$
pairs, these pairs contribute at most $\mathcal{O}(1)$ to the sum in~\eqref{eq:bt-edgesum}. Hence it remains to show that 
\begin{equation}\label{eq:bt-smallpairs}
\sum_{j<l} \1\{a_j\le r,\ a_l\le r\}(1-R_{jl}^2)\,\nu_{jl}\lesssim\log(2p).
\end{equation}
The idea is to thicken each edge $\mathcal E_{jl}=\{x\le a:\ x_j=a_j,\ x_l=a_l\}$ of the orthant $\{x\le a\}$ into an ``expanded edge''
$\tilde{\mathcal E}_{jl}\subset\R^p$, in such a way that the expanded edges are  disjoint and satisfy
\begin{equation}\label{eqn:nujlbound}
    (1-R_{jl}^2)\nu_{jl}\lesssim r^2\,\PP(Y\in\tilde{\mathcal E}_{jl}) 
\end{equation}
where $a_j,a_l\le r$ and $Y\sim \mathcal{N}(0,R)$. The disjointness of the expanded edges implies
$\sum_{j<l}\PP(Y\in\tilde{\mathcal E}_{jl})\le1$, and hence leads to~\eqref{eq:bt-smallpairs} once~\eqref{eqn:nujlbound} is established.
To define the expanded edge for a pair of indices $(j,l)$ with $j<l$, we let
\[
\tilde{\mathcal E}_{jl}=\mathcal E_{jl}+\big\{sRe_j+tRe_l:\ s,t\in(0,1/(2r)]\big\},
\]
which is understood in the sense of a Minkowski sum. (An illustration of these expanded edges is provided in Figure~\ref{fig:wedges}.)

We address the disjointness of the expanded edges first. We will show that if there are two points $\tilde x\in \tilde{\mathcal{E}}_{jl}$ and $\tilde x'\in\tilde{\mathcal{E}}_{j'l'}$ satisfying $\tilde x=\tilde x'$, then $(j,l)=(j',l')$. As an initial step, note that the equation $\tilde x=\tilde x'$ can be equivalently written as  $x+sRe_j+tRe_l=x'+s'Re_{j'}+t'Re_{l'}$, for some $x\in\mathcal E_{jl}$, some
$x'\in\mathcal E_{j'l'}$, and some $s,t,s',t'>0$. Rearranging the last equation so that $x-x'$ is on the left, and then applying $R^{-1}$ to both sides, we obtain
$
R^{-1}(x-x')=s'e_{j'}+t'e_{l'}-se_j-te_l.
$
Next, taking the inner product with $x-x'$ directly and noting that $x_j=a_j$, $x_l=a_l$ (and similarly for $x'$), we have
\begin{equation*}
    \begin{split}
        \ip{x-x'}{R^{-1}(x-x')}& \ = \ s'(x_{j'}-a_{j'})+t'(x_{l'}-a_{l'})+s(x_j'-a_j)+t(x_l'-a_l)
     \ \le \ 0,
    \end{split}
\end{equation*}
where the inequality holds due to $x\le a$ and $x'\le a$. On the other hand, since $R$ is positive definite, this forces $x=x'$, so that
$se_j+te_l=s'e_{j'}+t'e_{l'}$. In turn, since $j<l$ and $j'<l'$, this forces  $(j,l)=(j',l')$.  Thus the sets $\tilde{\mathcal E}_{jl}$ are 
disjoint.

Completing the proof is now reduced to showing~\eqref{eqn:nujlbound}. This will be done by directly calculating the probability $\PP(Y\in\tilde{\mathcal{E}}_{jl})$  and showing that it is at least $(1-R_{jl}^2)\nu_{jl}/r^2$, up to an absolute constant. Throughout this calculation, we consider a fixed pair of indices $j<l$ with $a_j,a_l\le r$.
An argument similar to the one in the previous paragraph shows that the map
$T_{jl}(x,s,t)=x+sRe_j+tRe_l$ is a bijection between $\mathcal E_{jl}\times(0,1/(2r)]^2$ and $\tilde{\mathcal E}_{jl}$.
Since every $x\in\mathcal E_{jl}$ has $x_j=a_j$ and $x_l=a_l$, we may parametrize $\mathcal E_{jl}$ by the remaining
coordinates $x_k$ with $k\ne j,l$. So, the change-of-variables formula gives
\begin{equation}\label{eq:bt-cov}
\PP(Y\in\tilde{\mathcal E}_{jl})=\int_{\tilde{\mathcal E}_{jl}}\phi_R(z)\,dz=\abs{\det J_{jl}}\int_{\mathcal E_{jl}}\int_0^{1/(2r)}\!\!\int_0^{1/(2r)}\phi_R\big(T_{jl}(x,s,t)\big)\,ds\,dt\,d\mathcal H^{p-2}(x),
\end{equation}
where $J_{jl}$ is the Jacobian matrix of $T_{jl}$, which does not depend on $(x,s,t)$ because $T_{jl}$ is affine. Since
$\nu_{jl}=\int_{\mathcal E_{jl}}\phi_R\,d\mathcal H^{p-2}$ by Lemma~\ref{lem:bt-stein}, it is enough to show $\phi_R(T_{jl}(x,s,t))\ge e^{-3/2}\phi_R(x)$ and  $\abs{\det J_{jl}}=1-R_{jl}^2$, because these two facts and~\eqref{eq:bt-cov} combine to give~\eqref{eqn:nujlbound}.

For the determinant formula, note that the $k$th column of $J_{jl}$ is equal to $e_k$ except when $k=j,l$, and in these two cases the columns are  $Re_j$
and $Re_l$. Based on this, the formula $|\det(J_{jl})|=1-R_{jl}^2$ can be obtained by reordering the columns of $J_{jl}$ so that it has an upper-left block equal to $I_{p-2}$ and it has a lower-right block with determinant $1-R_{jl}^2$. Finally, for the lower bound on the density, if we let $x\in\mathcal E_{jl}$ and $v=se_j+te_l$ with $s,t\in(0,1/(2r)]$, then
\[
%\begin{split}
\phi_R\big(T_{jl}(x,s,t)\big)
=\phi_R(x+Rv)
=\phi_R(x)\exp\big(-\langle x,v\rangle -\tfrac12 \langle v, Rv\rangle \big)
\ge e^{-3/2}\phi_R(x),
%\end{split}
\]
where the last bound is due to 
 $\ip{x}{v}=sa_j+ta_l\le1$, and
$\ip{v}{Rv}=s^2+2R_{jl}st+t^2\le1$, because $s,t\leq \frac{1}{2r}$.

\end{proof}

\section{Proof of Theorem~\ref{thm:boot}}\label{sec:bootproof}

Conditionally on $X_1,\dots,X_n$, the bootstrap sum $S_n^*$ is centered Gaussian with covariance matrix $\hat\Sigma$. So, once
$d(\hat\Sigma,\Sigma)$ is controlled, Theorem~\ref{thm:bt} can be applied conditionally on $X_1,\dots,X_n$ with $\hat\Sigma$ playing the role of $\tilde\Sigma$. This is handled by the next lemma.

\begin{lemma}\label{lem:bt-conc}
There is an absolute constant $C$ such that the following holds. Let $p,n\ge2$ and let $X_1,\dots,X_n$ be i.i.d.\ centered
random vectors in $\R^p$ with covariance matrix $\Sigma=\E(X_1X_1^\top)$, such that $\norm{\ip{X_1}{v}}_{\psi_2}\le B\sqrt{\ip{v}{\Sigma v}}$ for
all $v\in\mathcal V$, where $\mathcal V$ is the set in~\eqref{eq:bt-disc}. Also, let $\hat\Sigma=\frac1n\sum_{i=1}^nX_iX_i^\top$. Then, for every
$\alpha\in(0,1)$ such that $n\ge C\log(p/\alpha)$, the event
\[
d(\hat\Sigma,\Sigma)\le CB^2\sqrt{\frac{\log(p/\alpha)}n}
\]
holds with probability at least $1-\alpha$.
\end{lemma}

\begin{proof}
If $\ip{v}{\Sigma v}=0$ holds for some $v\in\mathcal{V}$, then $\ip{X_i}v=0$ almost surely for all $i\leq n$, and hence  $\ip{v}{(\hat\Sigma-\Sigma)v}=0$ almost surely. Otherwise, we have $\langle v,\Sigma v\rangle>0$, and in this case the variables $\zeta_i(v)=\ip{X_i}{v}/\sqrt{\ip{v}{\Sigma v}}$ satisfy $\E \zeta_i^2(v)=1$, $\norm{\zeta_i(v)}_{\psi_2}\le B$, and 
$$\frac{\ip{v}{(\hat\Sigma-\Sigma)v}}{\ip{v}{\Sigma v}}=\frac1n\sum_{i=1}^n (\zeta_i(v)^2-1). $$
By Bernstein's inequality~\citep[Corollary~2.8.3]{Ver18},  there is an absolute constant $c>0$ such that $\PP(\abs{\frac1n\sum_i(\zeta_i^2(v)-1)}>t)\le2\exp(-c n\min(t^2/B^4,t/B^2))$ for $t>0$. Taking a union bound over the $p^2$ vectors in $\mathcal V$ completes the proof.
\end{proof}

\begin{proof}[Proof of Theorem~\ref{thm:boot}]

Let $\Sigma=\E(X_1X_1^{\top})$ and $G\sim\mathcal{N}(0,\Sigma)$. By Theorem~\ref{thm:main} and the triangle inequality, it suffices to show that the following bound holds with probability at least $1-\alpha$,
\begin{equation}\label{eq:boot1}
\sup_{A\in\mathscr R}\big|\PP^*(S_n^*\in A)-\PP(G\in A)\big|\lesssim \,B^2\sqrt{\frac{\log(p/\alpha)}n}\,\log(p)\log(n).
\end{equation}
Note that we must have $B\geq 1$ because the $\psi_2$-norm upper bounds the $L^2$-norm. Letting $C$ be as in the statement of Lemma~\ref{lem:bt-conc}, we may assume without loss of generality that $n\geq C\log(p/\alpha)$, because otherwise~\eqref{eq:boot1} is trivially true. Thus, we may apply Lemma~\ref{lem:bt-conc}, and if we let $\delta=CB^2\sqrt{\log(p/\alpha)/n}$, then the event $d(\hat\Sigma,\Sigma)\le\delta$ holds with
probability at least $1-\alpha$. On this event, Theorem~\ref{thm:bt}, applied conditionally on $X_1,\dots,X_n$ with $\tilde\Sigma=\hat\Sigma$, gives
\[
\sup_{A\in\mathscr R}\big|\PP^*(S_n^*\in A)-\PP(G\in A)\big|\lesssim\delta\log(p)\log(1+1/\delta),
\]
and since our choice of $\delta$ implies $\log(1+1/\delta)\lesssim\log n$, the proof is complete.
\end{proof}

\section{Statement on AI use}\label{sec:AI}
Artificial intelligence was used in this work to develop proofs, check the arguments for correctness, and edit language. The authors confirm that the paper conforms with the Institute of Mathematical Statistics policy on use of AI. Essentially all of the text reflects the writing of the authors. The authors retain full responsibility for all aspects of the paper.

\bibliographystyle{spr-ims-nameyear}
\bibliography{refs}

\end{document}